\documentclass{article}
\usepackage{geometry}
\usepackage[T1]{fontenc}
\usepackage[utf8]{inputenc}
\usepackage{amsmath, amssymb}
\usepackage{textcomp}
\usepackage{eurosym}
\usepackage{graphicx}
\usepackage{subcaption}
\usepackage{booktabs}
\usepackage{multirow}
\usepackage{float}
\usepackage[table]{xcolor}
\usepackage{algorithm}
\floatname{algorithm}{Algorithm}
\usepackage{algpseudocode}
\usepackage[numbers]{natbib}
\usepackage[hidelinks]{hyperref}
\usepackage{url}

\newenvironment{keywords}{\vspace{0.4in}\begin{quote}\small \em
 {\bf Keywords\/}:}{\end{quote}}

\title{Evaluation of Electricity Market Clearing Mechanisms via Reinforcement Learning:
       Prices, Remuneration and Competitive Dynamics}

\author{Andrea Altamura \\
        Department of Electrical and Information Engineering, Polytechnic University of Bari \\
        {\tt a.altamura12@studenti.poliba.it}
        \and
        Fabrizio Lacalandra \\
        Regulatory Authority for Energy, Networks and Environment (ARERA)\footnote{This technical report has an eminently scientific character. The ideas and the conclusions expressed here \emph{do not} necessarily represent those of ARERA, where the second author works.} \\
        {\tt flacalandra@arera.it}
        \and
        Antonio Frangioni \\
        Department of Computer Science, University of Pisa, \\
        {\tt frangio@di.unipi.it}
        }

\date{\today\\}

\begin{document}
\maketitle

\begin{abstract}
    The Pay-as-Clear (PaC) mechanism currently used in the European electricity market can generate significant submarginal profits for renewable sources when the clearing price is determined by the marginal offers of gas-fired generation units and the cost of natural gas exceeds certain levels. This exposes consumers to high price volatility related to the cost of natural gas. This report analyzes the recently proposed Segmented Pay-as-Clear (SPaC) mechanism as a market alternative, evaluating its system cost-effectiveness through simulations based on Reinforcement Learning (Q-Learning) to model the strategic behavior of operators. Three market models are compared, the two classic Pay-as-Clear (PaC) and Pay-as-Bid (PaB) along with SPaC, under two scenarios: a simplified one based on the 2030 NECP objectives and one built on the portfolios of ten operators obtained from the GME's 2024 public offers. The results show that the SPaC market clearing mechanism reduces intramarginal profits and price volatility compared to PaC, while maintaining fair participation incentives for all operators, and is more robust than PaB to the exercise of market power in oligopolistic contexts. The developed framework can serve as a support tool for regulators and policymakers in the evaluation of proposals for market design reforms.
\end{abstract}


\tableofcontents
\listoffigures
\listoftables

\section{Introduction}
\label{sec:introduzione}
The European electricity market is based on the uniform marginal price mechanism (Pay-as-Clear, PaC), which has represented for decades a pillar of the internal energy market architecture. This model, in which all accepted production units are remunerated at the marginal unit's price, has ensured efficient resource allocation in the short term and investment signals consistent with long-term economic principles. However, the ongoing energy transition, and the recent energy crisis events, have highlighted some structural limitations of the approach regarding economic sustainability for consumers. It was expected that the growing penetration of renewable sources, characterized by marginal production costs close to zero, would progressively alter price formation dynamics. This has not happened significantly in Italy. In particular, even during periods of high renewable availability, the market price continues frequently to be determined by thermoelectric technologies, mainly fueled by natural gas, potentially generating significant inframarginal rents for renewable producers.

\smallskip
\noindent
This phenomenon emerged dramatically during the 2022 gas crisis, when gas prices---and consequently electricity prices---reached historically unprecedented levels, resulting in unsustainable costs for final consumers in the face of exceptional profits for producers with technologies other than natural gas. This is also worsened by a side effect of the Emission Trading System (ETS) mechanism: fossil fuel-fired thermoelectric plants, which often are the price-maker marginal units, correctly internalize in their sales offers the cost of $\text{CO}_2$ emission allowances they must pay for the ETS system, in addition to the gas price. Consequently, through the Pay-as-Clear mechanism, this cost is incorporated into the market marginal price and thus transferred into the remuneration of all accepted producers, including those from sources that bear no emission costs. This might be considered an implicit incentive to install new renewable capacity, but in fact for a long time other explicit incentive mechanisms have been in force. Today, therefore, a significant economic effect (or distortion?) arises\footnote{The price of the ETS allowance varied in the last year approximately between \EUR{} 70/Ton and \EUR{} 80/Ton which, considering 400 kg of CO2 emitted by a CCGT per MWh of electricity produced, results in a burden between \EUR{} 28 and \EUR{} 32 per MWh of electricity in hours when CCGTs are marginal.}: renewable technologies receive a price component associated with ETS allowances without bearing the corresponding burden, further increasing inframarginal rents and amplifying the overall cost for the final consumer.

\smallskip
\noindent
The direct link between electricity prices and natural gas costs has thus exposed consumers to a significant price increase, raising questions of distributive equity and economic sustainability of the energy transition, as well as national competitiveness of industries that significantly use electricity. This suggests the opportunity to rethink the spot electricity market design to decouple, at least partially, the remuneration of renewable sources from that of fossil sources, while preserving short-term allocative efficiency signals and long-term investment incentives.

\subsection{Objectives and contributions of the report}

This report aims to analyze and compare different electricity market clearing mechanisms, with a specific focus on a recent innovative proposal introduced in the literature: \emph{Segmented Pay-as-Clear (SPaC)} \cite{frangioniBilevelProgrammingApproach2024}. This model represents an original solution to the problem of decoupling the remuneration of groups of plants, e.g., those producing electricity using fossil fuels (mainly gas) from those that do not. Although more complex scenarios are also possible, in this work we will consider market segmentation into two coordinated sub-markets: one for Negligible Marginal Cost Sources (NMCS, mainly renewables), and one for Non-Negligible Marginal Cost Sources (NNMCS, mainly thermoelectric). The distinctive element of the SPaC mechanism lies in the fact that the demand split between the two sub-markets is determined endogenously by the model through a bespoke bilevel optimization problem, with the objective of minimizing the total system cost. This architecture allows preserving the benefits of marginal pricing, i.e., allocation efficiency and investment signals, while simultaneously reducing consumer exposure to gas price volatility.

\smallskip
\noindent
Our analysis is not based on theoretical simulations under perfect competition conditions, but explicitly explores the strategic behavior of market operators. To this end, we adopted an approach based on \emph{Reinforcement Learning (Q-Learning)}, which allows simulating adaptive agents capable of autonomously learning profitable bidding strategies through repeated interaction with the market environment. This agent-based approach allows analyzing the robustness of different market mechanisms (PaC, Pay-as-Bid and SPaC) with respect to the exercise of market power, without requiring unrealistic assumptions about perfect rationality or complete knowledge of information by operators. The main contribution of this work therefore consists in the---to the best of our knowledge, novel---application of the Reinforcement Learning approach to evaluation of the SPaC market clearing mechanism, assessing the ability of operators to learn optimal strategies within a segmented market and analyzing the implications for social welfare, profit distribution, and price formation. Remarkably, price remains unique on the demand side for each relevant period: this is relevant in that it would allow individual member states to adopt the SPaC approach within the EU-wide single market cleared by the so-called Euphemia model, without it necessarily being adopted by all and thereby require a change to the reference regulations. Further research would clearly be necessary to assess the global effect of such a change.

\subsection{Methodology and analysis scenarios}

The analysis is conducted through an agent-based simulation framework developed in the Julia programming language \cite{bezansonJuliaFreshApproach2017}, which integrates:
\begin{itemize}
 \item a market clearing module capable of solving the optimization problems associated with the different mechanisms (PaC, PaB and SPaC), including the management of the bilevel problem for the SPaC regime;
 \item a multi-agent environment in which each market operator is modeled as an autonomous RL agent learning bidding strategies through the Q-Learning algorithm;
 \item a data collection and analysis system that allows evaluating the emergent equilibrium outcomes in terms of prices, profits, social welfare, and market concentration indicators.
\end{itemize}
Simulations are conducted on two distinct scenarios:
\begin{enumerate}
 \item \emph{PNIEC 2030 Scenario}: a simplified scenario based on the objectives of the Italian National Integrated Energy and Climate Plan, with a technological mix representative of the target configuration for 2030 and a limited number of aggregated operators. This scenario allows isolating the effects of different market architectures in a controlled context.
 \item \emph{10-operator Scenario}: a scenario based on actual data from public offers of the Italian electricity market (MGP), which includes a pseudo-real structure of the generation fleet from which portfolios for 10 operators are derived. This scenario allows evaluating the feasibility and effectiveness of the different mechanisms under realistic market conditions.
\end{enumerate}

\subsection{Document structure}

The rest of the report is organized as follows:
\begin{itemize}
 \item Section \ref{sec:regimi-mercato} describes in detail the three analyzed market models: the traditional ones, PaC and PaB, and the new SPaC, formalizing their mathematical clearing models;
 \item Section \ref{sec:reinforcement-learning} introduces the Reinforcement Learning framework used to model operator strategies, describing the Q-Learning algorithm, the definition of the state and action space, and the implementation choices related to learning parameters;
 \item Section \ref{sec:simulazioni} presents the results of the numerical simulations conducted on the two scenarios, analyzing the convergence of the learning process, emergent strategies, price dynamics, and the distribution of surplus between producers and consumers.
\end{itemize}
\section{Market Models}
\label{sec:regimi-mercato}

In this section, the three main price formation schemes considered are briefly presented: PaC, PaB and SPaC. The objective is to define clearing rules, strategic incentives, and main economic implications for consumers and producers.

\subsection{Pay-as-Clear (PaC)}

In the uniform price mechanism of \emph{Pay-as-Clear} all accepted offers are valued at the same price, equal to the offer of the unit resulting as marginal (System Marginal Price, SMP). The clearing coincides with the maximization of the so-called \emph{Social Welfare} under other possible constraints, such as network constraints or those derived from the peculiar characteristics of certain markets (e.g., block offers): the equilibrium price equalizes the marginal benefit of demand with the marginal cost of supply~\citep{kirschenFundamentalsPowerSystem2004,schweppeSpotPricingElectricity1988}. This setup guarantees allocation efficiency in the short term and consistent signals for long-term investments, but can generate inframarginal rents for low marginal cost technologies (e.g., RES), more pronounced when the price-setting unit is thermoelectric. Transparency is high, while the outcome remains sensitive to the degree of actual competition and demand/supply conditions \citep{boscoStrategicBiddingVertically2012}.

\subsection{Pay-as-Bid (PaB)}

In \emph{Pay-as-Bid}, remuneration is discriminatory: each accepted unit is paid at its own offer price, while still following offer acceptance in economic merit order. This results in an incentive to include a markup over marginal cost to cover fixed costs, with the risk of not being accepted; the optimal strategy involves estimating the expected price of the acceptance threshold. In terms of efficiency and transparency, outcomes depend on the degree of competition and operational rules, with possible deviations from efficient allocation~\citep{ferrariCompetitionElectricityMarkets2005}. Several studies based on agent-based simulations show that PaC tends to induce offers closer to the actual marginal costs of units and, sometimes, higher overall profits compared to PaB; however, such results depend on the rationing policy and the competitive context \citep{aliabadiAgentbasedSimulationPower2017}.

\subsection{Segmented Pay-as-Clear (SPaC)}\label{ssec:SPaC}

\begin{figure}[htb]
	\centering
	\includegraphics[width=0.8\textwidth]{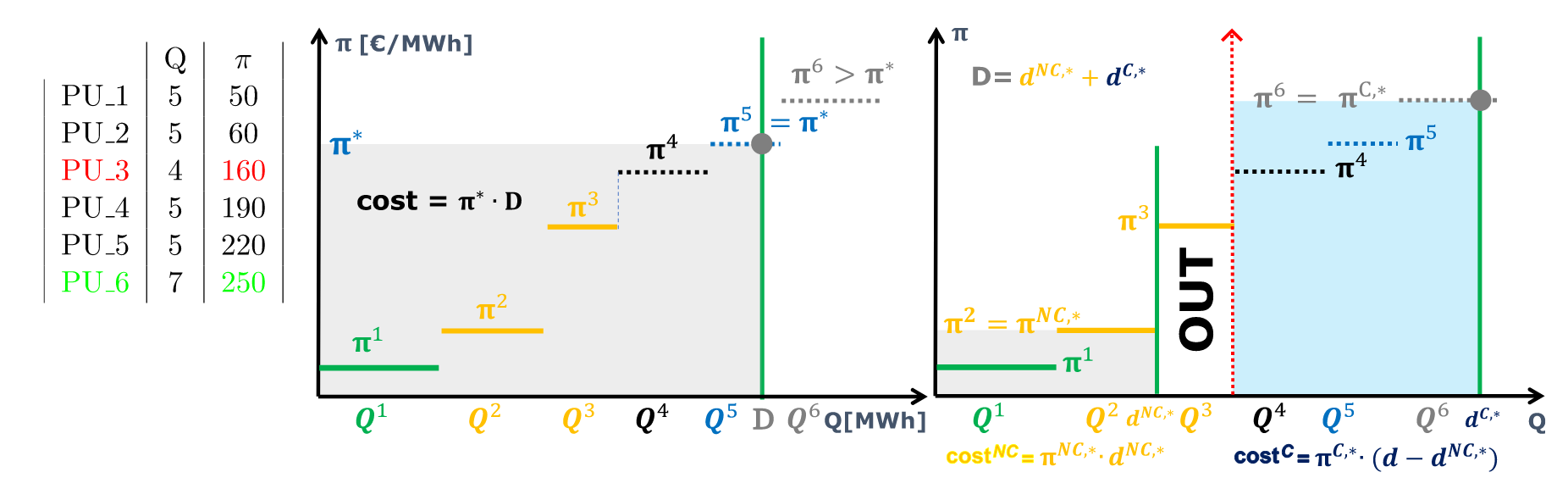}
	\caption{Illustrative example of the SPaC model \citep{frangioniBilevelProgrammingApproach2024}.}
	\label{fig:esempio-spac}
\end{figure}

The \emph{Segmented Pay-as-Clear} (SPaC) model, recently proposed in \citep{frangioniBilevelProgrammingApproach2024}, represents an innovative alternative. The fundamental intuition consists in preserving the benefits of the PaC mechanism (efficient short- and long-term price signals), introducing an optimal market segmentation that decouples technologies with negligible marginal costs (NMCS) from those with non-negligible marginal costs (NNMCS). The SPaC mechanism operates by creating two separate virtual markets solved simultaneously, see Figure \ref{fig:esempio-spac} for a graphical comparison example. The innovative element lies in the fact that the split of demand between these two markets is not fixed a priori, but is a decision variable of the model. In the case of rigid demand, this split is determined with the objective of minimizing the total procurement cost for the consumer, ensuring that the sum of accepted quantities in the two segments satisfies total demand.

\begin{figure}[htb]
	\centering
	\includegraphics[width=0.45\textwidth]{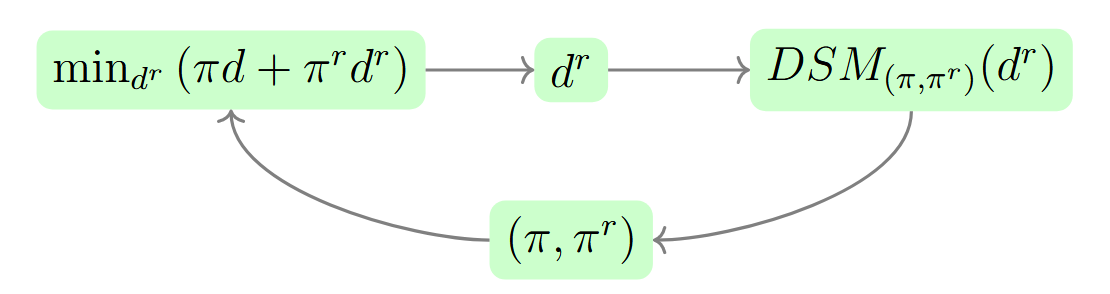}
	\caption{Bilevel structure of the SPaC model: the leader minimizes total cost by choosing $d_r$, the followers determine the marginal prices $\pi_r, \pi_g$ \citep{frangioniBilevelProgrammingApproach2024}.}
	\label{fig:struttura-bilevel}
\end{figure}

\noindent
From a mathematical point of view, assuming rigid demand $d > 0$, the problem is configured as a \emph{bilevel} optimization. Given a set $S$ of offers $(sp_j, sq_j)$ partitioned into $S = S_r \cup S_g$ (where $S_r$ is the set of NMCS and $S_g$ that of NNMCS), the model seeks the optimal demand split $d_r + d_g = d$ which, together with the prices, minimizes the overall system cost, i.e., the objective function \eqref{FO:Imp}. The problem is formally expressed in implicit form as follows:
\begin{align}
	\textstyle \min_{d_r, d_g} \quad & \pi_r d_r + \pi_g d_g \label{FO:Imp} \\
	\text{s.t.} \quad & d_r + d_g = d  \\
	& d_r \geq 0, \quad d_g \geq 0 \\
	& \text{with } \pi_r, \pi_g \text{ equilibrium prices of the submarkets:} \nonumber \\
	& \text{Clearing } S_r(d_r) \rightarrow \pi_r \\
	& \text{Clearing } S_g(d_g) \rightarrow \pi_g
\end{align}
where prices $\pi_r$ and $\pi_g$ emerge as dual variables from the respective social cost minimization problems (standard PaC clearing) constrained to the variable demand quotas $d_r$ and $d_g$. The \emph{bilevel} nature derives from the decision hierarchy: the \emph{leader} chooses the split $d_r$ to minimize total cost, while the \emph{followers} (the two segmented markets) react by determining the equilibrium prices $\pi_r$ and $\pi_g$ (Figure \ref{fig:struttura-bilevel}).

\noindent
An alternative, mathematically equivalent formulation considers a single market with an explicit constraint on the maximum salable quantity from NMCS technologies. For each non-negative value of the constant $d_r$, we have the linear problem \eqref{eq:POF-S}--\eqref{eq:mrkteq-S}:
\begin{align}
	\textstyle \min_{s_j} \quad &  \textstyle \sum_{j \in S} sp_j s_j \label{eq:POF-S}\\
	\text{s.t.} \quad & 0 \leq s_j \leq sq_j \quad \forall j \in S \label{eq:Sbnd-S}\\
	& \textstyle \sum_{j \in S_r} s_j \leq d_r \quad (\text{dual variable } \pi_r \leq 0) 
      \label{eq:UPbnd-S} \\
	& \textstyle \sum_{j \in S} s_j = d \quad\;\;\; (\text{dual variable } \pi)
      \label{eq:mrkteq-S}
\end{align}
The dual of the linear problem \eqref{eq:POF-S}--\eqref{eq:mrkteq-S}, say $DSM(d_r)$, generates two distinct remunerations: $\pi$ for NNMCS and a lower or equal $\pi + \pi_r$ for NMCS (since $\pi_r \leq 0$). Ultimately, the bilevel problem solved by SPaC is written in compact form \eqref{eq:bi-pi-argmin-S} where $d_r$ becomes a variable \cite{frangioniBilevelProgrammingApproach2024}.
\begin{equation}
 \textstyle
 \min \big\{ \, \pi d + \pi_r d_r  \;:\;
    \left( \pi,\pi_r\right) \in \mbox{argmax\;} DSM(d_r) \;,\; d_r \geq 0 \, \big\}
    \label{eq:bi-pi-argmin-S} 
\end{equation} 
The problem can be solved as an MPCC (\emph{Mathematical Program with Equilibrium Constraints}) by explicitly stating the KKT conditions of the linear problem. The complementarity conditions ensure that:
\begin{align}
	\eta_j (sq_j - s_j) &= 0 \quad \forall j \in S \\
	(sp_j - \eta_j - \pi - \pi_r) s_j &= 0 \quad \forall j \in S_r\\
	(sp_j - \eta_j - \pi) s_j &= 0 \quad \forall j \in S \setminus S_r
\end{align}
These guarantee that offers below the marginal remuneration of their own segment are accepted and those above are rejected, maintaining the logic of uniform remuneration within each group. It is proven in \citep{frangioniBilevelProgrammingApproach2024} that, in the optimal solution, if the remunerations of the two segments differ, the marginal offer setting the remuneration in the NMCS market is always fully accepted. This implies that optimal segmentation exploits discontinuities in the aggregated supply curve.

\smallskip
\noindent
The significant difference w.r.t.~PaB is that SPaC maintains the \emph{intra-segment} uniformity of remuneration. While PaB remunerates each unit at its offered price, incentivizing strategic \emph{markup}, SPaC preserves competition within homogeneous technological groups and between them. RES, for example, compete among themselves in the NMCS segment becoming \emph{price-maker} in their own market, but without directly suffering gas price volatility. This also reduces the risk of market power exercise: out-of-market offers from renewables risk not being accepted, with consequent shift of a quota of demand to the other NNMCS market, thus realizing a dynamic competition where strategic offering is more difficult. Furthermore, the fragmentation of RES operators makes collusion difficult, which remains subject to the Regulation REMIT surveillance anyway \citep{acerACERGuidanceREMIT2024}.

\section{Bidding Strategy and Q-Learning}
\label{sec:reinforcement-learning}

The analysis of bidding strategies in electricity markets PaC, PaB and SPaC poses significant challenges from both theoretical and computational points of view. In the SPaC case, in particular, the bilevel nature of the clearing problem, with endogenous demand split between segmented markets, makes it difficult to directly apply traditional equilibrium concepts (Nash, Cournot, Supply Function Equilibrium) proposed for conventional markets. Furthermore, the assumption of complete information on costs and competitors' strategies, already unrealistic in PaC and PaB regimes, becomes even more problematic in an innovative market context like SPaC, where the mechanism structure itself is new to operators.

\smallskip
\noindent
To address this complexity, an approach based on Reinforcement Learning (RL) is adopted, where operators are modeled as agents that dynamically learn profitable bidding strategies through repeated interaction with the market environment. In this framework, market clearing (PaC, PaB or SPaC) constitutes the "black-box" that, given the offer profile of operators, returns remunerations, prices, quantities, and profits; agents update their strategies based on the observed reward. The use of RL techniques in agent-based environments to evaluate new market rules has been explored in projects like ASSUME \cite{harderASSUMEAgentbasedSimulation2025}. However, to the best of our knowledge there is no previous research that analyzes the strategic behavior of operators within a segmented SPaC electricity market, making this analysis an original contribution in this area.

\subsection{Basic Elements of Reinforcement Learning}
\label{subsec:rl-basics}

Reinforcement Learning is a machine learning paradigm where an \emph{agent} learns to make optimal decisions through interaction with a dynamic \emph{environment} \cite{suttonReinforcementLearningIntroduction2018}.
Unlike supervised learning, where the agent learns from labeled examples, in RL the agent autonomously discovers the best actions through a trial-and-error process, receiving feedback in the form of \emph{reward} or penalties. The fundamental elements of an RL system are:
\begin{itemize}
 \item \emph{Agent}: the entity that makes decisions and learns from received feedback. In the context of electricity markets, each operator constitutes an autonomous agent with a portfolio of production units (PU).
 \item \emph{Environment}: the external system with which the agent interacts. In the specific case, the environment is represented by the market mechanism (PaC, PaB or SPaC) that, given the offer vector, determines clearing remunerations/prices, accepted quantities, and operator profits.
 \item \emph{State} ($s \in \mathcal{S}$): a representation of the current conditions of the environment that the agent can observe; in the implemented simulations, the state corresponds to the electricity demand level $D$ at a given moment.
 \item \emph{Action} ($a \in \mathcal{A}$): the decisions the agent can undertake in a given state. For market operators, the action consists in choosing the markups to apply to their offers for each technology in the portfolio.
 \item \emph{Reward} ($r$): the numerical feedback the agent receives after executing an action. In the implemented model, the reward corresponds to the profit realized by the operator in the market.
 \item \emph{Policy} ($\pi$): the strategy the agent uses to choose actions based on the observed state. The objective of RL is to find the optimal policy $\pi^*$ that maximizes expected long-term reward.
\end{itemize}
The RL process unfolds through sequential episodes: at each time step $t$, the agent observes state $s_t$, selects an action $a_t$ according to its policy, receives a reward $r_{t+1}$ and observes the new state $s_{t+1}$. The agent uses this experience to update its policy with the aim of maximizing cumulative future reward.

\subsection{Q-Learning}
\label{subsec:q-learning}

Q-Learning is a model-free RL algorithm that allows the agent to learn the value of state-action pairs without requiring an explicit model of the environment \cite{watkinsQLearning1992}. The algorithm is based on the value function $Q(s, a)$, which represents the expected reward the agent would obtain by executing action $a$ in state $s$ and subsequently following the optimal policy. The optimal Q-function $Q^*(s, a)$ satisfies the Bellman equation:
\begin{equation}
    Q^*(s, a) = \mathbb{E}\left[r + \gamma \max_{a'} Q^*(s', a') \mid s, a\right]
\end{equation}
where $r$ is the immediate reward, $\gamma \in [0, 1]$ is the discount factor balancing immediate and future rewards, and $s'$ is the next state.

\smallskip
\noindent
The Q-Learning algorithm implements an iterative update process based on the rule:
\begin{equation}
    Q(s, a) \leftarrow Q(s, a) + \alpha \left[r + \gamma \max_{a'} Q(s', a') - Q(s, a)\right]
\end{equation}
where $\alpha \in (0, 1]$ is the learning rate that determines the weight given to new information compared to already acquired information. In the model implemented for the considered electricity markets, since each market episode is independent and has no linked future states ($\gamma = 0$), the update rule simplifies to:
\begin{equation}
    Q(s, a) \leftarrow Q(s, a) + \alpha \left[r - Q(s, a)\right]
\end{equation}
where an adaptive learning rate is used:
\begin{equation}
    \alpha_t(s, a) = \frac{1}{\text{counts}(s, a)}
\end{equation}
with $\text{counts}(s, a)$ representing the number of times the state-action pair $(s, a)$ has been visited. This choice guarantees theoretical convergence of the algorithm and ensures that the influence of new observations progressively decreases as the agent gains experience.

\smallskip
\noindent
A crucial aspect of Q-Learning is balancing \emph{exploration} (exploring new actions) and \emph{exploitation} (exploiting acquired knowledge). This is managed through an $\varepsilon$-greedy strategy: with probability $\varepsilon$ the agent chooses a random action (exploration), while with probability $1-\varepsilon$ it selects the action with the highest Q-value (exploitation). In the implemented model, $\varepsilon$ decays exponentially over episodes ($t$) according to:
\begin{equation}
    \varepsilon(t) = \varepsilon_{\text{max}} \cdot e^{-\lambda t}
\end{equation}
where $\lambda$ is the decay rate calculated as:
\begin{equation}
    \lambda = -\frac{\ln(\varepsilon_{\text{min}} / \varepsilon_{\text{max}})}{T}
\end{equation}
with $T$ total number of training episodes. This approach guarantees intensive exploration in the initial learning phases, followed by progressive convergence towards exploitation of the learned policy.

\begin{algorithm}[htb]
    \caption{Q-Learning Operation for Each Operator}
    \label{alg:qlearning}
    \begin{algorithmic}[1]
        \State Initialize $Q(s,a) = 0$ for all states $s$ and actions $a$
        \State Initialize $\text{counts}(s,a) = 0$ for all states $s$ and actions $a$
        \For{episode = 1, 2, ..., $T$}
            \State Observe current state $s$ (demand level)
            \State Compute $\varepsilon = \varepsilon_{\text{max}} \cdot e^{-\lambda \cdot \text{episode}}$
            \For{each operator $i$}
                \If{$\text{rand}() < \varepsilon$}
                    \State Choose random action $a_i$ (exploration)
                \Else
                    \State Choose $a_i = \arg\max_a Q_i(s, a)$ (exploitation)
                \EndIf
            \EndFor
            \State Simulate market with actions $\{a_i\}$ and observe rewards $\{r_i\}$
            \For{each operator $i$}
                \State $\text{counts}_i(s, a_i) \leftarrow \text{counts}_i(s, a_i) + 1$
                \State $\alpha = 1 / \text{counts}_i(s, a_i)$
                \State $Q_i(s, a_i) \leftarrow Q_i(s, a_i) + \alpha [r_i - Q_i(s, a_i)]$
            \EndFor
        \EndFor
    \end{algorithmic}
\end{algorithm}

\subsection{Implementation for Electricity Markets}
\label{subsec:rl-implementation}

Implementing Q-Learning for electricity markets presents some specificities that differ from the classical problem. The strategic bidding problem can be modeled as a variant of the \emph{Multi-Armed Bandit Problem} \cite{suttonReinforcementLearningIntroduction2018}, where each "arm" of the slot machine corresponds to a different markup strategy.

\paragraph{State Space Definition}

The state space $\mathcal{S}$ is defined by electricity demand levels discretized in an interval $[D_{\text{min}}, D_{\text{max}}]$. In the implemented model:
\begin{equation}
    \mathcal{S} = \{D_{\text{min}}, D_{\text{min}} + \Delta D, ..., D_{\text{max}}\}
\end{equation}
where $D_{\text{min}} = 0.25 \cdot C_{\text{tot}}$ and $D_{\text{max}} = 0.80 \cdot C_{\text{tot}}$, with $C_{\text{tot}}$ total market capacity. Discretization uses 100 equidistant points to balance precision and computational complexity.

\paragraph{Action Space Definition}

The action space $\mathcal{A}_i$ for operator $i$ is defined by the combinations of markups applicable to the different technologies in the portfolio. If operator $i$ owns $n_i$ different technologies, the action is a tuple:
\begin{equation}
    a_i = (m_{i,1}, m_{i,2}, ..., m_{i,n_i}) \in \mathcal{M}^{n_i}
\end{equation}
where $\mathcal{M} = \{0, 5, 10, 20\}$ \% for PaC and SPaC markets, while $\mathcal{M} = \{0, 50, 100, 200\}$ \% for PaB market. The choice of higher markups for PaB serves to reflect operator strategic behavior in that regime.

\paragraph{Reward Function}

The reward for operator $i$ in episode $t$ corresponds directly to realized profit:
\begin{equation}
 \textstyle
 r_{i,t} = \sum_{j \in \mathcal{U}_i} (p_j - c_j) \cdot q_{j,t}
\end{equation}
where $\mathcal{U}_i$ is the set of production units of operator $i$, $p_j$ is the clearing price received by unit $j$, $c_j$ is its marginal cost, and $q_{j,t}$ is the dispatched quantity.

\paragraph{Training Process}

Training is performed separately for each market regime (PaC, PaB, SPaC), keeping the electric system model structure unchanged and modifying only the price formation rules; parallelization is used to accelerate the computational process. For each state $s \in \mathcal{S}$, agents perform $T = 2000$ learning episodes. The overall process therefore requires $|\mathcal{S}| \times T = 100 \times 2000 = 200{,}000$ simulations per regime. Convergence is monitored through analysis of profit evolution during training. Stabilization of Q-values and reduction of reward variability indicate that agents have learned consistent strategies.

\paragraph{Optimal Policy Extraction}

At the end of training, the optimal policy for operator $i$ in state $s$ is extracted as:
\begin{equation}
    \pi_i^*(s) = \arg\max \{ \, Q_i(s, a) \;:\; a \in \mathcal{A}_i \, \}
\end{equation}
For intermediate states not present in the discretization, interpolation based on Euclidean distance is used: the policy is determined by selecting the discrete state closest to the observed demand value. This approach allows simulating operator behavior in realistic scenarios where demand varies continuously, maintaining consistency with strategies learned during training on discrete states.

\section{Simulations and Results}
\label{sec:simulazioni}

\subsection{Setup, Hypotheses and Scenarios}

In this section, operator behavior is analyzed through numerical simulations, with the objective of evaluating the impact of market rules and bidding strategies on main economic indicators. The following \emph{hypotheses} are made:
\begin{enumerate}
 \item Single-zone market without network, in order to focus purely on the economic side, i.e., market structure.
 \item Rigid and known demand, in order to focus on supply dynamics and equilibrium prices without introducing additional complexity related to demand flexibility.
 \item Operators do not know the costs and capacities of competitors, reflecting the market reality where information is incomplete and operators must base their strategies on historical data and observations.
\end{enumerate}
It should be noted that neglecting demand elasticity (an assumption often adopted in literature) risks overestimating the impact of bidding strategies on equilibrium price, since in a real market demand could reduce in response to high prices, mitigating---at least, partially---the effects of strategic offers. Furthermore, under these conditions, market optimization problems no longer maximize social welfare according to Samuelson's theory, but just \emph{minimize total market cost}.

\smallskip
\noindent
Two distinct scenarios are considered: the first, based on the assumptions of the PNIEC 2030, studies a simplified market with two operators and representative technological portfolios. The second scenario rather uses real data from public offers of 2024 to build a simplified representation of the Italian market with ten operators, each characterized by a portfolio obtained through clustering. In both cases, results of the three market regimes (Pay-as-Bid, Pay-as-Clear, Segmented Pay-as-Clear) are compared, both under perfect competition conditions (offers at marginal cost), and in the presence of strategic behaviors learned by agents through RL. The analysis focuses on how different market rules and bidding strategies influence the consumer-side price, total production cost, and profit distribution among operators, providing quantitative indications on allocation efficiency and participation incentives in the new segmented model.

\smallskip
\noindent
Simulations were implemented in Julia~\cite{bezansonJuliaFreshApproach2017}; modeling of bilevel optimization problems for the SPaC market was performed via the BilevelJuMP package~\cite{garciaBilevelJuMPjlModelingSolving2022} thanks to the work of Fabrizio Lacalandra~\cite{frangioniBilevelProgrammingApproach2024}. The code used for the simulations object of this report is publicly available on GitHub \cite{report_code}.

\subsection{Scenario 1: PNIEC 2030}

According to the Italian National Integrated Energy and Climate Plan (PNIEC) 2030, renewable sources will cover approximately 39.4\% of gross final energy consumption in Italy, while the remaining share will be guaranteed by traditional sources (natural gas, oil, coal). Thus, two symmetric market operators are considered, with comparable portfolios: approximately 40\% of capacity from NMCS plants and 60\% from NNMCS plants. Each operator owns five production units (PU) of different types, identical in capacity between the two subjects. The overall power for each operator is set at 1~GW, thus obtaining a reference market with total capacity equal to 2~GW. The only difference between operators lies in the marginal cost values of the units, chosen based on literature \cite{StudyLevelizedCost} and differentiated to reflect a minimal degree of heterogeneity. The portfolios of the two operators in the PNIEC 2030 scenario are reported in Table~\ref{tab:portafoglio-operatori-pniec}.

\begin{table}[H]
    \centering
    \caption{Operator portfolios in the PNIEC 2030 scenario. Marginal cost ranges are taken from \cite{StudyLevelizedCost}.}
    \label{tab:portafoglio-operatori-pniec}
    \resizebox{0.6\textwidth}{!}{
    \begin{tabular}{lllccl}
        \toprule
        Operator & PU & Type & Marginal cost [\euro/MWh] & Capacity [MW] \\
        \midrule
        OpA & PV    & FCMT  & 4.2  & 120 \\
        OpA & WIND  & FCMT  & 5.0  & 120 \\
        OpA & HYDRO & FCMT  & 12.0 & 160 \\
        OpA & GAS   & FCMNT & 94.0 & 420 \\
        OpA & COAL  & FCMNT & 149.0 & 180 \\
        \midrule
        OpB & PV    & FCMT  & 2.7  & 120 \\
        OpB & WIND  & FCMT  & 2.0  & 120 \\
        OpB & HYDRO & FCMT  & 20.0 & 160 \\
        OpB & GAS   & FCMNT & 69.0 & 420 \\
        OpB & COAL  & FCMNT & 139.0 & 180 \\
        \bottomrule
    \end{tabular}
    }
\end{table}

\smallskip
\noindent
The market is thus constituted by 10 production units overall, equally distributed between the two operators. Each operator holds a 50\% market share, thus none has a dominant position; this also holds in the two distinct sub-markets (NMCS and NNMCS) of SPaC.

\subsubsection{Simulation of Markets at Marginal Costs}

In this first analysis, the three market regimes (PaB, PaC and SPaC) are simulated assuming all operators offer at their marginal cost, as expected under perfect competition conditions. A rigid demand equal to $D = 1~\text{GW}$ is considered, corresponding to 50\% of total market capacity in the PNIEC 2030 scenario. Numerical results of the simulations are summarized in Table~\ref{tab:risultati-mc-1gw}, where the total system cost and the unique national price (PUN) resulting in the three market regimes are reported. 

\begin{figure}[htb]
    \centering
    \begin{subfigure}{0.48\textwidth}
        \centering
        \includegraphics[width=\textwidth]{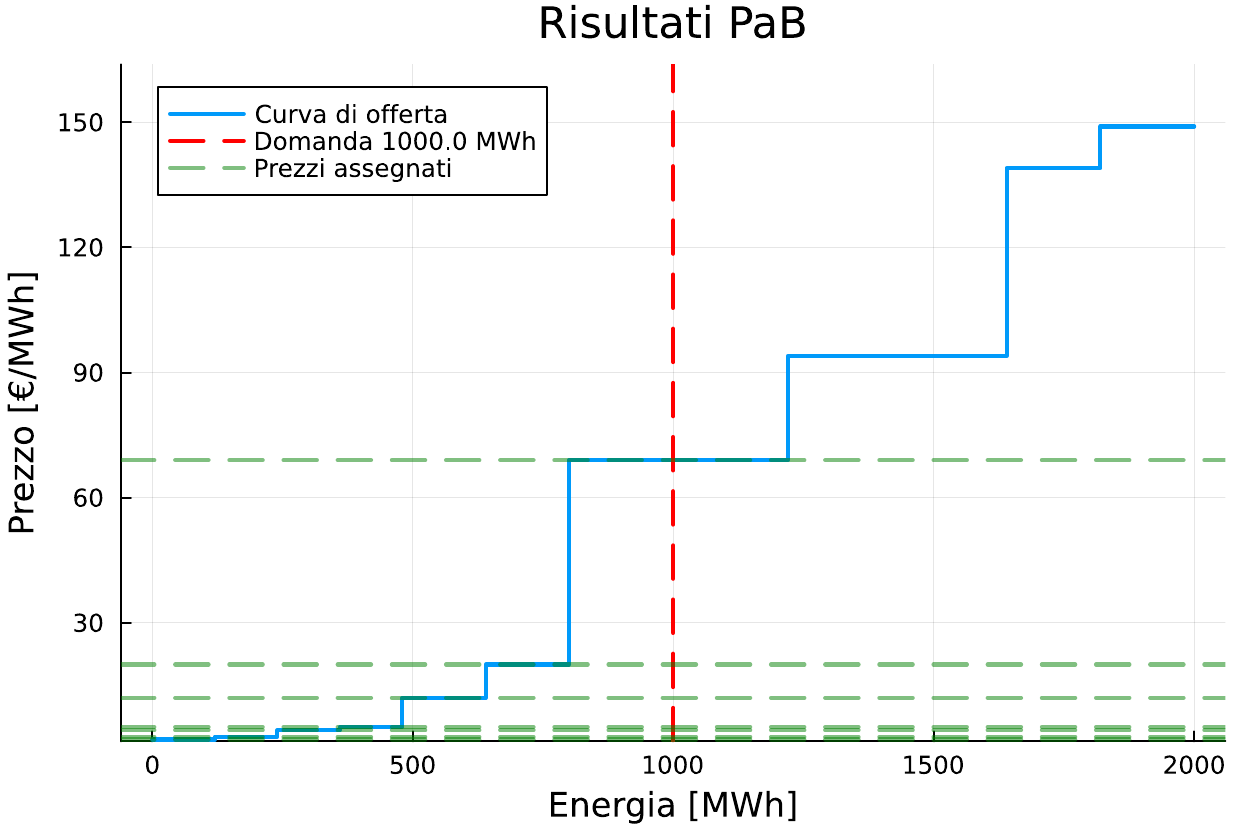}
        \caption{PaB: offers at marginal cost}
        \label{fig:risultati-pab-mc}
        \vspace{0.5cm}
        \includegraphics[width=\textwidth]{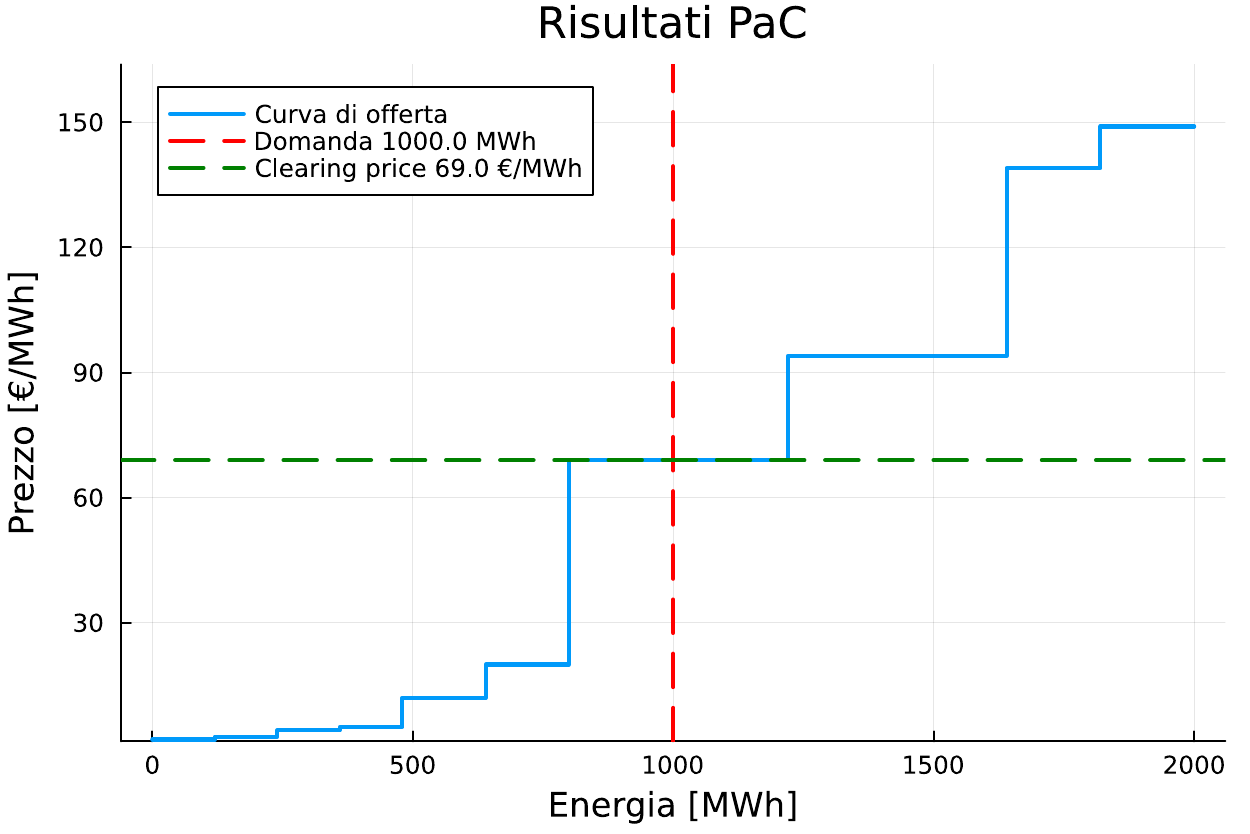}
        \caption{PaC: offers at marginal cost}
        \label{fig:risultati-pac-mc}
    \end{subfigure}%
    \hfill
    \begin{subfigure}{0.48\textwidth}
        \centering
        \includegraphics[height=12cm, width=\textwidth]{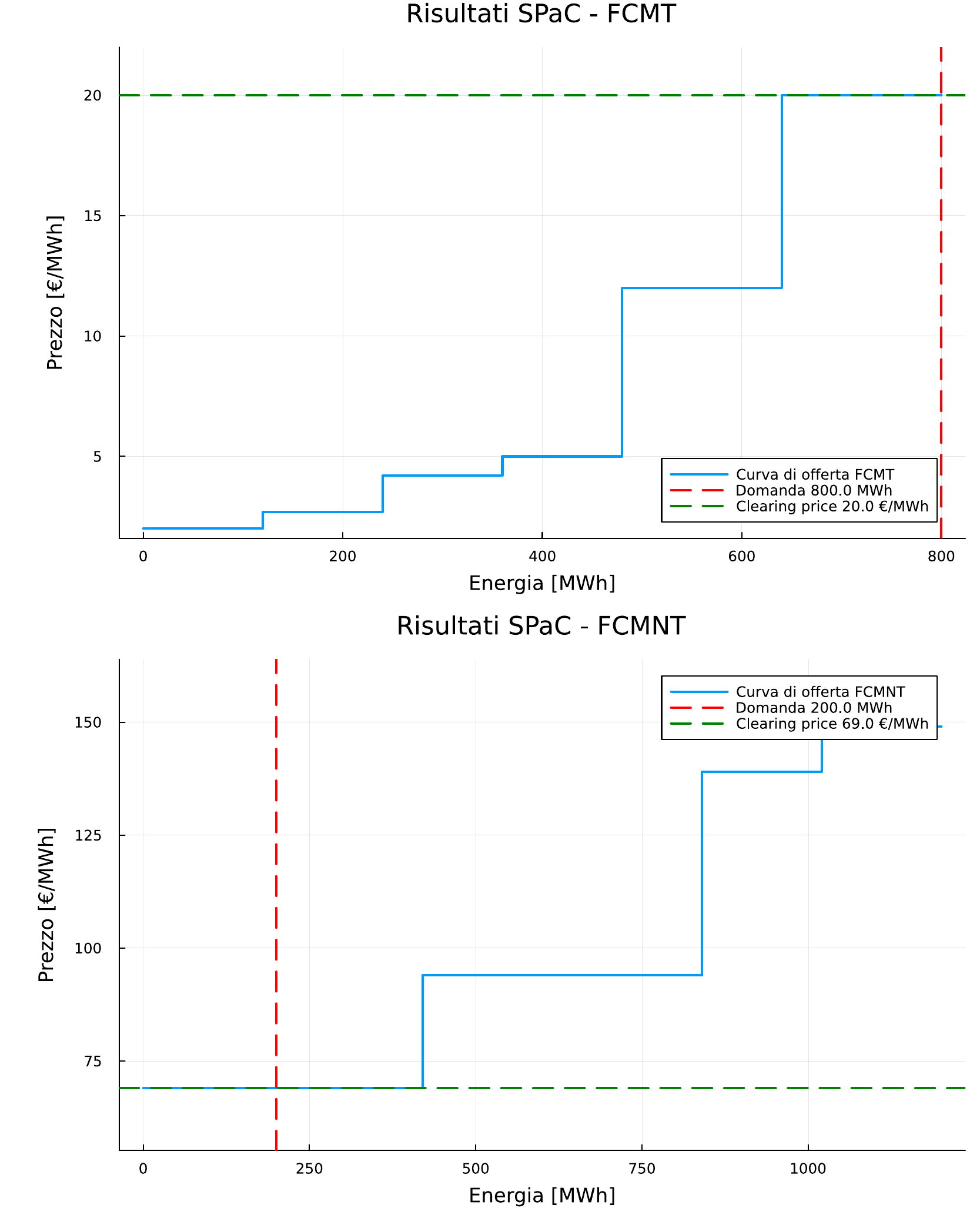}
        \caption{SPaC: offers at marginal cost}
        \label{fig:risultati-spac-mc}
    \end{subfigure}
    \caption{Supply curves and results of simulations with offers at marginal cost: (a) PaB, (b) PaC, (c) SPaC.}
    \label{fig:risultati-mc-confronto}
\end{figure}

\begin{table}[H]
    \centering
    \caption{Results of simulations of the three markets at marginal costs for $D = 1~\text{GW}$ (PNIEC 2030 scenario).}
    \label{tab:risultati-mc-1gw}
    \begin{tabular}{lcc}
        \toprule
        \textbf{Market} & \textbf{Total cost [\EUR{}]} & \textbf{PUN [\EUR{}/MWh]} \\
        \midrule
        PaB  & 20\,588 & 20.59 \\
        PaC  & 69\,000 & 69.00 \\
        SPaC & 29\,800 & 29.80 \\
        \bottomrule
    \end{tabular}
\end{table}

\noindent
Firstly, it is observed that the PaB regime leads to the lowest total cost (20\,588~\EUR{}) and PUN (20.59~\EUR{}/MWh). This result is consistent with theory, since in PaB each unit is remunerated at the offered price and not at the system marginal price: operators, offering at marginal cost, receive no extra profit. However, this configuration is completely unrealistic, since in practice operators have an incentive to apply strategic markup to avoid operational losses and guarantee economic sustainability of their assets. 

\smallskip
\noindent
The PaC market, on the contrary, shows a significantly higher system cost (69\,000~\EUR{}) and a PUN equal to 69.0~\EUR{}/MWh. In this case, all accepted units receive the price of the last marginal technology (gas), with the result that even renewable and hydroelectric sources with low marginal costs are remunerated at values much higher than their actual costs. This is the known "renewables paradox".

\smallskip
\noindent
The SPaC market is in an intermediate position, with a total cost equal to 29\,800~\EUR{} and an average PUN of 29.80~\EUR{}/MWh.  The segmentation between NMCS and NNMCS allows preserving marginal price signals for each technological category, while simultaneously reducing overall expenditure compared to traditional PaC. In particular, renewable and hydroelectric technologies are remunerated at 20~\EUR{}/MWh, in line with their actual marginal costs, while fossil sources remain valued at the higher level of 69~\EUR{}/MWh. This result suggests that SPaC is able to limit inframarginal rents of low-cost technologies.

\smallskip
\noindent
From the point of view of economic efficiency, it emerges that SPaC captures about 71\% of the benefits in terms of cost reduction compared to the transition from PaC to PaB (reduction of 39\,200~\EUR{} vs. theoretical maximum reduction of 48\,412~\EUR{}). This is particularly significant considering that SPaC maintains differentiated price signals for different technologies, a crucial aspect for correctly orienting future investments in the electricity system. These results highlight the different trade-offs between economic efficiency and quality of price signals generated by the three market regimes, providing the starting point for the subsequent analysis where strategic behaviors by operators is introduced.

\subsubsection{Simulation of Markets with Random Fixed Markups}

In this second analysis, operators no longer offer at marginal cost: instead, each operator applies random fixed markups in the range 0--20\% for each technology in their portfolio, to be applied to the respective marginal cost. This simulation aims to evaluate how the three market regimes respond to non-optimized strategic behaviors, representative of a situation where operators apply arbitrary markups without knowing the optimal strategy. The applied markup values are reported in Table~\ref{tab:markup-casuali}. The choice of markups reflects a strategic asymmetry between the two operators: OpA adopts a more ``aggressive'' strategy with generally higher markups (average: 10.6\%), while OpB applies a more conservative approach (average: 8.0\%), except for photovoltaic technology where it is more ``aggressive''.

\begin{table}[H]
    \centering
    \caption{Random markups applied.}
    \label{tab:markup-casuali}
    \begin{tabular}{lccccc}
    \toprule
    Operator & COAL & GAS & HYDRO & PV & WIND \\
    \midrule
    OpA & 0\% & 12\% & 16\% & 15\% & 10\% \\
    OpB & 0\% & 7\% & 11\% & 18\% & 4\% \\
    \bottomrule
    \end{tabular}
\end{table}

\begin{figure}[htb]
    \centering
    \begin{subfigure}{0.48\textwidth}
    \centering
    \includegraphics[width=\textwidth]{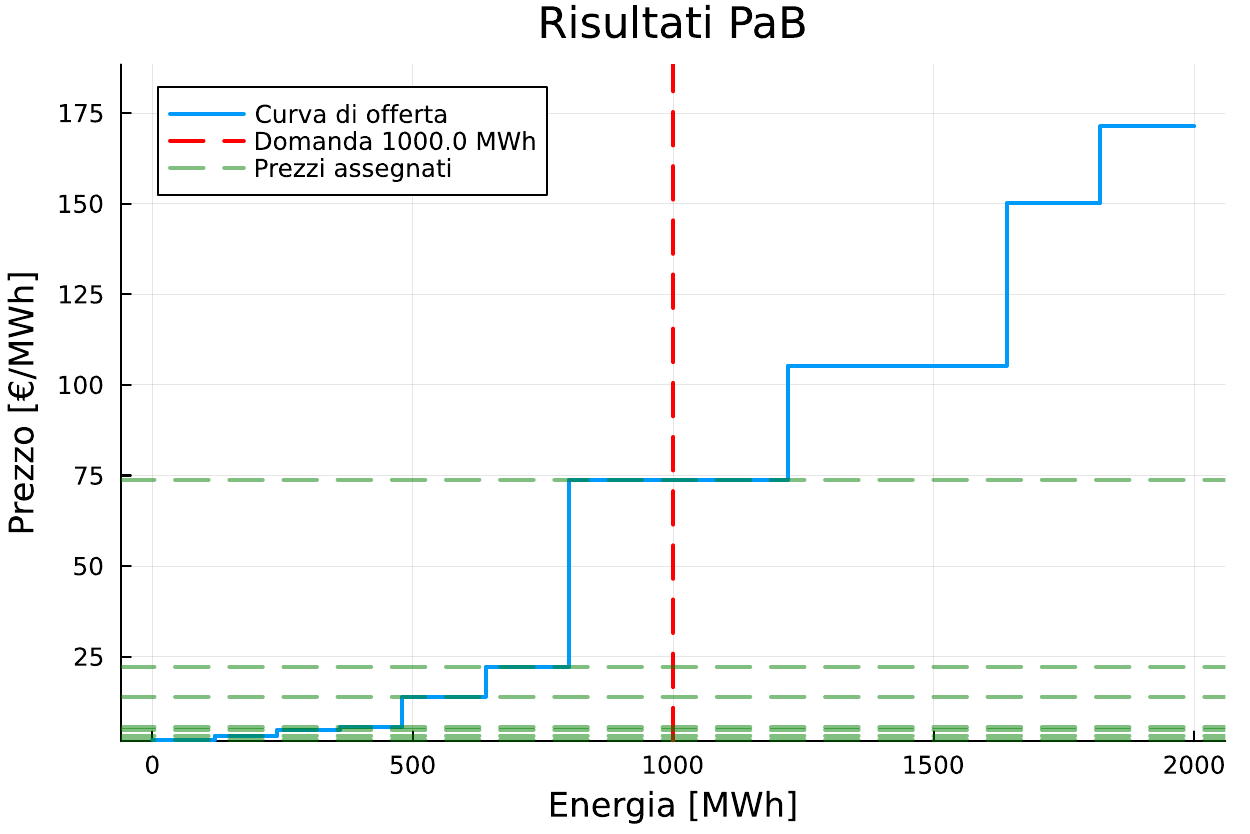}
    \caption{PaB: offers with random markups}
    \label{fig:risultati-pab-markup}
    \vspace{0.5cm}
    \includegraphics[width=\textwidth]{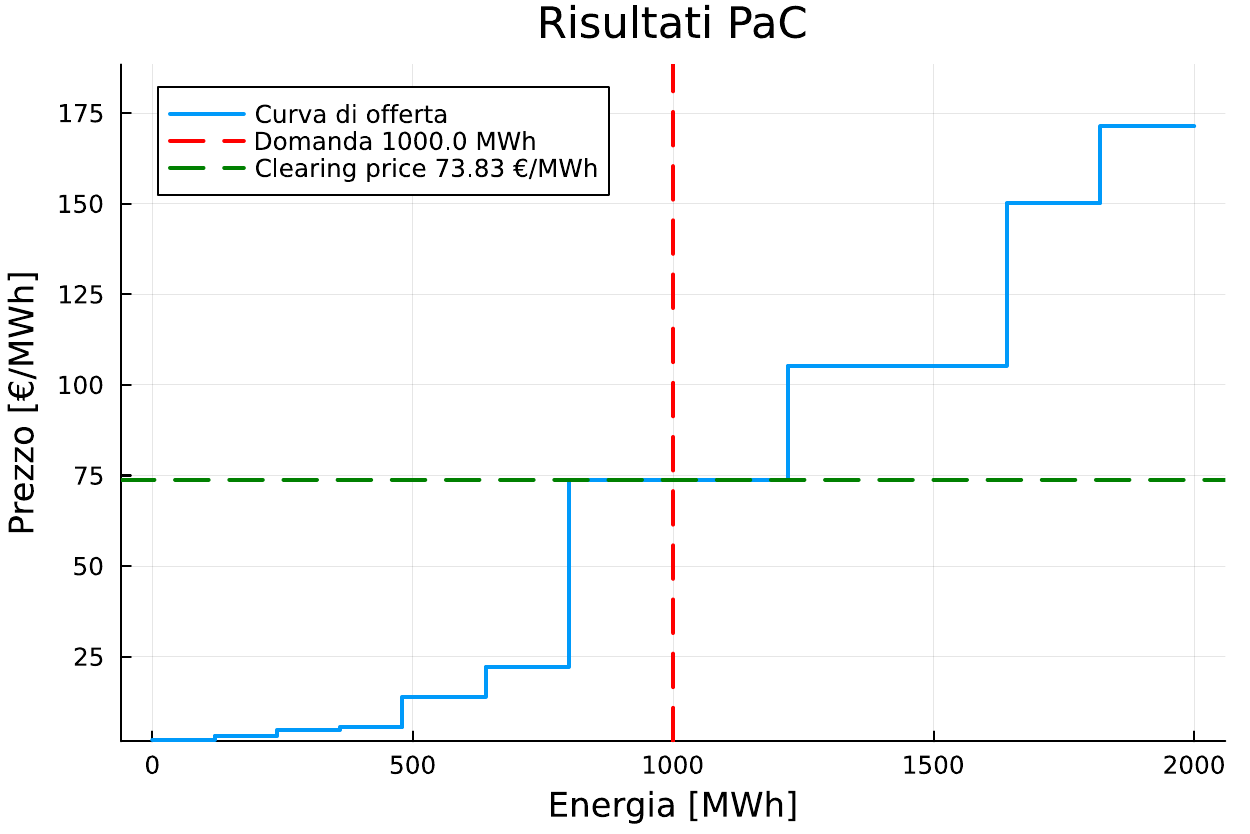}
    \caption{PaC: offers with random markups}
    \label{fig:risultati-pac-markup}
    \end{subfigure}%
    \hfill
    \begin{subfigure}{0.48\textwidth}
    \centering
    \includegraphics[height=12cm, width=\textwidth]{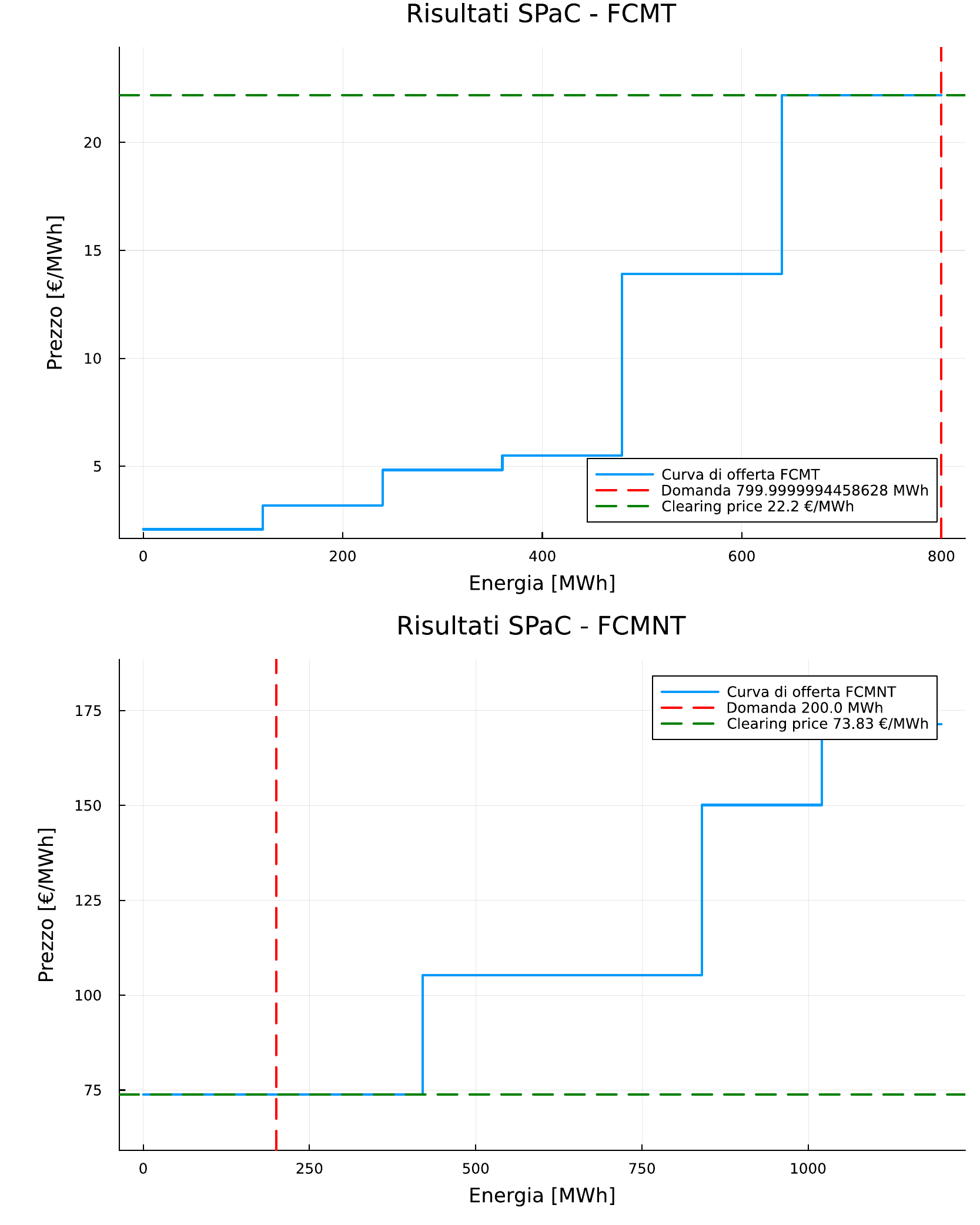}
    \caption{SPaC: offers with random markups}
    \label{fig:risultati-spac-markup}
    \end{subfigure}
    \caption{Results of the three markets with offers at random markups (Table \ref{tab:markup-casuali}): left PaB and PaC, right SPaC.}
    \label{fig:risultati-markup-confronto}
\end{figure}

\noindent
Results of simulations of the three markets with random fixed markups for a demand of 1~GW are shown in Figure \ref{fig:risultati-markup-confronto} and reported in Table~\ref{tab:risultati-markup-1gw}. 

\begin{table}[H]
    \centering
    \caption{Results of simulations of the three markets with random markups (Table \ref{tab:markup-casuali}) for D = 1~GW}
    \label{tab:risultati-markup-1gw}
    \begin{tabular}{lcc}
        \toprule
        \textbf{Market} & \textbf{Total cost [\EUR{}]} & \textbf{PUN [\EUR{}/MWh]} \\
        \midrule
        PaB  & 22\,417 & 22.42 \\
        PaC  & 73\,830 & 73.83 \\
        SPaC & 32\,526 & 32.53 \\
        \bottomrule
    \end{tabular}
\end{table}

\smallskip
\noindent
In the SPaC case, the split between categories is as follows: 800~MWh accepted from NMCS units (total cost 17\,760~\EUR{} and PUN 22.2~\EUR{}/MWh), 200~MWh from NNMCS units (total cost 14\,766~\EUR{} and PUN 73.83~\EUR{}/MWh). As expected, introducing positive markups causes a generalized increase in both total costs and prices in all three regimes. However, the extent of such increases varies significantly between different price formation mechanisms. Table~\ref{tab:confronto-costi-markup} summarizes variations compared to the marginal cost case.

\begin{table}[H]
    \centering
    \resizebox{0.95\textwidth}{!}{
    \begin{tabular}{lcccc}
    \toprule
    \textbf{Market} & \textbf{Cost (marg.) [\EUR{}]} & \textbf{Cost (random markup) [\EUR{}]} & \textbf{$\Delta$C [\%]} & \textbf{$\Delta$ PUN [\EUR{}/MWh]} \\
    \midrule
    PaB  & 20\,588 & 22\,417 & +8.9\% & +1.83 \\
    PaC  & 69\,000 & 73\,830 & +7.0\% & +4.83 \\
    SPaC & 29\,800 & 32\,526 & +9.1\% & +2.73 \\
    \bottomrule
    \end{tabular}
    }
    \caption{Comparison between results at marginal costs and with random markups (D = 1~GW)}
    \label{tab:confronto-costi-markup}
\end{table}

\noindent
The PaB regime shows the greatest sensitivity to applied markups (+8.9\%), a phenomenon that reflects the direct nature of price formation in this mechanism: each increase in offer costs translates immediately into a proportional increase in total expenditure. Conversely, PaC presents the smallest relative variation (+7.0\%), since the effect of markups is partially attenuated by the fact that all accepted units still receive the marginal unit price.

\smallskip
\noindent
The SPaC market registers an increase (+9.1\%) similar to PaB, highlighting how segmentation maintains a direct relationship between offer strategies and resulting costs within each segment. It is particularly significant that in the NMCS segment the price increased from 20.0 to 22.2~\EUR{}/MWh (+11\%), while in the NNMCS segment the increase was from 69.0 to 73.83~\EUR{}/MWh (+7\%) in line with markups applied by operators.

\smallskip
\noindent
From the point of view of relative efficiency, SPaC maintains its intermediate position even in the presence of strategic behaviors: its total cost results 56\% lower than PaC (against 57\% in the marginal cost case) and 45\% higher than PaB (same as the previous case). This stability of relative ratios suggests that allocation benefits of SPaC persist even when operators adopt sub-optimal strategies.

\subsubsection{Training of Operators for a Demand Level}

In the third analysis, operators are trained via RL to learn optimal bidding strategies in the three market regimes, considering a fixed demand level equal to 1~GW. The objective is to determine optimal markups that maximize each operator's profit, taking into account strategic interactions and specific price formation rules. The action space is defined by discrete markups $m \in \{0\%, 5\%, 10\%, 20\%\}$ for PaC and SPaC markets, while for PaB higher values are considered given the clearing mechanism nature. Training was conducted for 2,000 episodes with Q-learning algorithm and $\varepsilon$-greedy policy, using exponential decay of the exploration parameter from $\varepsilon_{\text{max}} = 1.0$ to $\varepsilon_{\text{min}} = 0.05$.

\begin{figure}[hp]
    \centering
    \begin{subfigure}{\textwidth}
        \centering
        \begin{subfigure}{0.48\textwidth}
            \centering
            \includegraphics[height=5cm, width=\textwidth]{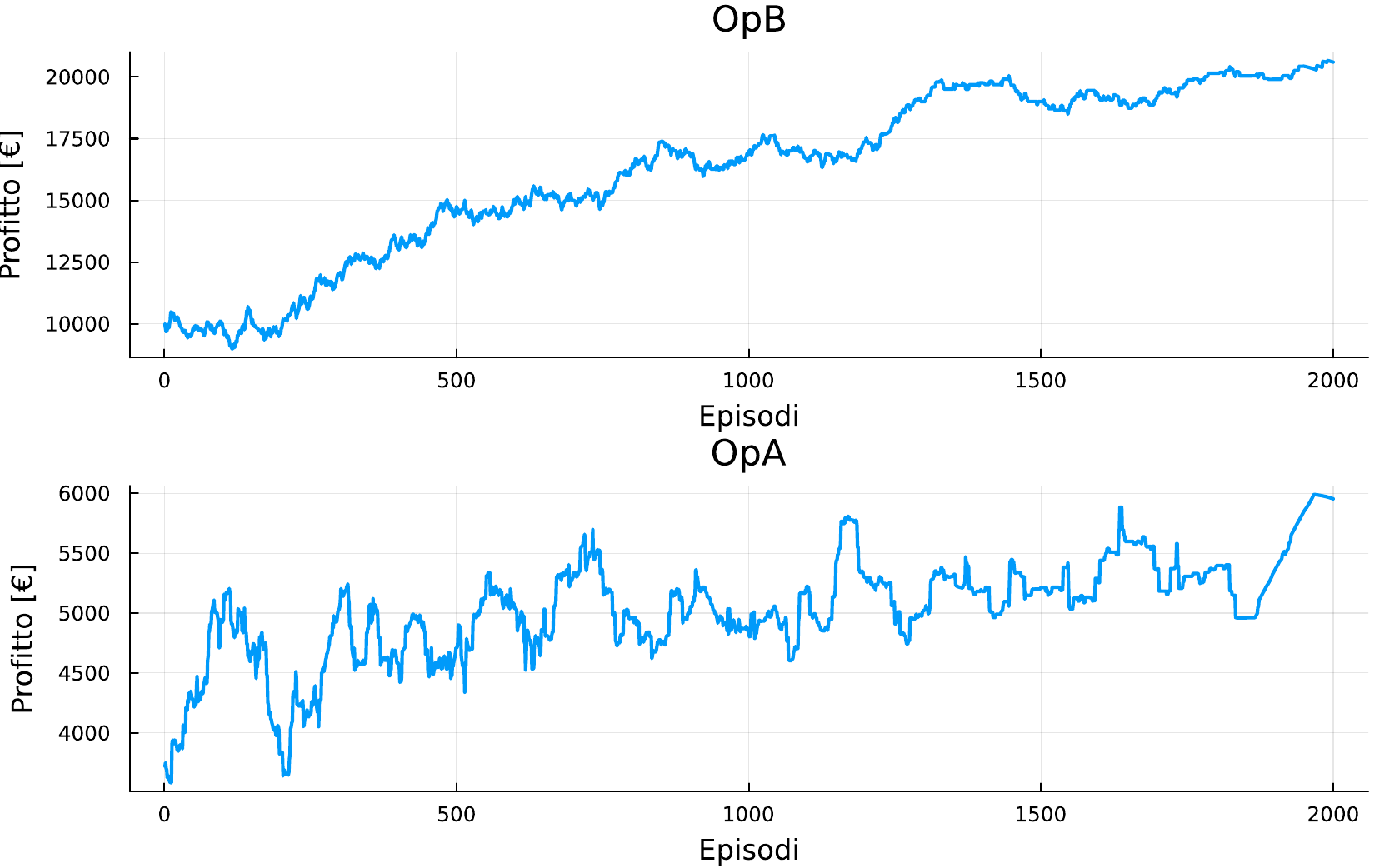}
            \caption*{Operator profits}
        \end{subfigure}%
        \hfill
        \begin{subfigure}{0.48\textwidth}
            \centering
            \includegraphics[height=5cm, width=\textwidth]{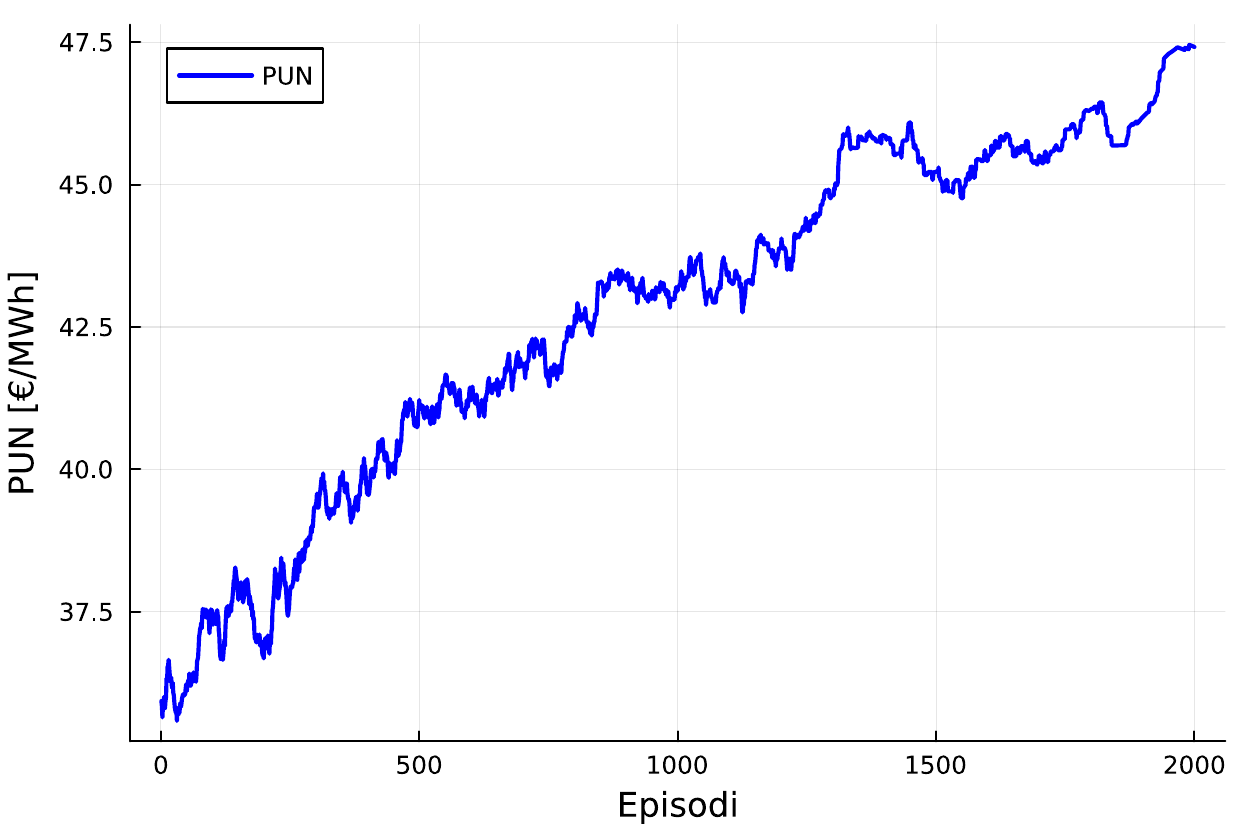}
            \caption*{PUN}
        \end{subfigure}
        \caption{PaB market.}
    \end{subfigure}
    
    \vspace{0.5cm}
    
    \begin{subfigure}{\textwidth}
        \centering
        \begin{subfigure}{0.48\textwidth}
            \centering
            \includegraphics[height=5cm, width=\textwidth]{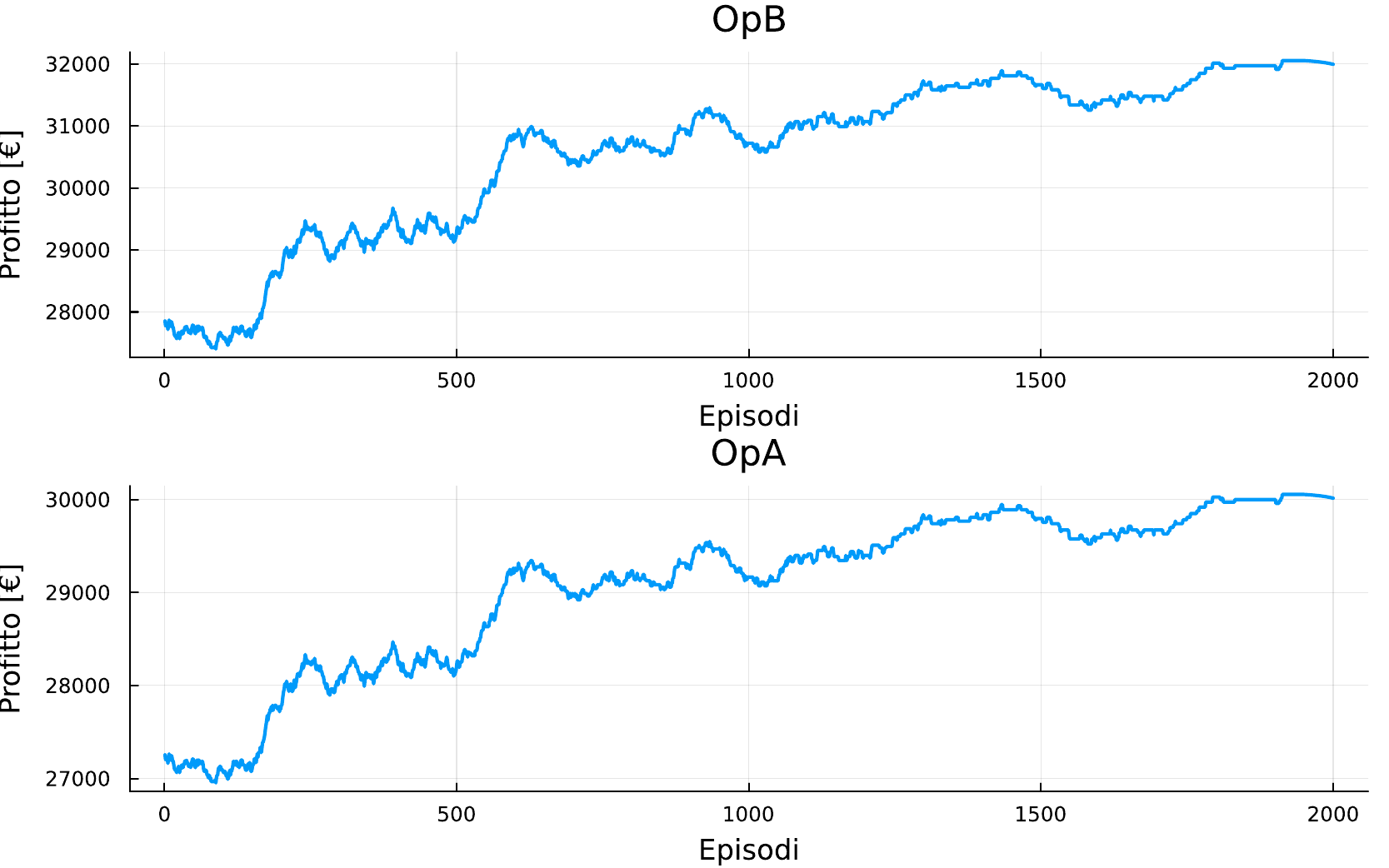}
            \caption*{Operator profits}
        \end{subfigure}%
        \hfill
        \begin{subfigure}{0.48\textwidth}
            \centering
            \includegraphics[height=5cm, width=\textwidth]{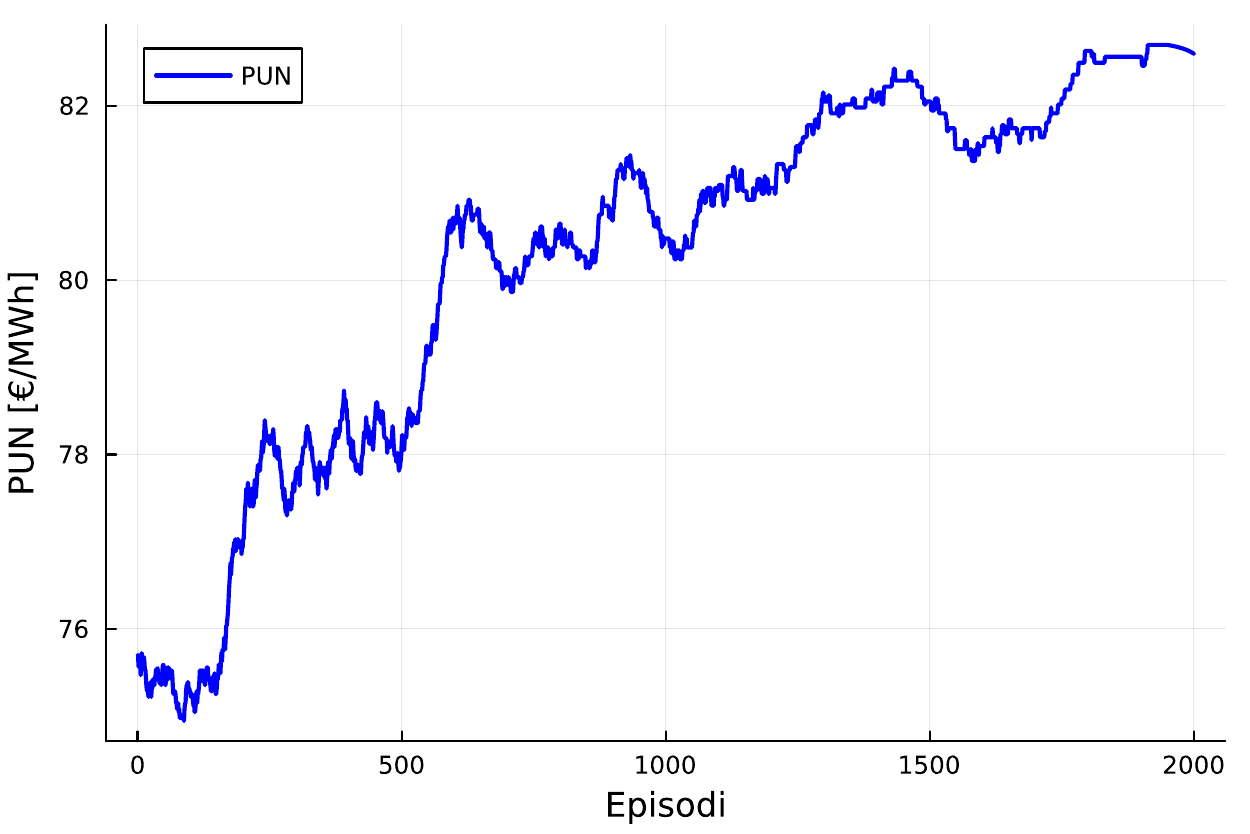}
            \caption*{PUN}
        \end{subfigure}
        \caption{PaC market.}
    \end{subfigure}
    
    \vspace{0.5cm}
    
    \begin{subfigure}{\textwidth}
        \centering
        \begin{subfigure}{0.48\textwidth}
            \centering
            \includegraphics[height=5cm, width=\textwidth]{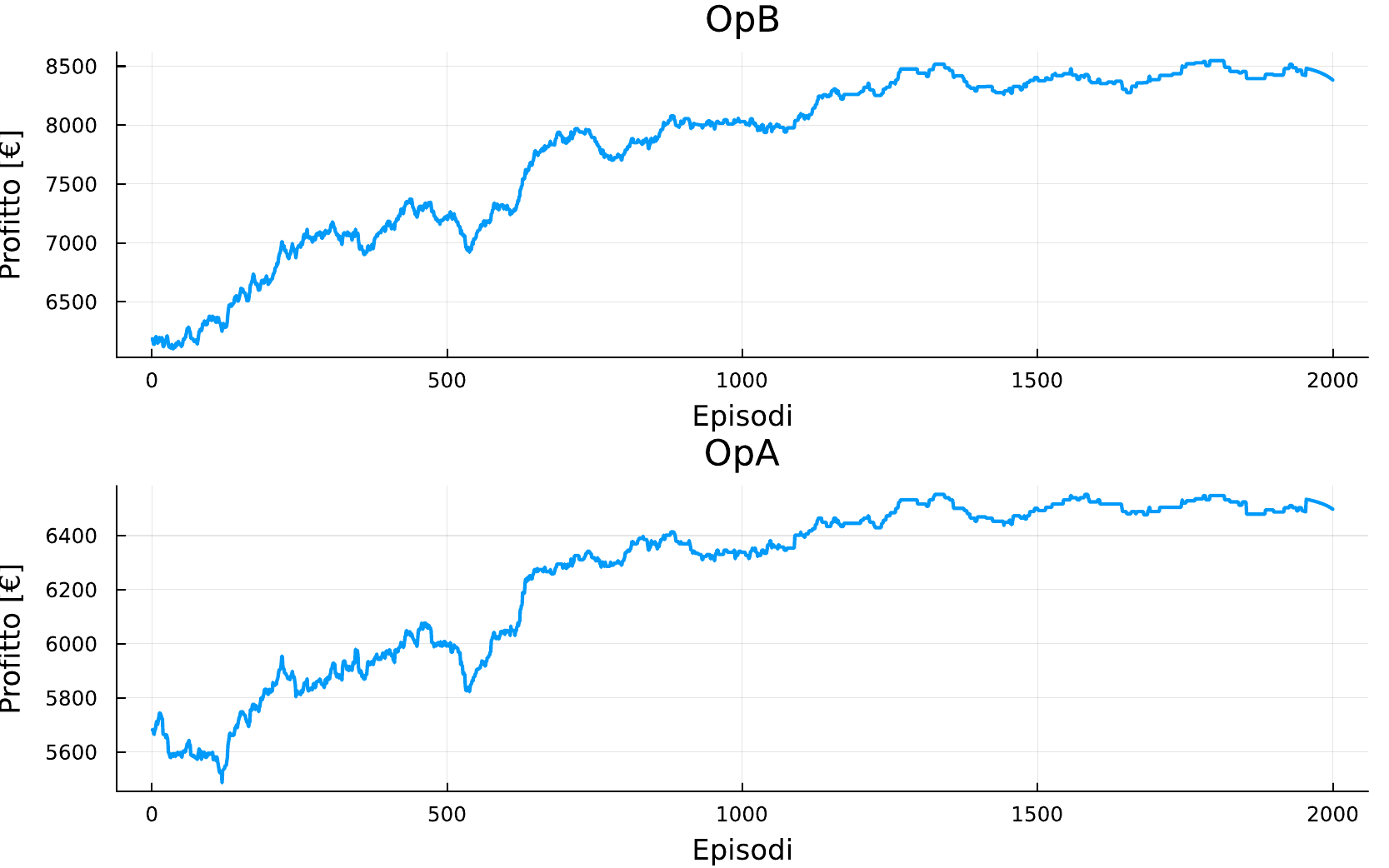}
            \caption*{Operator profits}
        \end{subfigure}%
        \hfill
        \begin{subfigure}{0.48\textwidth}
            \centering
            \includegraphics[height=5cm, width=\textwidth]{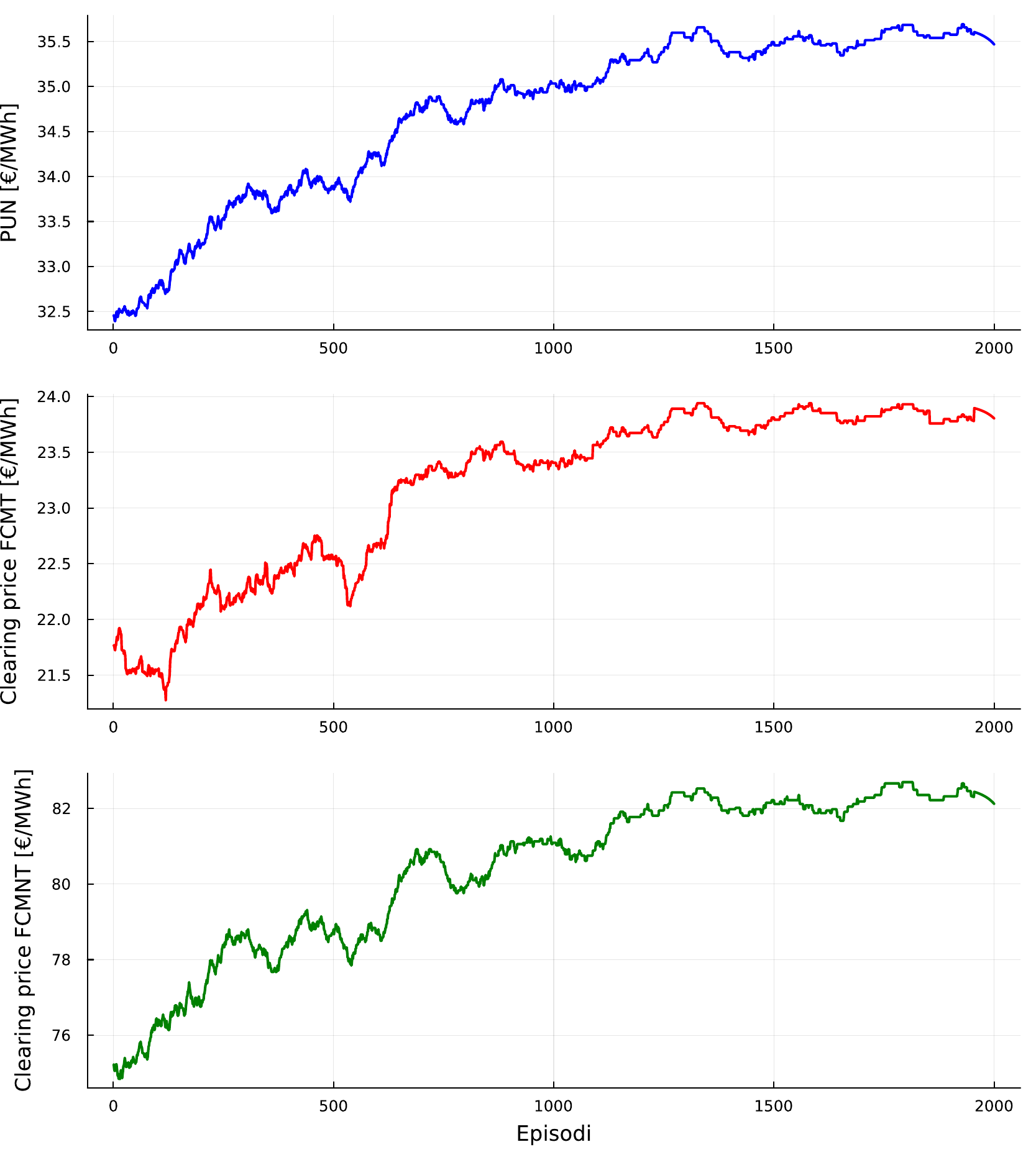}
            \caption*{PUN}
        \end{subfigure}
        \caption{SPaC market.}
    \end{subfigure}
            
    \caption{Monitoring of training for the three markets and demand 1.0~GW. Left: operator profits, right: PUN.}
    \label{fig:monitoraggio-training-rl-tutti}
\end{figure}

\smallskip
\noindent
The graphs in Figure \ref{fig:monitoraggio-training-rl-tutti} show that, after an initial exploration phase ($\varepsilon = 1$), operators converge towards stable strategies that seek to maximize their profits.

\smallskip
\noindent
In the \emph{PaB} case operators learn extremely high markups (up to 200\%) for renewable and hydroelectric technologies, while applying lower (or zero) markups to fossil technologies. This strategy reflects the incentive to maximize remuneration of low marginal cost units that have higher probability of being accepted.

\smallskip
\noindent
In the \emph{PaC} case optimal markups are significantly lower (0--20\%), with a more balanced distribution among technologies. Operators learn to balance offer competitiveness with the need to obtain positive margins, considering that all accepted units receive the marginal price. The PUN follows a growing trend, although less steep than in PaB, stabilizing after about 1500 episodes. The outcome is consistent with marginalist clearing logic: operators quickly learn that marginal fossil sources can apply higher markups, while negligible cost technologies tend to be offered with values close to marginal cost. The result is a more predictable equilibrium but with still significant rents for marginal units.

\smallskip
\noindent
Finally, in \emph{SPaC} case (Figure~\ref{fig:monitoraggio-training-rl-tutti}c), training shows a more articulated dynamic: operator profits grow less and are distributed heterogeneously among different technologies. The PUN stabilizes at intermediate values, lower compared to PaC and PaB, highlighting the segmentation effect. The SPaC model indeed limits markup possibilities on renewable and hydroelectric sources, which remain remunerated at prices close to costs, while allowing lower profit margins to fossil technologies. This translates into an equilibrium where operators do not maximize profits as in PaB, but the system overall benefits from cost reduction and greater allocation efficiency.

\smallskip
\noindent
The optimal markup strategies (in percentage relative to marginal cost) learnt by the operators at the end of the training process are reported in Table \ref{tab:policy-ottimale-rl} for each technological cluster and market regime.

\begin{table}[H]
    \centering
    \caption{Optimal policy learned by operators via RL for each market and technology; values are the optimal technology-specific markups ($m$).}
    \label{tab:policy-ottimale-rl}
    \resizebox{0.6\textwidth}{!}{
    \begin{tabular}{llccccc}
    \toprule
    Market & Operator & $m_{\text{HYDRO}}$ [\%] & $m_{\text{WIND}}$ [\%] & $m_{\text{PV}}$ [\%] & $m_{\text{GAS}}$ [\%] & $m_{\text{COAL}}$ [\%] \\  \midrule
    \multirow{2}{*}{PaB} 
        & OpB & 200.0 & 200.0 & 100.0 & 100.0 & 50.0 \\
        & OpA & 200.0 & 200.0 & 200.0 & 200.0 & 0.0 \\
    \midrule
    \multirow{2}{*}{PaC} 
        & OpB & 10.0 & 0.0 & 0.0 & 20.0 & 0.0 \\
        & OpA & 0.0 & 10.0 & 20.0 & 20.0 & 5.0 \\
    \midrule
    \multirow{2}{*}{SPaC} 
        & OpB & 20.0 & 0.0 & 5.0 & 20.0 & 20.0 \\
        & OpA & 10.0 & 20.0 & 0.0 & 0.0 & 0.0 \\
    \bottomrule
    \end{tabular}
    }
\end{table}

\noindent
Notably, that optimal strategies vary significantly between different market regimes. Results of simulations with learned strategies (Table~\ref{tab:riepilogo-rl-singolo-stato}) show significant cost increases compared to previous marginal cost simulations.

\begin{table}[H]
    \centering
    \caption{Summary RL results for a single demand level ($D = 1$~GW).}
    \label{tab:riepilogo-rl-singolo-stato}
    \begin{tabular}{lcc}
    \toprule
    \textbf{Market} & \textbf{Total cost [\euro]} & \textbf{PUN [\euro/MWh]} \\
    \midrule
    PaB  & 47\,640 & 47.64 \\
    PaC  & 82\,800 & 82.80 \\
    SPaC & 35\,760 & 35.76 \\
    \bottomrule
    \end{tabular}
\end{table}

\noindent
The comparison highlights total cost increases of 131\% for PaB (from 20,588 to 47,640~\EUR{}), 20\% for PaC (from 69,000 to 82,800~\EUR{}) and 20\% for SPaC (from 29,800 to 35,760~\EUR{}). In this case, SPaC is the lowest cost market, with a total cost 57\% lower than PaC and 25\% lower than PaB. This confirms the effectiveness of the segmentation mechanism even in the presence of optimized strategic behaviors.

\subsubsection{Training of Operators on Multiple Demand Levels}

\begin{figure}[hp]
    \centering
    \begin{subfigure}{0.8\textwidth}
        \centering
        \includegraphics[scale=0.2]{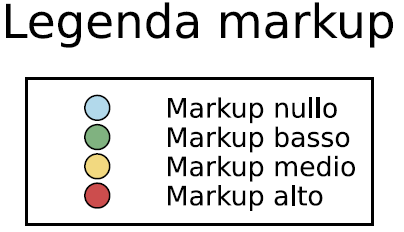}
    \end{subfigure}
    
    \vspace{1cm} 
    
    \begin{subfigure}{0.9\textwidth}
        \centering
        \includegraphics[height=0.28\textheight]{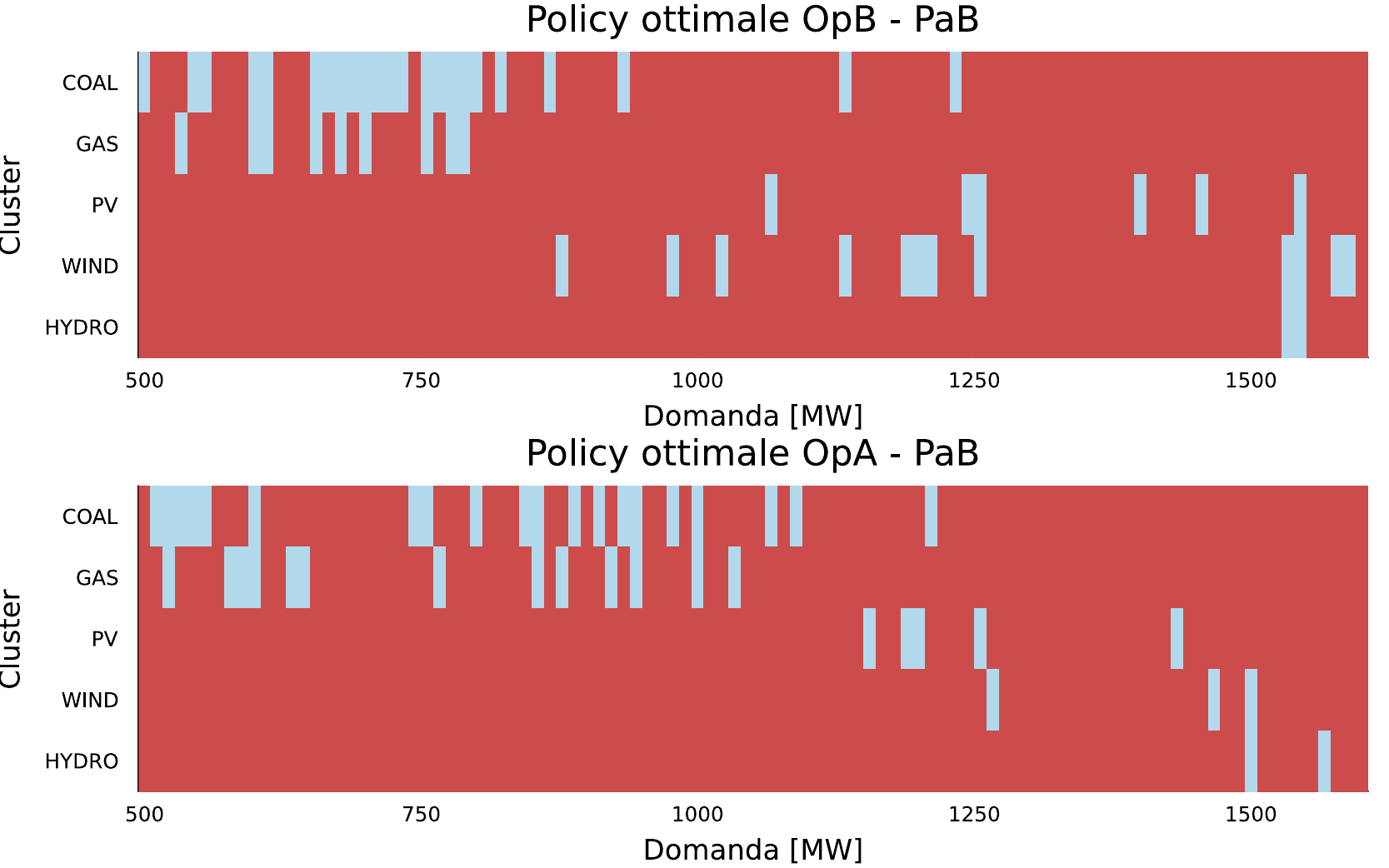}
        \label{fig:policy-pab}
    \end{subfigure}
    \vspace{0.3cm}
    \begin{subfigure}{0.9\textwidth}
        \centering
        \includegraphics[height=0.28\textheight]{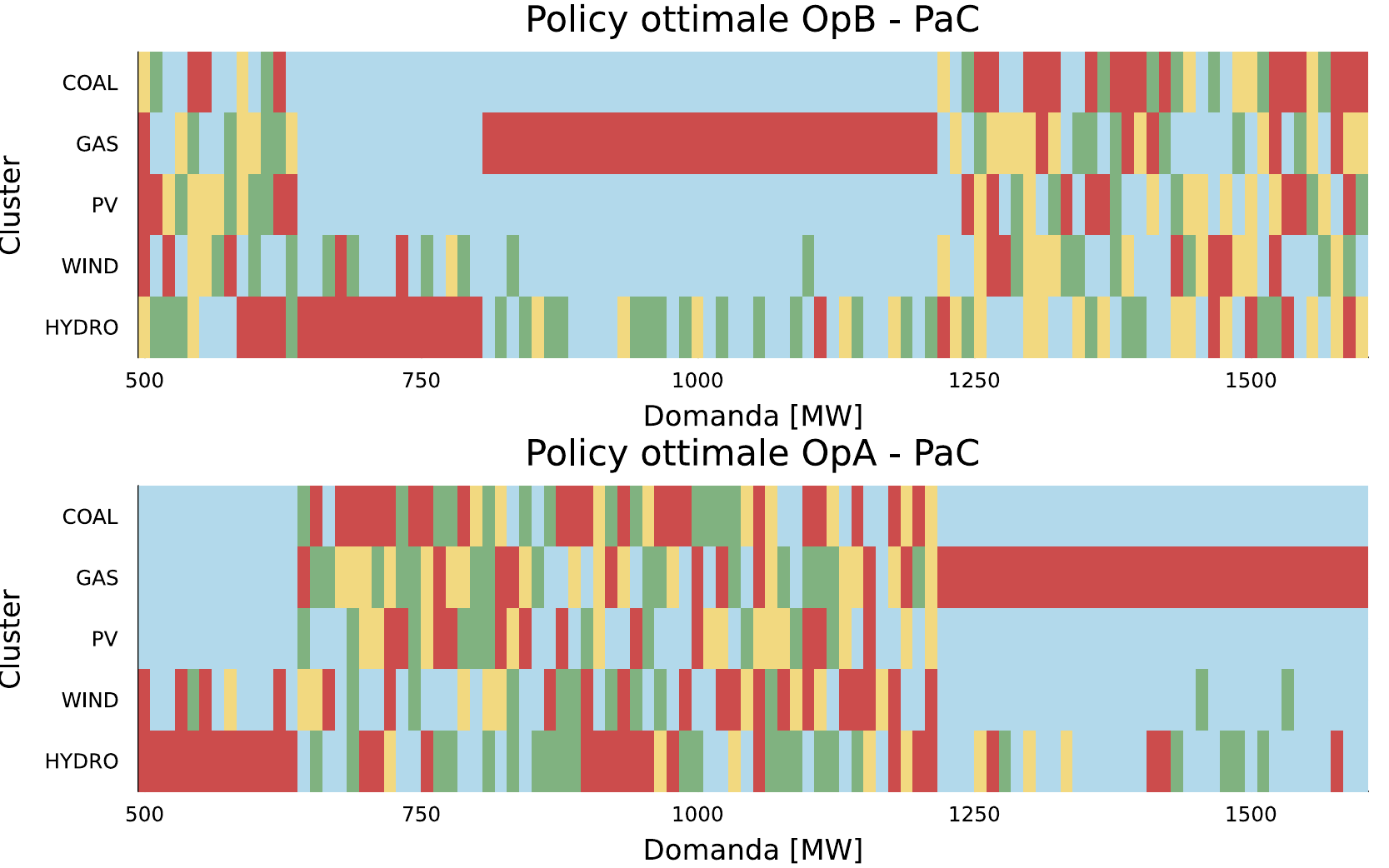}
        \label{fig:policy-pac}
    \end{subfigure}
    \vspace{0.3cm}
    \begin{subfigure}{0.9\textwidth}
        \centering
        \includegraphics[height=0.28\textheight]{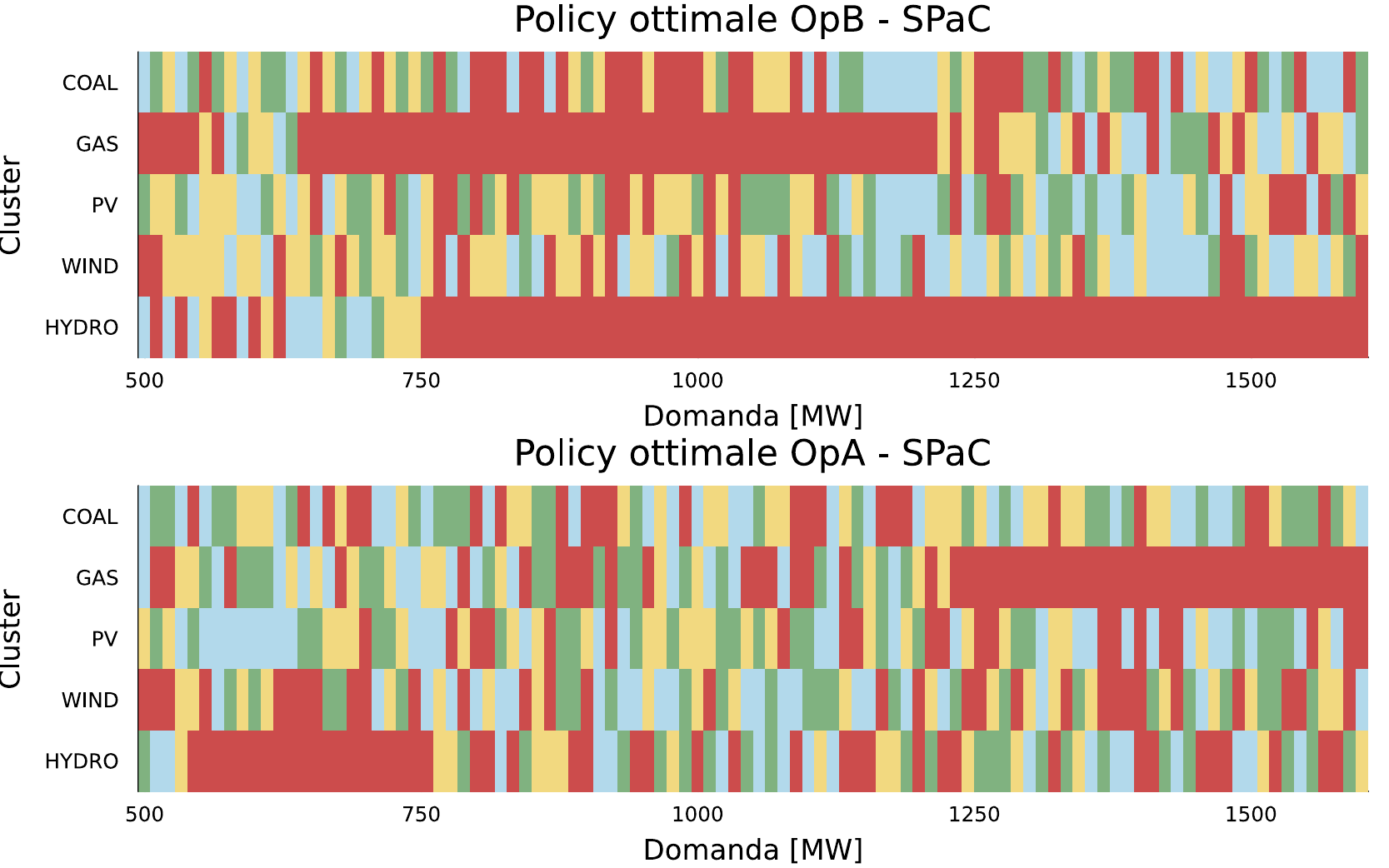}
        \label{fig:policy-spac}
    \end{subfigure}
    \caption{Optimal policy learned by operators via RL for the three markets on multiple demand levels.}
    \label{fig:policy-rl-multidemand}
\end{figure}

The previous approach is extended to train operators on multiple demand levels in order to simulate a typical day in the MGP. Demand is varied between \emph{25\% and 80\% of total market capacity} (2~GW) with 100 intermediate equidistant steps, thus increasing the state space from 1 to 100. Hyperparameters used are the same as previous training. The training process is repeated for each market and for each demand level (state), for a total of $2,000 \times 100 = 200,000$ episodes solved in parallel thanks to Julia multithreading on 13 cores of an Intel i7-13700H processor. The simulation can be performed in parallel for each demand level, since states (demand levels) are independent from each other due to assumptions made, and each operator chooses its action (markup) only as a function of the current state. This solution, besides accelerating training, allows reducing memory requirements to store the Q-table since each new value replaces the previous one. Training with this setup required about 15 minutes of which 14 minutes only for the SPaC market, the most complex to solve with Gurobi. Policies obtained for each operator in the three markets are represented in Figure~\ref{fig:policy-rl-multidemand}.

\smallskip
\noindent
Figure \ref{fig:policy-rl-multidemand} shows that in the PaB market operators tend to apply very high markups for all technologies, regardless of demand level. In the PaC market, instead, the two operators strategically apply high markups only for technologies that result marginal. Indeed, it can be seen that for low demand levels both apply high markups for hydroelectric, while this behavior shifts to gas technologies when demand increases, since these then become marginal. The SPaC market has a similar behavior to PaC, with the difference that high markups are applied simultaneously to hydroelectric and gas technologies for low demand levels, i.e., to the two technologies that are marginal in the respective segments (NMCS and NNMCS).

\subsubsection{Simulation of a Typical Day in the MGP}

\begin{figure}[htb]
    \centering
    \includegraphics[width=0.6\textwidth]{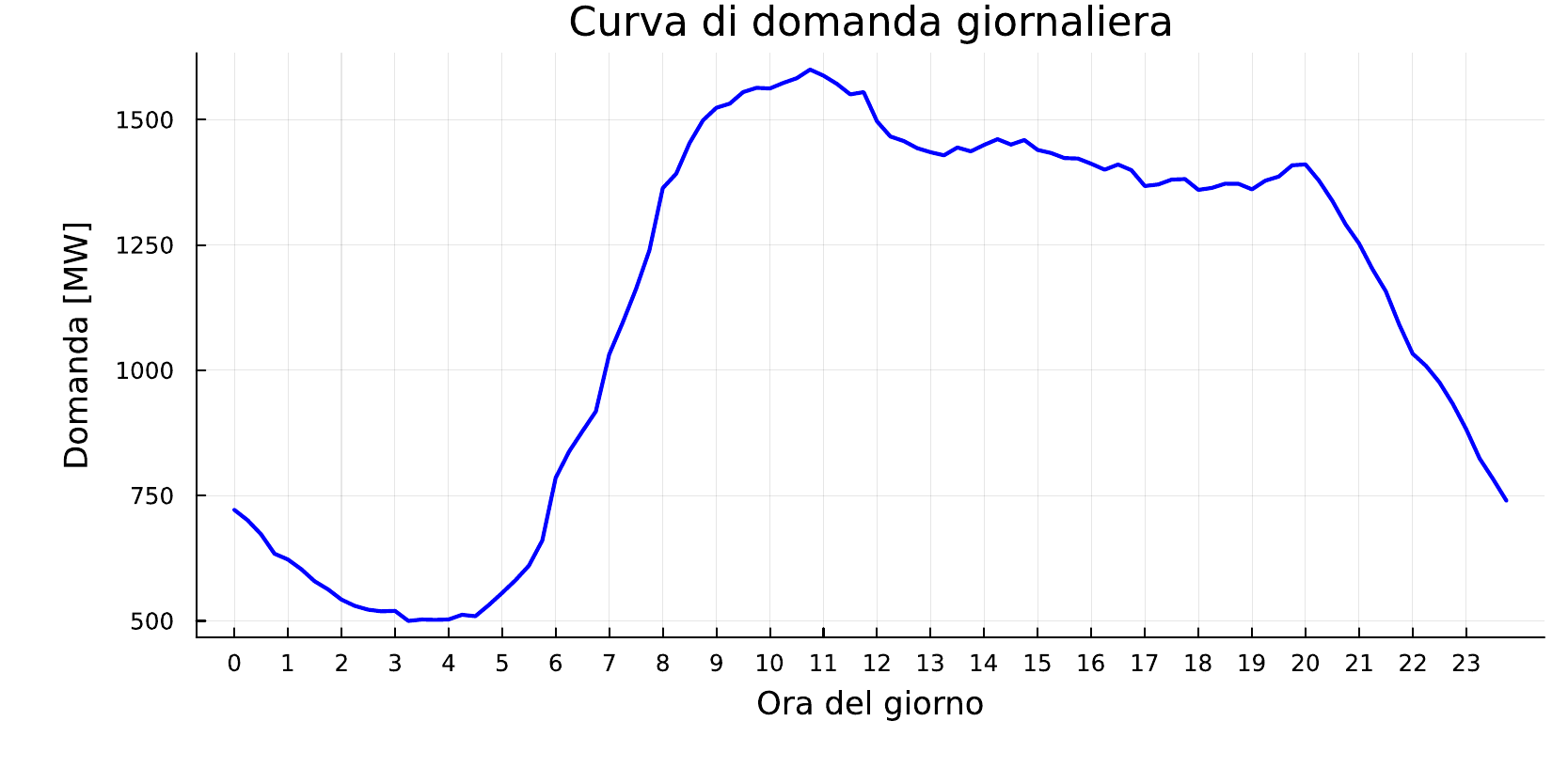}
    \caption{Load curve of September 12, 2025 scaled between 25\% and 80\% of market capacity.}
    \label{fig:curva-carico}
\end{figure}

Once the optimal policy for each market and operator is extracted, simulation of a typical day in the MGP is performed using the load curve of September 12, 2025 extracted from the Terna website (Figure~\ref{fig:curva-carico}). The curve is scaled between $y_{\min} = 25$\% and $y_{\max} = 80$\% of total market capacity (2~GW) via the formula
\[
y' = \frac{y - y_{\min}}{y_{\max} - y_{\min}} \cdot (y'_{\max} - y'_{\min}) + y'_{\min}
\;.
\]

\begin{figure}[htb]
    \centering
    \includegraphics[width=0.8\textwidth]{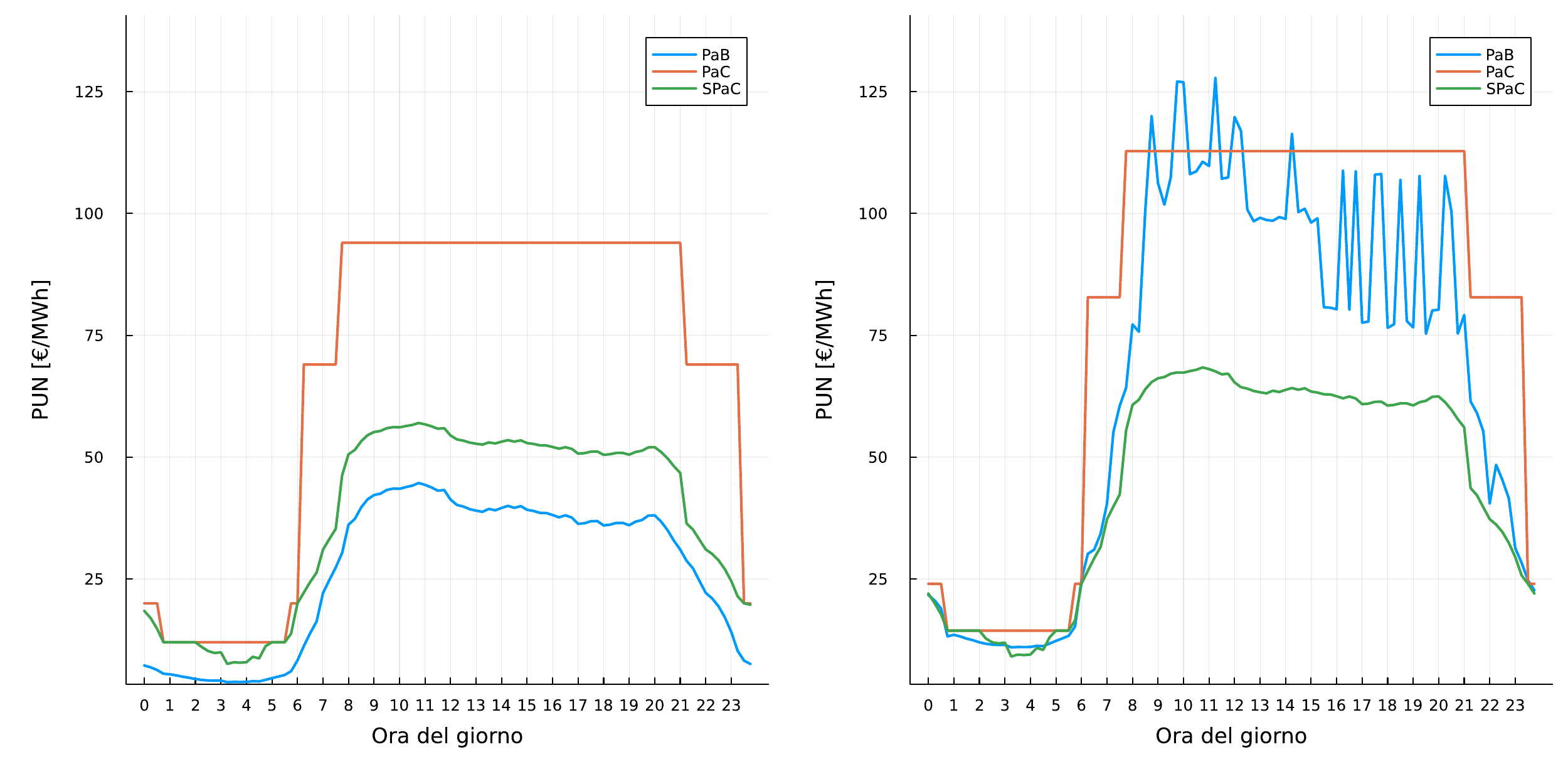}
    \caption{PUN evolution in a typical day, comparison between marginal strategies (left) and optimal policy (right). PNIEC 2030 scenario.}
    \label{fig:pun-giornaliero-rl-vs-mc}
\end{figure}

\noindent
For each iteration, i.e., for each quarter hour of the day, the corresponding demand level is selected from the load curve (Figure~\ref{fig:curva-carico}) and clearing of the three markets is executed using offers at marginal cost and offers based on the policy learned via RL (or the closest one in case of intermediate demand levels not present in the Q-table). The resulting PUN for each market is shown in Figure~\ref{fig:pun-giornaliero-rl-vs-mc}.

\noindent
In the marginal cost case, as already seen, the PaB market attains the lowest PUN which perfectly follows demand trend as expected. The higher PUN of the PaC market definitely stands out during peak hours, when gas units become marginal. The SPaC market attains an intermediate PUN, lower than PaC and higher than PaB. In the case with strategies learned via RL, the PaB market shows a much higher PUN compared to the marginal cost case, since operators try to maximize their profits by applying high markups. This is also reflected in the PaC market, but to a lesser extent in the SPaC market which manages to limit the PUN increase. 

\paragraph{Total Costs}

Analyzing in detail total costs for various demand levels (Table~\ref{tab:pniec-costi-mercati-domanda}), distinct patterns are observed between the two considered offer strategies. In the case of offers at marginal cost, total costs follow an increasing trend consistent with technological merit order.

\begin{table}[htb]
    \caption{PNIEC scenario: total costs in the three markets for different demand levels.}
    \label{tab:pniec-costi-mercati-domanda}
    \resizebox{0.95\textwidth}{!}{
        \begin{tabular}{@{}lrrrrrrrrrrrr@{}}
        \toprule
        Market   Capacity {[}MW{]} &
        \multicolumn{1}{l}{2000} &
        \multicolumn{1}{l}{2000} &
        \multicolumn{1}{l}{2000} &
        \multicolumn{1}{l}{2000} &
        \multicolumn{1}{l}{2000} &
        \multicolumn{1}{l}{2000} &
        \multicolumn{1}{l}{2000} &
        \multicolumn{1}{l}{2000} &
        \multicolumn{1}{l}{2000} &
        \multicolumn{1}{l}{2000} &
        \multicolumn{1}{l}{2000} &
        \multicolumn{1}{l}{2000} \\
        D {[}MW{]} &
        \multicolumn{1}{l}{581} &
        \multicolumn{1}{l}{660} &
        \multicolumn{1}{l}{838} &
        \multicolumn{1}{l}{918} &
        \multicolumn{1}{l}{1095} &
        \multicolumn{1}{l}{1239} &
        \multicolumn{1}{l}{1392} &
        \multicolumn{1}{l}{1499} &
        \multicolumn{1}{l}{1532} &
        \multicolumn{1}{l}{1564} &
        \multicolumn{1}{l}{1573} &
        \multicolumn{1}{l}{1600} \\
        D/Market Capacity {[}\%{]} &
        \multicolumn{1}{l}{29} &
        \multicolumn{1}{l}{33} &
        \multicolumn{1}{l}{42} &
        \multicolumn{1}{l}{46} &
        \multicolumn{1}{l}{55} &
        \multicolumn{1}{l}{62} &
        \multicolumn{1}{l}{70} &
        \multicolumn{1}{l}{75} &
        \multicolumn{1}{l}{77} &
        \multicolumn{1}{l}{78} &
        \multicolumn{1}{l}{79} &
        \multicolumn{1}{l}{80} \\ \midrule
        \multicolumn{13}{c}{\textbf{Strategy: marginal-cost offers}} \\
        Cost tot PaB   {[}\euro{]} &
        2,883.08 &
        3,996.16 &
        9,398.75 &
        14,934.90 &
        27,131.85 &
        37,599.89 &
        51,954.56 &
        61,961.56 &
        65,105.00 &
        68,067.53 &
        68,989.08 &
        71,488.00 \\
        Cost tot PaC {[}\euro{]} &
        6,975.08 &
        13,208.16 &
        57,810.75 &
        63,346.90 &
        75,543.85 &
        116,511.89 &
        130,866.56 &
        140,873.56 &
        144,017.00 &
        146,979.53 &
        147,901.08 &
        150,400.00 \\
        Cost tot SPaC {[}\euro{]} &
        6,975.08 &
        9,088.14 &
        18,610.75 &
        24,146.90 &
        36,343.85 &
        57,311.89 &
        71,666.56 &
        81,673.56 &
        84,817.00 &
        87,779.53 &
        88,701.08 &
        91,200.00 \\ \midrule
        \multicolumn{13}{c}{\textbf{Strategy: RL policy}} \\
        Cost tot PaB   {[}\euro{]} &
        7,434.16 &
        10,068.47 &
        25,261.50 &
        31,462.67 &
        60,307.70 &
        79,623.78 &
        105,435.13 &
        179,848.69 &
        155,961.83 &
        198,682.59 &
        170,021.05 &
        176,988.00 \\
        Cost tot PaC {[}\euro{]} &
        8,370.09 &
        15,849.79 &
        69,372.90 &
        76,016.28 &
        90,652.62 &
        139,814.27 &
        157,039.88 &
        169,048.28 &
        172,820.41 &
        176,375.43 &
        177,481.29 &
        180,480.00 \\
        Cost tot SPaC {[}\euro{]} &
        8,370.09 &
        10,905.77 &
        22,332.90 &
        28,976.28 &
        43,612.62 &
        68,774.27 &
        85,999.88 &
        98,008.28 &
        101,780.41 &
        105,335.43 &
        106,441.29 &
        109,440.00 \\
        Cost   PaB/Cost SPaC &
        \multicolumn{1}{l}{\cellcolor[HTML]{63BE7B}0.9} &
        \multicolumn{1}{l}{\cellcolor[HTML]{70C17B}0.9} &
        \multicolumn{1}{l}{\cellcolor[HTML]{BED880}1.1} &
        \multicolumn{1}{l}{\cellcolor[HTML]{ADD37F}1.1} &
        \multicolumn{1}{l}{\cellcolor[HTML]{FFDA81}1.4} &
        \multicolumn{1}{l}{\cellcolor[HTML]{C8DB80}1.2} &
        \multicolumn{1}{l}{\cellcolor[HTML]{E1E282}1.2} &
        \multicolumn{1}{l}{\cellcolor[HTML]{F9756E}1.8} &
        \multicolumn{1}{l}{\cellcolor[HTML]{FDB97B}1.5} &
        \multicolumn{1}{l}{\cellcolor[HTML]{F8696B}1.9} &
        \multicolumn{1}{l}{\cellcolor[HTML]{FCAA78}1.6} &
        \multicolumn{1}{l}{\cellcolor[HTML]{FCA677}1.6} \\
        Cost PaC/Cost SPaC &
        \multicolumn{1}{l}{\cellcolor[HTML]{63BE7B}1.0} &
        \multicolumn{1}{l}{\cellcolor[HTML]{C6DA80}1.5} &
        \multicolumn{1}{l}{\cellcolor[HTML]{F8696B}3.1} &
        \multicolumn{1}{l}{\cellcolor[HTML]{FB9774}2.6} &
        \multicolumn{1}{l}{\cellcolor[HTML]{FEC97E}2.1} &
        \multicolumn{1}{l}{\cellcolor[HTML]{FECE7F}2.0} &
        \multicolumn{1}{l}{\cellcolor[HTML]{FFE182}1.8} &
        \multicolumn{1}{l}{\cellcolor[HTML]{FFEA84}1.7} &
        \multicolumn{1}{l}{\cellcolor[HTML]{FCEA83}1.7} &
        \multicolumn{1}{l}{\cellcolor[HTML]{F6E883}1.7} &
        \multicolumn{1}{l}{\cellcolor[HTML]{F5E883}1.7} &
        \multicolumn{1}{l}{\cellcolor[HTML]{F1E783}1.6} \\ \bottomrule
        \end{tabular}
    }
\end{table}

\noindent
The PaB market presents the lowest costs, varying from 2,883~\EUR{} (demand 29\%) to 71,488~\EUR{} (demand 80\%). The PaC market shows significantly higher costs, from 6,975~\EUR{} to 150,400~\EUR{}, highlighting the effect of inframarginal rents. The SPaC market positions itself intermediately, with costs from 6,975~\EUR{} to 91,200~\EUR{}.

\smallskip
\noindent
Introducing strategies learned via RL involves generalized cost increases in all three regimes. However, the extent of such increases varies considerably: PaB registers the most dramatic increases (from +158\% to +148\% depending on demand level), reflecting operators' ability to fully exploit the pay-as-bid mechanism. PaC shows more contained increases (+20-30\%)---but starting from a much higher value---, while SPaC presents intermediate increases (+20-50\%). Comparison of the last two rows of the table highlights SPaC superiority: the PaB/SPaC ratio varies from 0.9 to 1.9, while the PaC/SPaC ratio oscillates between 1.0 and 3.1. This confirms that SPaC maintains a substantial advantage in terms of economic efficiency, especially at low demand levels where aggressive behavior in NMCS units leads the bilevel problem to consider offers in the NNMCS market to reduce costs.

\paragraph{Operator Profits}

Table~\ref{tab:pniec-profitti-domanda} reveals interesting dynamics in profit distribution among different market regimes. In the marginal cost case, PaB generates no profits for operators (offers at actual cost), while PaC and SPaC produce profits of 78,912~\EUR{} and 19,712~\EUR{} respectively at highest demand levels.

\begin{table}[htb]
    \caption{PNIEC scenario: operator profits in the three markets for different demand levels.}
    \label{tab:pniec-profitti-domanda}
    \resizebox{0.95\textwidth}{!}{
        \begin{tabular}{@{}lrrrrrrrrrrrr@{}}
        \toprule
        Market Capacity {[}MW{]} &
        \multicolumn{1}{l}{2000} &
        \multicolumn{1}{l}{2000} &
        \multicolumn{1}{l}{2000} &
        \multicolumn{1}{l}{2000} &
        \multicolumn{1}{l}{2000} &
        \multicolumn{1}{l}{2000} &
        \multicolumn{1}{l}{2000} &
        \multicolumn{1}{l}{2000} &
        \multicolumn{1}{l}{2000} &
        \multicolumn{1}{l}{2000} &
        \multicolumn{1}{l}{2000} &
        \multicolumn{1}{l}{2000} \\
        D {[}MW{]} &
        \multicolumn{1}{l}{581} &
        \multicolumn{1}{l}{660} &
        \multicolumn{1}{l}{838} &
        \multicolumn{1}{l}{918} &
        \multicolumn{1}{l}{1095} &
        \multicolumn{1}{l}{1239} &
        \multicolumn{1}{l}{1392} &
        \multicolumn{1}{l}{1499} &
        \multicolumn{1}{l}{1532} &
        \multicolumn{1}{l}{1564} &
        \multicolumn{1}{l}{1573} &
        \multicolumn{1}{l}{1600} \\
        D/Market Capacity {[}\%{]} &
        \multicolumn{1}{l}{29} &
        \multicolumn{1}{l}{33} &
        \multicolumn{1}{l}{42} &
        \multicolumn{1}{l}{46} &
        \multicolumn{1}{l}{55} &
        \multicolumn{1}{l}{62} &
        \multicolumn{1}{l}{70} &
        \multicolumn{1}{l}{75} &
        \multicolumn{1}{l}{77} &
        \multicolumn{1}{l}{78} &
        \multicolumn{1}{l}{79} &
        \multicolumn{1}{l}{80} \\ \midrule
        \multicolumn{13}{c}{\textbf{Strategy: marginal-cost offers}} \\
        Operator profits PaB {[}\euro{]} &
        - &
        - &
        - &
        - &
        - &
        - &
        - &
        - &
        - &
        - &
        - &
        - \\
        Operator profits PaC   {[}\euro{]} &
        4,092.00 &
        9,212.00 &
        48,412.00 &
        48,412.00 &
        48,412.00 &
        78,912.00 &
        78,912.00 &
        78,912.00 &
        78,912.00 &
        78,912.00 &
        78,912.00 &
        78,912.00 \\
        Operator profits SPaC   {[}\euro{]} &
        4,092.00 &
        4,092.00 &
        9,212.00 &
        9,212.00 &
        9,212.00 &
        19,712.00 &
        19,712.00 &
        19,712.00 &
        19,712.00 &
        19,712.00 &
        19,712.00 &
        19,712.00 \\ \midrule
        \multicolumn{13}{c}{\textbf{Strategy: RL policy}} \\
        Operator profits PaB {[}\euro{]} &
        4,551.08 &
        6,072.32 &
        15,862.75 &
        13,576.00 &
        33,175.85 &
        42,023.89 &
        53,480.56 &
        117,887.13 &
        88,159.22 &
        122,515.06 &
        99,367.37 &
        104,500.00 \\
        Operator profits PaC   {[}\euro{]} &
        5,487.02 &
        11,853.63 &
        59,974.15 &
        61,081.38 &
        63,520.77 &
        102,214.38 &
        105,085.31 &
        107,086.71 &
        107,715.40 &
        108,307.91 &
        108,492.22 &
        108,992.00 \\
        Operator profits SPaC   {[}\euro{]} &
        5,487.02 &
        5,909.63 &
        12,934.15 &
        14,041.38 &
        16,480.77 &
        31,174.38 &
        34,045.31 &
        36,046.71 &
        36,675.40 &
        37,267.91 &
        37,452.22 &
        37,952.00 \\
        Profits PaB/Profits SPaC &
        \multicolumn{1}{l}{\cellcolor[HTML]{63BE7B}0.8} &
        \multicolumn{1}{l}{\cellcolor[HTML]{83C77C}1.0} &
        \multicolumn{1}{l}{\cellcolor[HTML]{A3D07E}1.2} &
        \multicolumn{1}{l}{\cellcolor[HTML]{79C47C}1.0} &
        \multicolumn{1}{l}{\cellcolor[HTML]{FED881}2.0} &
        \multicolumn{1}{l}{\cellcolor[HTML]{B7D67F}1.3} &
        \multicolumn{1}{l}{\cellcolor[HTML]{DBE081}1.6} &
        \multicolumn{1}{l}{\cellcolor[HTML]{F96B6C}3.3} &
        \multicolumn{1}{l}{\cellcolor[HTML]{FDB67A}2.4} &
        \multicolumn{1}{l}{\cellcolor[HTML]{F8696B}3.3} &
        \multicolumn{1}{l}{\cellcolor[HTML]{FBA176}2.7} &
        \multicolumn{1}{l}{\cellcolor[HTML]{FB9874}2.8} \\
        Profits PaC/Profits SPaC &
        \multicolumn{1}{l}{\cellcolor[HTML]{63BE7B}1.0} &
        \multicolumn{1}{l}{\cellcolor[HTML]{B3D57F}2.0} &
        \multicolumn{1}{l}{\cellcolor[HTML]{F8696B}4.6} &
        \multicolumn{1}{l}{\cellcolor[HTML]{FA8070}4.4} &
        \multicolumn{1}{l}{\cellcolor[HTML]{FCA677}3.9} &
        \multicolumn{1}{l}{\cellcolor[HTML]{FED280}3.3} &
        \multicolumn{1}{l}{\cellcolor[HTML]{FFE183}3.1} &
        \multicolumn{1}{l}{\cellcolor[HTML]{FFEA84}3.0} &
        \multicolumn{1}{l}{\cellcolor[HTML]{FDEA83}2.9} &
        \multicolumn{1}{l}{\cellcolor[HTML]{FBE983}2.9} &
        \multicolumn{1}{l}{\cellcolor[HTML]{FAE983}2.9} &
        \multicolumn{1}{l}{\cellcolor[HTML]{F8E983}2.9} \\ \bottomrule
        \end{tabular}
    }
\end{table}

\noindent
With adoption of RL strategies, all three markets generate significant profits for operators. PaB becomes the most profitable, with profits reaching 122,515~\EUR{} at demand peaks. PaC maintains high profits (up to 108,992~\EUR{}), while SPaC generates lower but still substantial profits (up to 37,952~\EUR{}). It is particularly significant to observe that in the SPaC market operators still receive significant compensation for technologies in the NMCS sub-market (renewables and hydroelectric). This happens because, despite segmentation limiting inframarginal rents, NMCS technologies still benefit from prices above marginal costs when demand in their segment is high, guaranteeing economic sustainability to investments in clean technologies. Ratios in the last two rows confirm the balance of the SPaC mechanism: PaB/SPaC profit ratios vary from 0.8 to 3.3, while PaC/SPaC profit ratios oscillate between 1.0 and 4.6.

\paragraph{PUN Evolution}

Table~\ref{tab:pniec-pun-mercati-domanda} shows PUN evolution in the three market regimes. With marginal cost strategies, PaB PUN grows gradually from 5.0~\EUR{}/MWh (demand 29\%) to 44.7~\EUR{}/MWh (demand 80\%), faithfully following merit order based on marginal costs. PaC presents a step-like trend, with pronounced jumps when more expensive technologies activate: from 12.0~\EUR{}/MWh to 94.0~\EUR{}/MWh. SPaC shows a more gradual evolution, from 12.0~\EUR{}/MWh to 57.0~\EUR{}/MWh, benefiting from segmentation.

\smallskip
\noindent
Adoption of RL strategies involves significant PUN increases in all regimes. PaB registers the largest increase in relative terms, with prices reaching 127.1~\EUR{}/MWh (+184\% compared to marginal case). PaC presents more contained but still substantial increases, reaching 112.8~\EUR{}/MWh (+20\%). SPaC shows the most limited increase, with a maximum of 68.4~\EUR{}/MWh (+20\%). From the consumer point of view, SPaC clearly is the most convenient regime both in absolute and relative terms. Even with optimized strategies, it maintains significantly lower prices compared to the other two regimes, offering the best compromise between economic efficiency and sustainability of investments in the electricity system. 

\begin{table}[htb]
    \caption{PNIEC scenario: PUN in the three markets for different demand levels.}
    \label{tab:pniec-pun-mercati-domanda}
    \resizebox{0.95\textwidth}{!}{
        \begin{tabular}{@{}lrrrrrrrrrrrr@{}}
        \toprule
        Market   Capacity {[}MW{]} &
        \multicolumn{1}{l}{2000} &
        \multicolumn{1}{l}{2000} &
        \multicolumn{1}{l}{2000} &
        \multicolumn{1}{l}{2000} &
        \multicolumn{1}{l}{2000} &
        \multicolumn{1}{l}{2000} &
        \multicolumn{1}{l}{2000} &
        \multicolumn{1}{l}{2000} &
        \multicolumn{1}{l}{2000} &
        \multicolumn{1}{l}{2000} &
        \multicolumn{1}{l}{2000} &
        \multicolumn{1}{l}{2000} \\
        D {[}MW{]} &
        \multicolumn{1}{l}{581} &
        \multicolumn{1}{l}{660} &
        \multicolumn{1}{l}{838} &
        \multicolumn{1}{l}{918} &
        \multicolumn{1}{l}{1095} &
        \multicolumn{1}{l}{1239} &
        \multicolumn{1}{l}{1392} &
        \multicolumn{1}{l}{1499} &
        \multicolumn{1}{l}{1532} &
        \multicolumn{1}{l}{1564} &
        \multicolumn{1}{l}{1573} &
        \multicolumn{1}{l}{1600} \\
        D/Market Capacity {[}\%{]} &
        \multicolumn{1}{l}{29} &
        \multicolumn{1}{l}{33} &
        \multicolumn{1}{l}{42} &
        \multicolumn{1}{l}{46} &
        \multicolumn{1}{l}{55} &
        \multicolumn{1}{l}{62} &
        \multicolumn{1}{l}{70} &
        \multicolumn{1}{l}{75} &
        \multicolumn{1}{l}{77} &
        \multicolumn{1}{l}{78} &
        \multicolumn{1}{l}{79} &
        \multicolumn{1}{l}{80} \\ \midrule
        \multicolumn{13}{c}{\textbf{Strategy: marginal-cost offers}} \\
        PUN PaB   {[}\euro/MWh{]} &
        \cellcolor[HTML]{FCFCFF}5.00 &
        \cellcolor[HTML]{FCFBFE}6.10 &
        \cellcolor[HTML]{FCF2F5}11.20 &
        \cellcolor[HTML]{FCEAED}16.30 &
        \cellcolor[HTML]{FCDCDF}24.80 &
        \cellcolor[HTML]{FBD3D5}30.30 &
        \cellcolor[HTML]{FBC7CA}37.30 &
        \cellcolor[HTML]{FBC1C3}41.30 &
        \cellcolor[HTML]{FBBFC1}42.50 &
        \cellcolor[HTML]{FBBDBF}43.50 &
        \cellcolor[HTML]{FBBCBF}43.80 &
        \cellcolor[HTML]{FBBBBD}44.70 \\
        PUN PaC {[}\euro/MWh{]} &
        \cellcolor[HTML]{FCF1F4}12.00 &
        \cellcolor[HTML]{FCE4E7}20.00 &
        \cellcolor[HTML]{FA9395}69.00 &
        \cellcolor[HTML]{FA9395}69.00 &
        \cellcolor[HTML]{FA9395}69.00 &
        \cellcolor[HTML]{F8696B}94.00 &
        \cellcolor[HTML]{F8696B}94.00 &
        \cellcolor[HTML]{F8696B}94.00 &
        \cellcolor[HTML]{F8696B}94.00 &
        \cellcolor[HTML]{F8696B}94.00 &
        \cellcolor[HTML]{F8696B}94.00 &
        \cellcolor[HTML]{F8696B}94.00 \\
        PUN SPaC {[}\euro/MWh{]} &
        \cellcolor[HTML]{FCF1F4}12.00 &
        \cellcolor[HTML]{FCEEF1}13.80 &
        \cellcolor[HTML]{FCE0E3}22.20 &
        \cellcolor[HTML]{FCD9DC}26.30 &
        \cellcolor[HTML]{FBCED1}33.20 &
        \cellcolor[HTML]{FBB8BB}46.20 &
        \cellcolor[HTML]{FAB0B2}51.50 &
        \cellcolor[HTML]{FAABAD}54.50 &
        \cellcolor[HTML]{FAA9AC}55.40 &
        \cellcolor[HTML]{FAA8AB}56.10 &
        \cellcolor[HTML]{FAA8AA}56.40 &
        \cellcolor[HTML]{FAA7A9}57.00 \\ \midrule
        \multicolumn{13}{c}{\textbf{Strategy: RL policy}} \\
        PUN PaB   {[}\euro/MWh{]} &
        \cellcolor[HTML]{FCFCFF}12.80 &
        \cellcolor[HTML]{FCF9FC}15.20 &
        \cellcolor[HTML]{FCE6E9}30.20 &
        \cellcolor[HTML]{FCE1E4}34.30 &
        \cellcolor[HTML]{FBC6C9}55.10 &
        \cellcolor[HTML]{FBBABD}64.20 &
        \cellcolor[HTML]{FAACAE}75.70 &
        \cellcolor[HTML]{F97375}120.00 &
        \cellcolor[HTML]{F98A8C}101.80 &
        \cellcolor[HTML]{F8696B}127.10 &
        \cellcolor[HTML]{F98284}108.10 &
        \cellcolor[HTML]{F97F81}110.60 \\
        PUN PaC {[}\euro/MWh{]} &
        \cellcolor[HTML]{FCFAFD}14.40 &
        \cellcolor[HTML]{FCEEF1}24.00 &
        \cellcolor[HTML]{FAA2A5}82.80 &
        \cellcolor[HTML]{FAA2A5}82.80 &
        \cellcolor[HTML]{FAA2A5}82.80 &
        \cellcolor[HTML]{F97C7E}112.80 &
        \cellcolor[HTML]{F97C7E}112.80 &
        \cellcolor[HTML]{F97C7E}112.80 &
        \cellcolor[HTML]{F97C7E}112.80 &
        \cellcolor[HTML]{F97C7E}112.80 &
        \cellcolor[HTML]{F97C7E}112.80 &
        \cellcolor[HTML]{F97C7E}112.80 \\
        PUN SPaC {[}\euro/MWh{]} &
        \cellcolor[HTML]{FCFAFD}14.40 &
        \cellcolor[HTML]{FCF8FB}16.50 &
        \cellcolor[HTML]{FCEBEE}26.70 &
        \cellcolor[HTML]{FCE4E7}31.60 &
        \cellcolor[HTML]{FCDADD}39.80 &
        \cellcolor[HTML]{FBC6C8}55.50 &
        \cellcolor[HTML]{FBBDC0}61.80 &
        \cellcolor[HTML]{FBB9BB}65.40 &
        \cellcolor[HTML]{FBB8BA}66.40 &
        \cellcolor[HTML]{FBB6B9}67.40 &
        \cellcolor[HTML]{FBB6B9}67.60 &
        \cellcolor[HTML]{FBB5B8}68.40 \\ \bottomrule
        \end{tabular}
    }
\end{table}

\subsubsection{Synthetic Results}

Table~\ref{tab:ris-sintetici-PNIEC} summarizes overall results of the simulation of a typical day in the MGP for the PNIEC 2030 scenario.

\begin{table}[H]
    \caption{Synthetic results, simulation of a typical MGP day, PNIEC scenario}
    \label{tab:ris-sintetici-PNIEC}
    \centering
    \resizebox{0.95\textwidth}{!}{%
    \begin{tabular}{@{}lllllll@{}}
        \toprule
        & \multicolumn{3}{c}{\textbf{Strategy: marginal   costs}}           & \multicolumn{3}{c}{\textbf{Strategy: RL policy}}                   \\ \midrule
        \textbf{Market} & Market cost {[}\EUR{}{]} \footnotemark & Operator profits {[}\EUR{}{]} & Profit/Cost & Market cost {[}\EUR{}{]} & Operator profits {[}\EUR{}{]} & Profit/Cost \\
        PaB  & 3,435,303.40 & -            & 0\%  & 8,543,866.32  & 4,979,805.83 & 58\% \\
        PaC  & 8,569,055.40 & 5,133,752.00 & 60\% & 10,282,866.49 & 6,847,563.08 & 67\% \\
        SPaC & 4,757,393.73 & 1,286,392.00 & 27\% & 5,706,729.83  & 2,250,017.18 & 39\% \\ \bottomrule
    \end{tabular}%
    }
\end{table}
\footnotetext{By \emph{Market cost} we mean the total expenditure sustained by consumers for purchasing energy in the market which, given demand inelasticity, coincides with aggregated producer revenues.}

\noindent
This data confirm the distinctive characteristics of the three analyzed market regimes discussed throughout previous simulations. The \emph{PaB} market is the most efficient in terms of total cost when operators offer at marginal cost (3.43~M\EUR{}), but becomes the least efficient with optimized strategies (8.54~M\EUR{}), registering a 149\% increase. The \emph{PaC} market maintains high costs in both cases (from 8.57~M\EUR{} to 10.28~M\EUR{}, +20\%), highlighting relative stability due to the marginal pricing mechanism. The \emph{SPaC} market confirms itself as the most balanced intermediate solution, with costs varying from 4.76~M\EUR{} to 5.71~M\EUR{} (+20\%).

\smallskip
\noindent
This is confirmed by the Profit/Cost ratio, which measures what percentage of total market cost is transferred to operators as profit, while the remaining share represents the actual production cost of energy. In the case of marginal cost strategies, PaB presents a Profit/Cost ratio of 0\%, indicating that the entire amount paid by consumers covers exclusively production costs, without generating rents for operators. PaC shows a ratio of 60\%, highlighting how a majority share of expenditure goes to constitute inframarginal rents for producers. SPaC reaches an intermediate equilibrium with 27\%, significantly limiting wealth transfer towards producers compared to traditional PaC.

\smallskip
\noindent
With introduction of RL strategies, all regimes see increased profit component. PaB reaches 58\%, demonstrating how operators manage to extract significant value from the pay-as-bid mechanism when adopting strategic behaviors. PaC increases to 67\%, confirming its intrinsic nature favorable to producers. SPaC rises to 39\%, still maintaining the most favorable ratio for consumers. Profit distribution among operators (Table~\ref{tab:dettaglio-profitti-PNIEC}) also reveals an interesting asymmetry: operator OpB results systematically more profitable compared to OpA in all market regimes, due to the lower marginal costs of its production units (Table~\ref{tab:portafoglio-operatori-pniec}).


\begin{table}[H]
    \caption{Operator profit breakdown: simulation of a representative day of the MGP, PNIEC scenario.}
    \label{tab:dettaglio-profitti-PNIEC}
    \centering
    \resizebox{0.5\textwidth}{!}{%
    \begin{tabular}{@{}llll@{}}
        \toprule
        \multicolumn{1}{r}{\multirow{2}{*}{Operator}} & \multicolumn{3}{c}{Profits {[}\euro{]}}                                         \\
        \multicolumn{1}{r}{}                           & \multicolumn{1}{c}{PaB} & \multicolumn{1}{c}{PaC} & \multicolumn{1}{c}{SPaC} \\ \midrule
        \multicolumn{4}{c}{\textbf{Strategy: marginal-cost offers}} \\
        OpB    & -              & 2,827,656.00   & 907,176.00     \\
        OpA    & -              & 2,306,096.00   & 379,216.00     \\ \midrule
        \multicolumn{4}{c}{\textbf{Strategy: RL policy}}         \\
        OpB    & 2,912,542.71   & 3,805,330.54   & 1,511,414.13   \\
        OpA    & 2,067,263.12   & 3,042,232.54   & 738,603.05     \\ \bottomrule
    \end{tabular}%
    }
\end{table}

\subsection{Scenario 2: 10 operators from GME public offers}


\subsubsection{Clustering of Public Offers (2024)}

\begin{figure}[htb]
    \centering
    \includegraphics[width=0.5\textwidth]{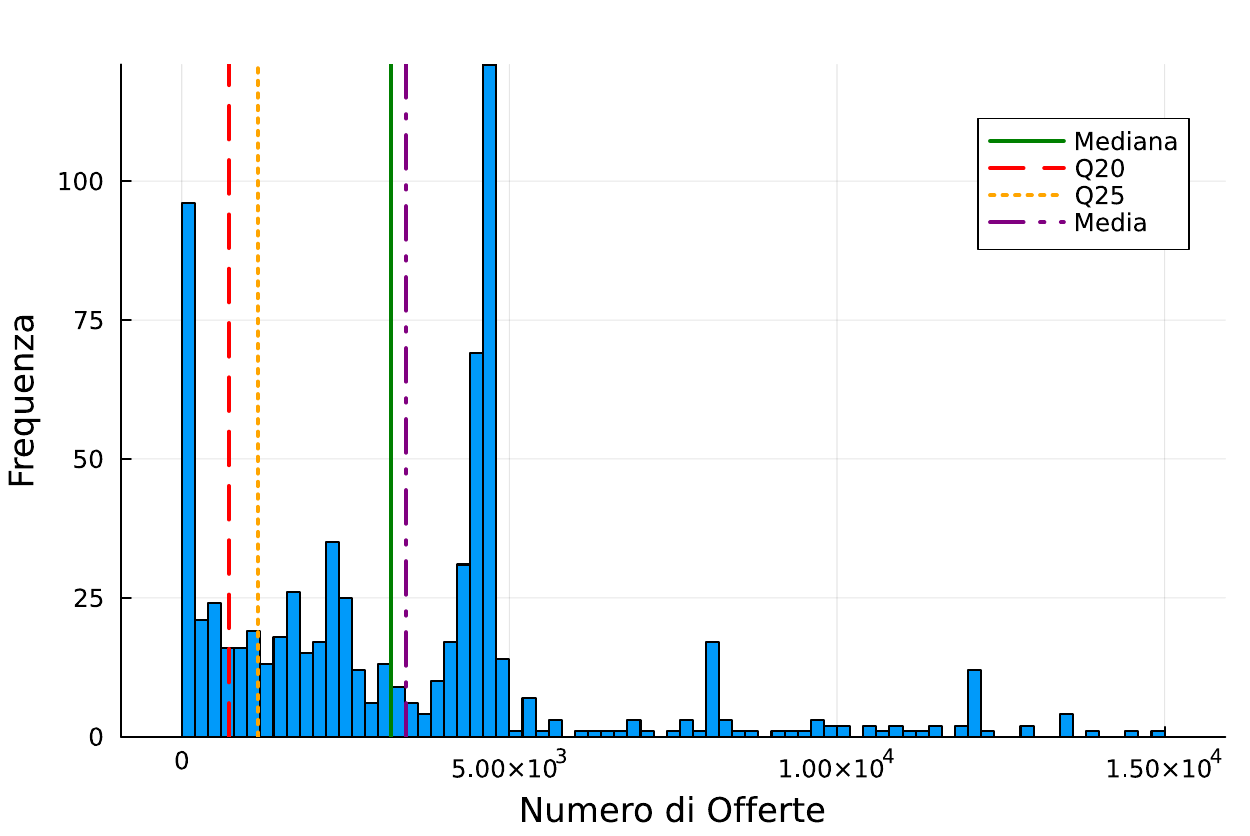}
    \caption{Distribution of number of public offers per production unit (PU)}
    \label{fig:offerte-per-up-distribution}
\end{figure}

Starting from public offers of the MGP of the year 2024, available on the GME website \cite{GestoreMercatiEnergetici}, portfolios for 10 operators to be used in simulations were extracted. The objective is to group production units that present offers into homogeneous technological clusters, allowing assignment to each operator of a structured portfolio on which to apply differentiated markup strategies per cluster via RL algorithms. The approach of inferring technological characteristics of plants through analysis of offer patterns is not new in literature, e.g., using clustering techniques on price-quantity pairs for automatic plant classification \cite{durvasuluClassificationGeneratorsParticipating2017,linAutomatedClassificationPower2020}. Our methodology involves the following data preparation steps:
\begin{enumerate}
 \item all public offers of 2024 were downloaded and filtered to keep exclusively those validated by GME controls, including only offers accepted and rejected by market outcomes;
 \item a filter was applied to extract only offers of main operators of the Italian electricity market according to the Annex of the ARERA Decision \cite{areraDelibera30220252025};
 \item no geographical distinctions were made between offers to simplify the analysis.
\end{enumerate}
From the dataset thus obtained, offers of 10 operators were extracted, for a total of 743 distinct UPs (Table~\ref{tab:numero-up-top10}). Furthermore, to avoid categorizing UPs with few offers, all those with less than the 20th percentile of number of offers were excluded (720 offers, Figure \ref{fig:offerte-per-up-distribution}) thus yielding a final set of 595 distinct UPs.

\begin{table}[H]
    \centering
    \caption{Number of PUs in the 2024 public offers for 10 operators.}
    \label{tab:numero-up-top10}
    \resizebox{0.5\textwidth}{!}{
    \begin{tabular}{l r r}
    \toprule
    Operator & Number of PUs & Number of PUs after filtering \\
    \midrule
    OP1     & 260 & 143 \\
    OP3     & 124 & 119 \\
    OP2     & 75 & 64  \\
    OP8     & 10 & 10  \\
    OP10    & 85 & 79  \\
    OP9     & 74 & 72  \\
    OP5     & 47 & 42  \\
    OP6     & 34 & 31  \\
    OP7     & 17 & 17  \\
    OP4     & 18 & 18  \\
    \midrule
    Total                 & 743 & 595 \\
    \bottomrule
    \end{tabular}
    }
\end{table}


\begin{figure}[htb]
    \centering
    \includegraphics[width=0.5\textwidth]{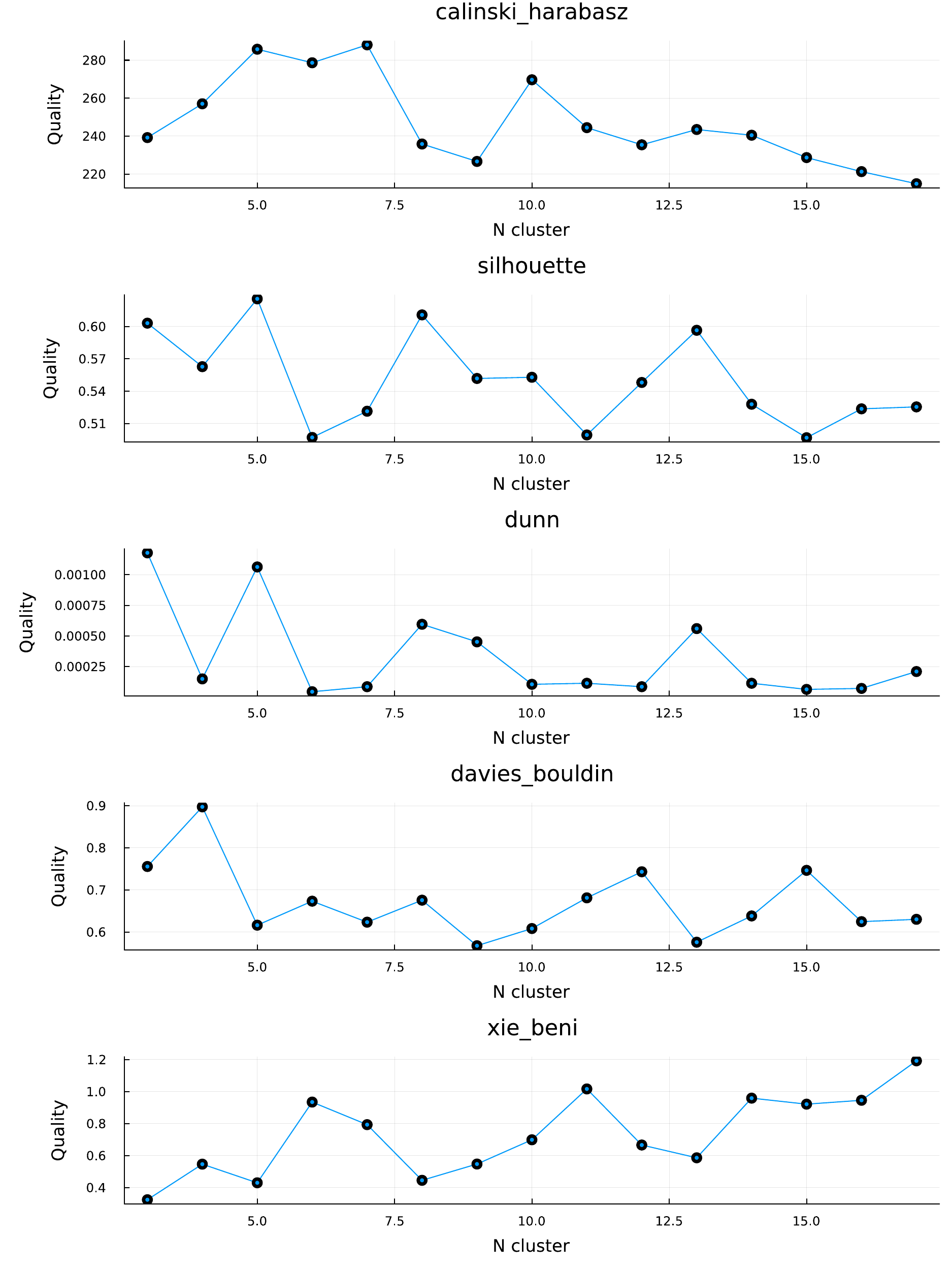}
    \caption{Quality indices for k-means clustering}
    \label{fig:kmeans-quality-indices}
\end{figure}

\noindent
For unsupervised clustering, the k-means algorithm of the \textit{Clustering.jl} Julia package was used \cite{KmeansClusteringjl2025}. Feature choice to characterize bidding behavior of each of the 595 UPs is based on eight indicators capturing different operational and economic aspects:
\begin{itemize}
 \item $P_\mathrm{w}$: weighted average offered price, calculated as average of offered prices weighted by corresponding quantities.
 \item $Q_\mathrm{max}$: maximum offered capacity, representative of plant size.
 \item $CV = \frac{\sigma(Q)}{\mu(Q)}$: coefficient of variation of offered quantity, representative of plant flexibility or intermittency.
 \item $\mathrm{Flex}\ Q = \frac{\max(Q) - \min(Q)}{\max(Q)}$: operational flexibility index, measuring modulation capability between technical minimum and maximum.
 \item $V_\mathrm{h} = \sigma\left(\left\{\overline{Q}_h\right\}_{h=0}^{23}\right)$: hourly variability, calculated as standard deviation of average offered quantity for each hour of the day.
 \item $V_\mathrm{m} = \sigma\left(\left\{\overline{Q}_m\right\}_{m=1}^{12}\right)$: seasonal variability, calculated as standard deviation of average offered quantity for each month.
 \item $h_\mathrm{night} = 100 \cdot \frac{N_\mathrm{night}}{N_\mathrm{tot}}$: percentage of offers in night hours (22:00--6:00), indicative of base load plants.
 \item $h_\mathrm{peak} = 100 \cdot \frac{N_\mathrm{picco}}{N_\mathrm{tot}}$: percentage of offers in peak hours (9:00--12:00 and 18:00--20:00), representative of plants used to follow demand.
\end{itemize}
All features were normalized via z-score to guarantee they have equivalent weight in the clustering algorithm:
\begin{equation}
    z = \frac{x - \mu}{\sigma}
\end{equation}
where $x$ is the feature value, $\mu$ the sample mean and $\sigma$ the standard deviation. Optimal cluster number was determined via evaluation of different clustering quality indices (Figure~\ref{fig:kmeans-quality-indices}).

\noindent
A multi-criteria normalization procedure was adopted: for each $k$ value, five validation indices were calculated (silhouette, calinski-harabasz, dunn, davies-bouldin, xie-beni) and normalized on 0--1 scale, where 1 always represents best performance. For indices where maximum is preferable (silhouette, calinski-harabasz, dunn) min-max normalization was applied: $(x - \min(x)) / (\max(x) - \min(x))$. For those where minimum is preferable (davies-bouldin, xie-beni) the max-min was rather used: $(\max(x) - x) / (\max(x) - \min(x))$. The optimal value $k = 5$ was identified as that maximizing the average of normalized scores, resulting in PU subdivision reported in Table~\ref{tab:clustering-up-statistiche}.

\begin{figure}[htb]
    \centering
    \includegraphics[width=0.5\textwidth]{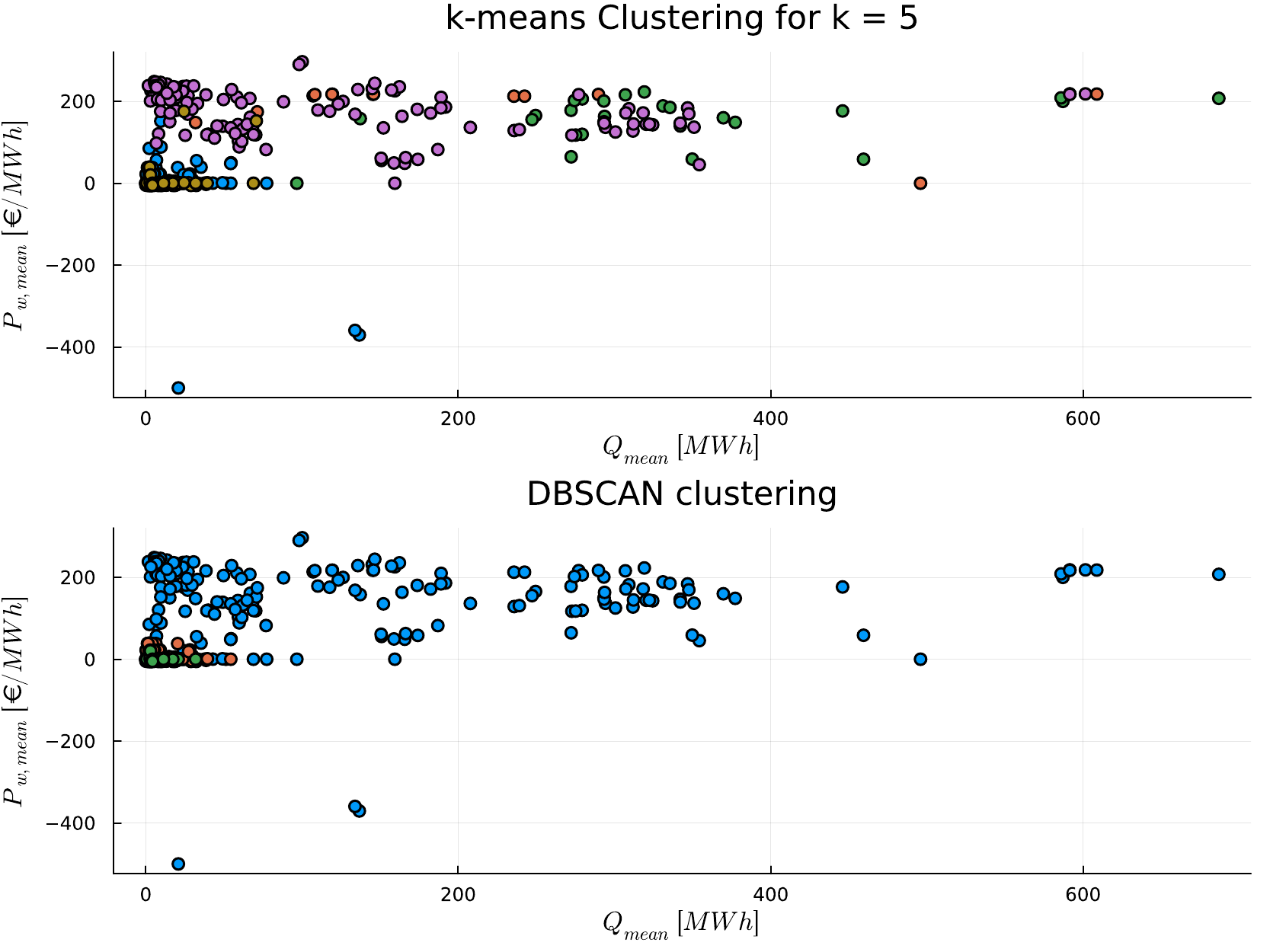}
    \caption{Clustering results with k-means and DBSCAN for comparison, projection on $P_\mathrm{w}$ and $Q_\mathrm{mean}$ of each PU.}
    \label{fig:clustering-results}
\end{figure}

\begin{table}[H]
    \centering
    \caption{Descriptive statistics of the clusters obtained from clustering of PUs (average feature values for each cluster).}
    \label{tab:clustering-up-statistiche}
    \resizebox{0.8\textwidth}{!}{
    \begin{tabular}{cccccccccc}
        \toprule
        Cluster & Number of PUs & $P_\mathrm{w}$ [\euro/MWh] & $Q_\mathrm{max}$ [MW] & CV & Flex $Q$ & $V_\mathrm{h}$ & $V_\mathrm{m}$ & $h_{night} $ & $h_{peak}$ \\
        \midrule
        1 & 329 & 1.94 & 28.99 & 0.80 & 0.98 & 1.51 & 2.96 & 36.34 & 29.88 \\
        2 & 32  & 97.27 & 133.45 & 0.07 & 0.30 & 0.18 & 7.45 & 35.63 & 29.98 \\
        3 & 25  & 156.03 & 744.56 & 0.63 & 0.98 & 65.35 & 100.20 & 32.29 & 34.03 \\
        4 & 123 & 174.05 & 198.61 & 0.43 & 0.97 & 5.22 & 23.25 & 36.75 & 30.39 \\
        5 & 86  & 3.79 & 23.79 & 0.68 & 0.97 & 4.76 & 2.27 & 2.82 & 45.79 \\
        \bottomrule
    \end{tabular}
    }
\end{table}

\noindent
Then, each PU was attributed a capacity as the average offered quantity $Q_\mathrm{mean}$ calculated over all public offers of 2024, and a marginal cost corresponding to the weighted average price $P_\mathrm{w}$. The final result is a structured dataset of 595 UPs classified in 5 technological clusters, each characterized by specific capacity and marginal cost\footnotemark.
\footnotetext{Notably, since marginal costs are estimated from the values derived from public offers of 2024, prices are expected to deviate from those actually observed in the market for the same year. Indeed, the approach is not considering the impact in the market offers of operating and maintenance costs, variable fuel costs for the offer period, costs associated to possible ETS certificates for CO$_2$ emissions and, where present, markups applied by operators. Thus, the final prices computed by the simulations cannot be taken as representative of realistic market prices. However, the analysis in this report just requires reasonably realistic reference marginal costs on which to apply differentiated markup strategies to compare the three market regimes.}

\smallskip
\noindent
Since no official public data on technologies of each PU are available, classification into NMCS and NNMCS was performed via a criterion based exclusively on estimated marginal costs. A threshold of 110~\EUR{}/MWh was adopted to discriminate the two categories: units with marginal cost below are classified as NMCS, those with cost above as NNMCS. This threshold was chosen since it allows separating renewable and hydroelectric technologies (typically characterized by marginal costs below 20~\EUR{}/MWh) from fossil technologies (natural gas and coal, with costs above 100~\EUR{}/MWh), consistent with values reported in literature \cite{StudyLevelizedCost}. Notably, this distinction is necessary for application of the SPaC mechanism, which explicitly requires information on technology type. The effect of this assumption on the validity of results represents a methodological limitation of this study, which could be overcome should official technological information of UPs become available.

\begin{figure}[htb]
    \centering
    \includegraphics[width=0.5\textwidth]{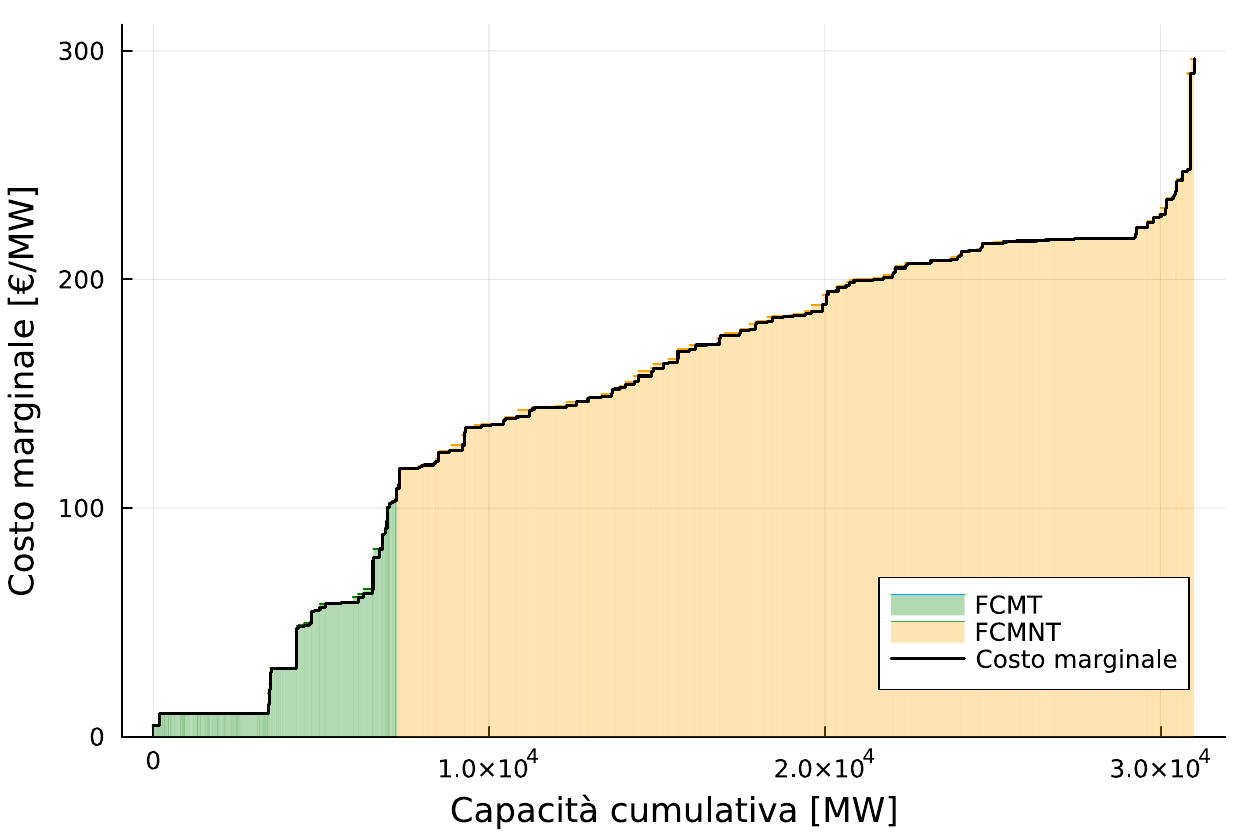}
    \caption{Aggregated marginal cost curve of UPs classified in the 5 clusters.}
    \label{fig:curva-costi-marginali}
\end{figure}

\subsubsection{Operator Portfolios for Market Simulations}

The clustering approach provided pseudo-realistic portfolios. Analysis of production capacity distribution, summarized in Table~\ref{tab:statistiche-operatori}, reveals a concentrated market with 595 UPs for a total capacity of 30,999~MW. Market structure presents a marked competitive asymmetry: OP1 emerges as dominant operator with a capacity of 13,767~MW, equivalent to 44\% of total and distributed over 143~UPs. This is particularly significant if compared with distribution of the remaining nine operators, whose market shares are between 2\% and 12\%. This configuration reflects an oligopolistic structure, where the presence of a dominant operator can significantly influence competitive dynamics and price formation.

\begin{table}[htb]
    \centering
    \caption{Summary statistics of the simulated operators: total capacity, market share, and number of PUs.}
    \label{tab:statistiche-operatori}
    \resizebox{0.7\textwidth}{!}{
    \begin{tabular}{lrrr}
    \toprule
    Operator & \multicolumn{1}{c}{Total capacity [MW]} & \multicolumn{1}{c}{Market share [\%]} & \multicolumn{1}{c}{Number of PUs} \\
    \midrule
    OP1 & 13,767 & 44 & 143 \\
    OP2 & 3,688  & 12 &  64 \\
    OP3 & 3,639  & 12 & 119 \\
    OP4 & 3,213  & 10 &  18 \\
    OP5 & 1,367  &  4 &  42 \\
    OP6 & 1,347  &  4 &  31 \\
    OP7 & 1,246  &  4 &  17 \\
    OP8 & 1,072  &  3 &  10 \\
    OP9 & 1,004  &  3 &  72 \\
    OP10 &   655  &  2 &  79 \\
    \midrule
    Total                 & 30,999 & 100 & 595 \\
    \bottomrule
    \end{tabular}
    }
\end{table}


\smallskip
\noindent
Analysis of technological segmentation reveals an even more articulated market structure, which is strategically relevant for the SPaC mechanism. The split between NNMCS and NMCS technologies highlights a clear predominance of NNMCS plants, which represent 23,724~MW (77\% of total), compared to 7,276~MW (23\%) of NMCS plants. Tables~\ref{tab:statistiche-fcmt} and \ref{tab:statistiche-fcmnt} reveal profoundly different competitive dynamics in the two segments. In the NMCS segment, characterized by lower operational flexibility and lower marginal costs, a rebalancing of competitive positions is observed: OP1 loses its dominant position, holding only 7\% of NMCS capacity with 530~MW. In this segment, leadership of OP3 emerges, which with 2,062~MW controls 28\% of NMCS capacity, followed by OP5 (19\%) and OP6 (17\%). Such redistribution of market power in the NMCS segment suggests greater competitiveness in this compartment, with at least three operators holding significant shares.

\begin{table}[H]
    \centering
    \caption{Capacity breakdown by NMCS technology.}
    \label{tab:statistiche-fcmt}
    \resizebox{0.7\textwidth}{!}{
    \begin{tabular}{lrrr}
    \toprule
    Operator & \multicolumn{1}{c}{NMCS capacity [MW]} & \multicolumn{1}{c}{NMCS share [\%]} & \multicolumn{1}{c}{Number of PUs} \\
    \midrule
    OP3  & 2,062 & 28 & 111 \\
    OP5  & 1,367 & 19 &  42 \\
    OP6  & 1,222 & 17 &  27 \\
    OP4  &   540 &  7 &   7 \\
    OP1  &   530 &  7 &  53 \\
    OP9  &   500 &  7 &  62 \\
    OP10 &   443 &  6 &  58 \\
    OP2  &   435 &  6 &  43 \\
    OP7  &   146 &  2 &  13 \\
    OP8  &    31 &  0 &   6 \\
    \midrule
    Total                 & 7,276 & 100 & 422 \\
    \bottomrule
    \end{tabular}
    }
\end{table}

\begin{table}[H]
    \centering
    \caption{Capacity breakdown by NNMCS technology.}
    \label{tab:statistiche-fcmnt}
    \resizebox{0.7\textwidth}{!}{
    \begin{tabular}{lrrr}
    \toprule
    Operator & \multicolumn{1}{c}{NNMCS capacity [MW]} & \multicolumn{1}{c}{NNMCS share [\%]} & \multicolumn{1}{c}{Number of PUs} \\
    \midrule
    OP1  & 13,237 & 56 &  90 \\
    OP2  &  3,253 & 14 &  21 \\
    OP4  &  2,674 & 11 &  11 \\
    OP3  &  1,577 &  7 &   8 \\
    OP7  &  1,100 &  5 &   4 \\
    OP8  &  1,041 &  4 &   4 \\
    OP9  &    504 &  2 &  10 \\
    OP10 &    212 &  1 &  21 \\
    OP6  &    126 &  1 &   4 \\
    \midrule
    Total                 & 23,724 & 100 & 173 \\
    \bottomrule
    \end{tabular}
    }
\end{table}   

\smallskip
\noindent
In the NNMCS segment, characterized by high marginal cost and more flexible technologies, OP1 instead confirms its dominant position with 13,237~MW, equivalent to 56\% of NNMCS capacity. This concentration results even more marked compared to the overall market, with OP2 placing second with 14\% (3,253~MW) and OP4 third with 11\% (2,674~MW). The asymmetric structure of the two segments presents significant strategic implications for SPaC mechanism implementation. In the NMCS segment, presence of at least three operators with significant shares (OP3 28\%, OP5 19\%, OP6 17\%) suggests possibility of more intense competitive dynamics, with potential benefits in terms of allocative efficiency and price containment. Vice versa, in the NNMCS segment, dominance of OP1 could translate into greater market power, especially during high demand hours when NNMCS plants result marginal for market clearing. This asymmetric oligopolistic configuration between the two segments represents one of the most interesting aspects of the simulated scenario, as it allows evaluating effectiveness of the SPaC mechanism in heterogeneous competitive contexts and analyzing how technological segmentation can influence operator strategies and market outcomes.

\subsubsection{Market Simulation at Marginal Cost}
       
\begin{figure}[htb]
    \centering
    \begin{subfigure}{0.48\textwidth}
    \centering
    \includegraphics[width=\textwidth]{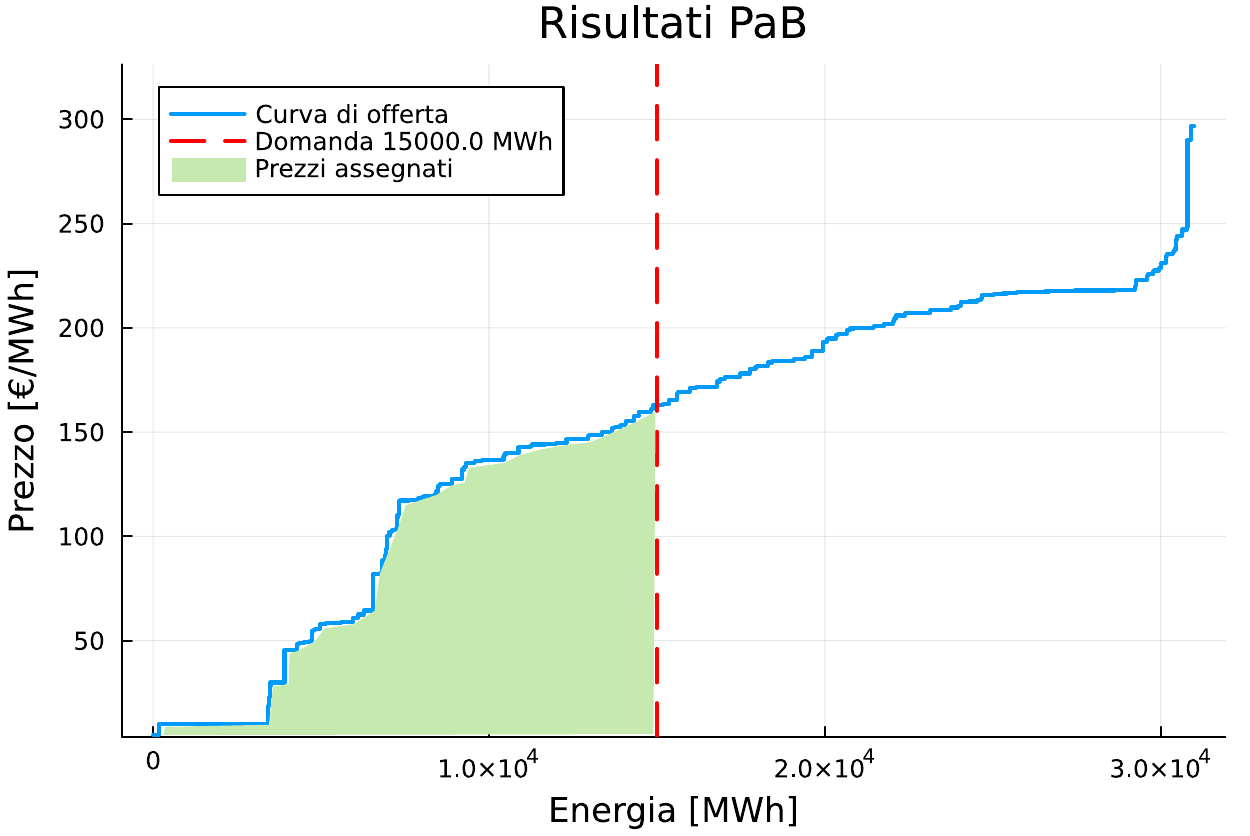}
    \caption{PaB market: offers at marginal cost.}
    \label{fig:risultati-pab-mc-top10}
    \vspace{0.5cm}
    \includegraphics[width=\textwidth]{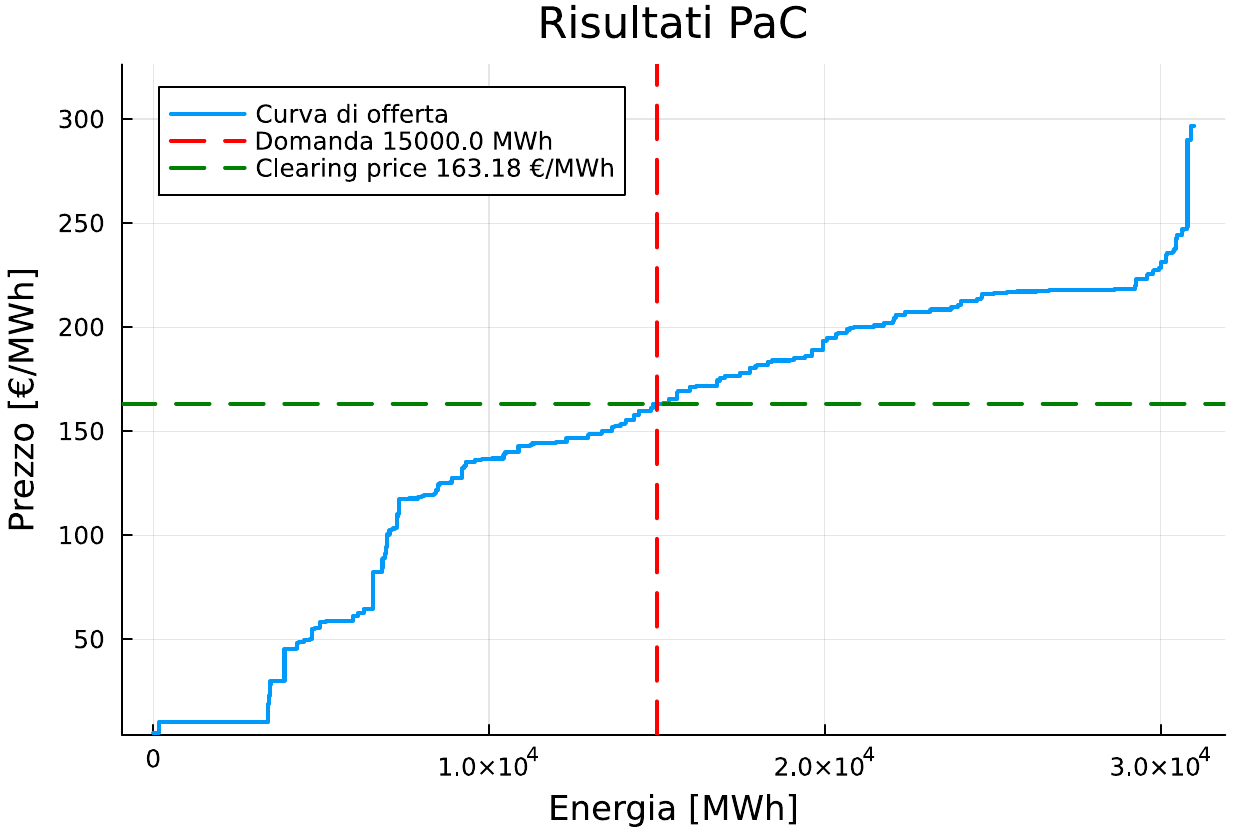}
    \caption{PaC market: offers at marginal cost.}
    \label{fig:risultati-pac-mc-top10}
    \end{subfigure}%
    \hfill
    \begin{subfigure}{0.48\textwidth}
    \centering
    \includegraphics[height=12cm, width=\textwidth]{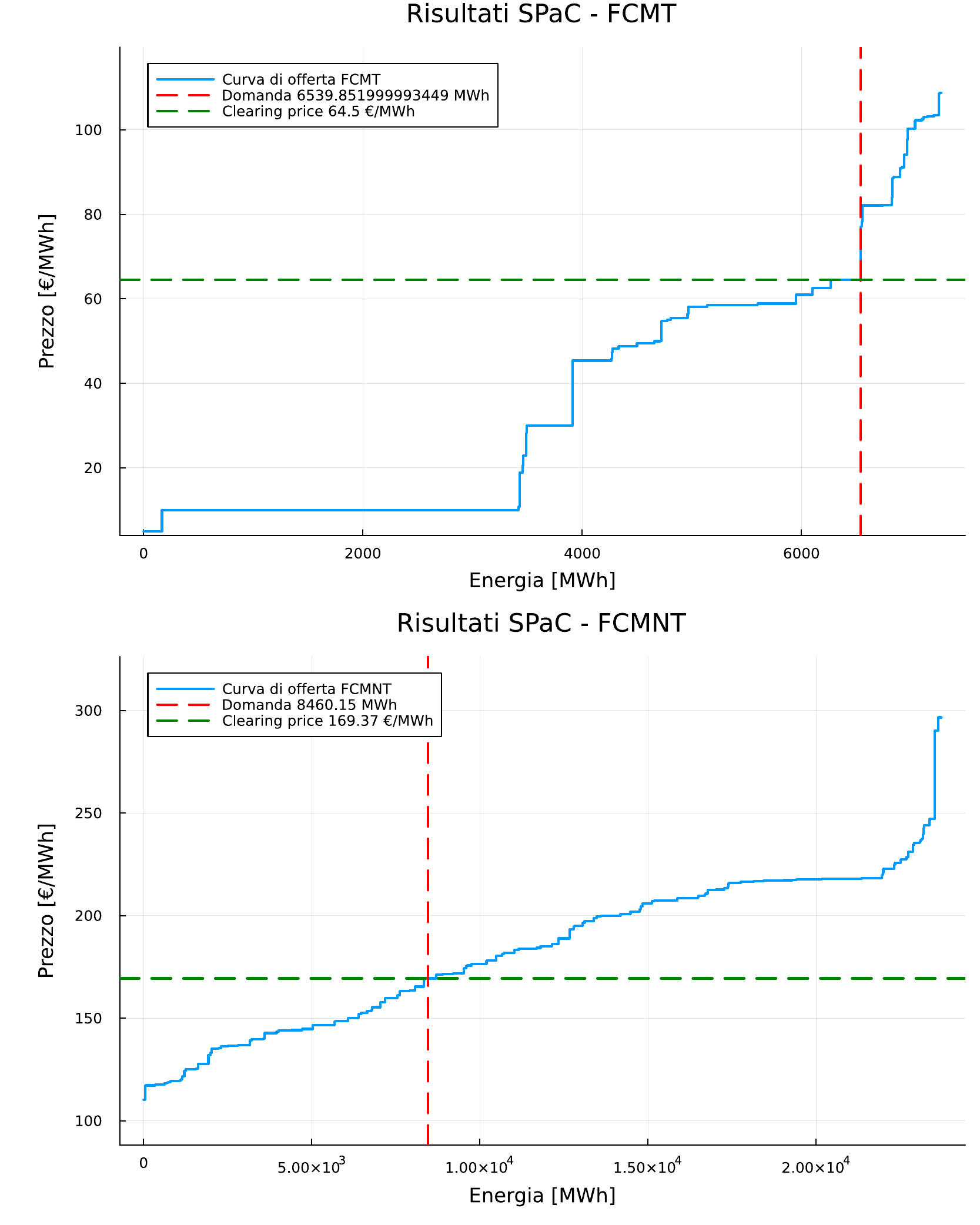}
    \caption{SPaC market: offers at marginal cost.}
    \label{fig:risultati-spac-mc-top10}
    \end{subfigure}
    \caption{Results of market simulations with offers at marginal cost for the 10 operators: left PaB and PaC, right SPaC.}
    \label{fig:risultati-mc-confronto-top10}
\end{figure}

Figure~\ref{fig:risultati-mc-confronto-top10} shows how presence of a dominant operator (OP1 with 44\% of market) significantly influences supply curves, while Table~\ref{tab:risultati-mc-15gw-top10op} quantifies economic impact for a demand level of 15~GW.

\begin{table}[H]
     \centering
     \caption{Simulation results for the three markets under marginal-cost offers for $D = 15$~GW (10-operator scenario).}
     \label{tab:risultati-mc-15gw-top10op}
     \begin{tabular}{lcc}
      \toprule
      \textbf{Market} & \textbf{Total cost [\euro]} & \textbf{PUN [\euro/MWh]} \\
      \midrule
      PaB  & 1\,341\,050 & 89.40 \\
      PaC  & 2\,447\,700 & 163.18 \\
      SPaC & 1\,854\,720 & 123.65 \\
      \bottomrule
     \end{tabular}
\end{table}

\noindent
Quantitative results confirm efficiency order already observed in the PNIEC scenario, but with significantly different absolute values due to real market structure. PaB maintains the most contained cost (1,341,050~\EUR{}, PUN 89.40~\EUR{}/MWh), followed by SPaC (1,854,720~\EUR{}, PUN 123.65~\EUR{}/MWh) and PaC (2,447,700~\EUR{}, PUN 163.18~\EUR{}/MWh). Compared to PNIEC scenario results at same demand, higher PUNs are observed in all three markets, reflecting presence of plants with higher marginal costs in real operator portfolios. Particularly significant is the 38.3\% gap between PaB and SPaC, highlighting how technological segmentation generates an intermediate price premium even under competitive conditions. Analysis of supply curves further reveals influence of market concentration: OP1, dominating the NNMCS segment with 56\% of capacity, largely determines system marginal price in the SPaC market. Vice versa, in the NMCS segment, more distributed competitiveness (with OP3 at 28\%, OP5 at 19\% and OP6 at 17\%) translates into a more fragmented supply curve, potentially more sensitive to operator strategies. These results constitute the baseline for evaluating effectiveness of markup strategies and learning via RL in a realistic oligopolistic context, where asymmetric market concentration can amplify or mitigate effects of strategic behaviors.

\subsubsection{Training of Operators on Multiple Demand Levels}

\begin{figure}[htbp]
      \centering
      \includegraphics[width=0.2\textwidth]{figures/chapter_3/Legenda_markup_policy.png}
      \vspace{0.5cm}
      \begin{subfigure}{\textwidth}
           \centering
           \begin{subfigure}{0.32\textwidth}
                \centering
                \includegraphics[width=\textwidth]{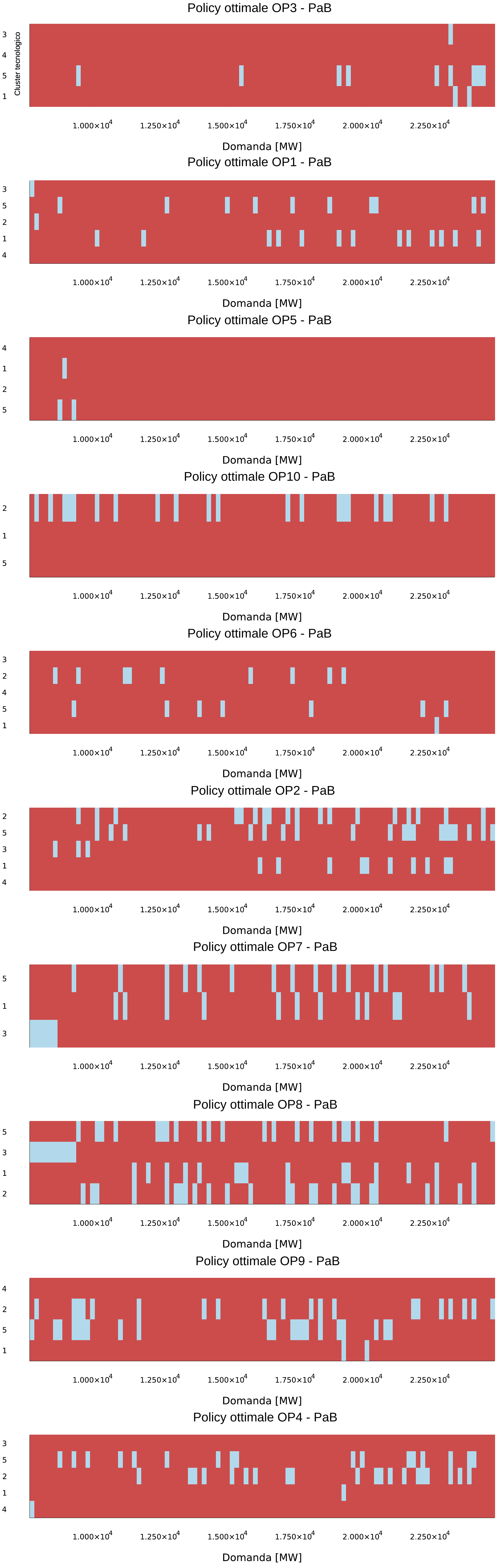}
                \caption{PaB}
           \end{subfigure}
           \hfill
           \begin{subfigure}{0.32\textwidth}
                \centering
                \includegraphics[width=\textwidth]{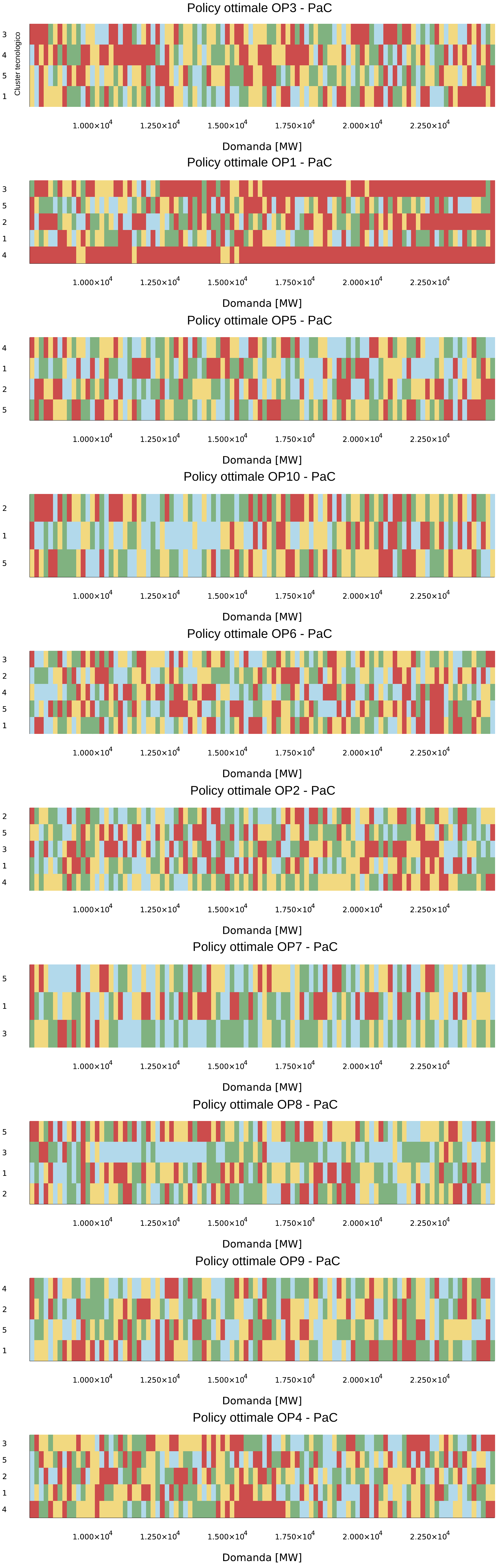}
                \caption{PaC}
           \end{subfigure}
           \hfill
           \begin{subfigure}{0.32\textwidth}
                \centering
                \includegraphics[width=\textwidth]{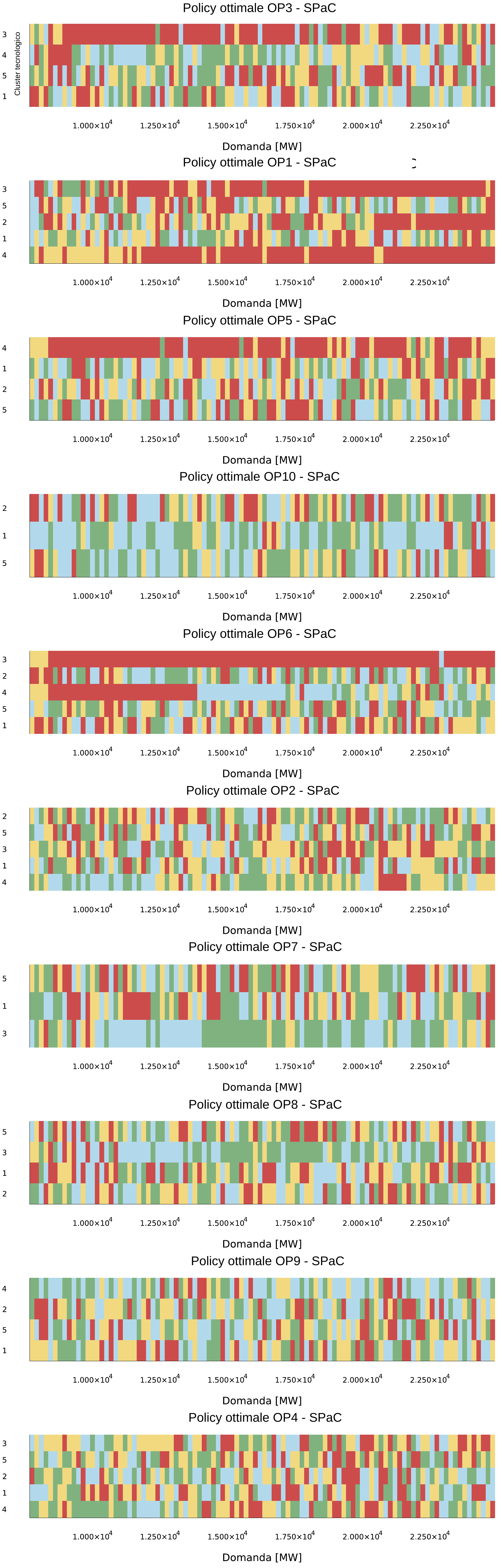}
                \caption{SPaC}
           \end{subfigure}
      \end{subfigure}
      \caption{Optimal policy learned by operators for each technological cluster via RL in the three markets on multiple demand levels (10-operator scenario). Above the markup legend, below policies for PaB, PaC and SPaC.}
      \label{fig:policy-rl-multidemand-top10}
\end{figure}

In the oligopolistic context of ten operators, training via RL was conducted on 100 demand levels between 25\% and 80\% of total market capacity (30,999~MW), for a total of 200,000 episodes distributed among the three market mechanisms. Presence of operators with heterogeneous technological portfolios and significantly different market shares introduces additional strategic complexities compared to the symmetric case. Figure~\ref{fig:policy-rl-multidemand-top10} shows optimal policies learned by the ten operators.

\smallskip
\noindent
In the PaB market, all operators adopt predominantly aggressive markup strategies, confirming what observed in the PNIEC scenario.

\smallskip
\noindent
In the PaC market, high markup behaviors by OP1 are detected for units belonging to technological clusters 3 and 4 (characterized by high marginal costs) over most demand levels, a clear sign that these technologies are often marginal. When the market approaches saturation, high markups also emerge for technological cluster 2, which presents lower marginal costs compared to clusters 3 and 4; this behavior could be interpreted as an economic withholding strategy of capacity to influence market price.

\smallskip
\noindent
In the SPaC regime, a more articulated strategy is observed: for low demand levels, operators OP3, OP5 and OP6 (the main ones in the NMCS segment) adopt aggressive markups for technologies of clusters 3 and 4, while OP1 is more conservative, probably to not compromise sale of NNMCS capacity of which it holds the majority. With increasing demand, OP1 adopts the same the aggressive strategy already observed in PaC, while also OP6 and OP5 apply high markups for technologies of cluster 3, which are marginal in the NMCS market.

\subsubsection{Simulation of a Typical Day in the MGP}

\begin{figure}[htb]
    \centering
    \includegraphics[width=0.8\textwidth]{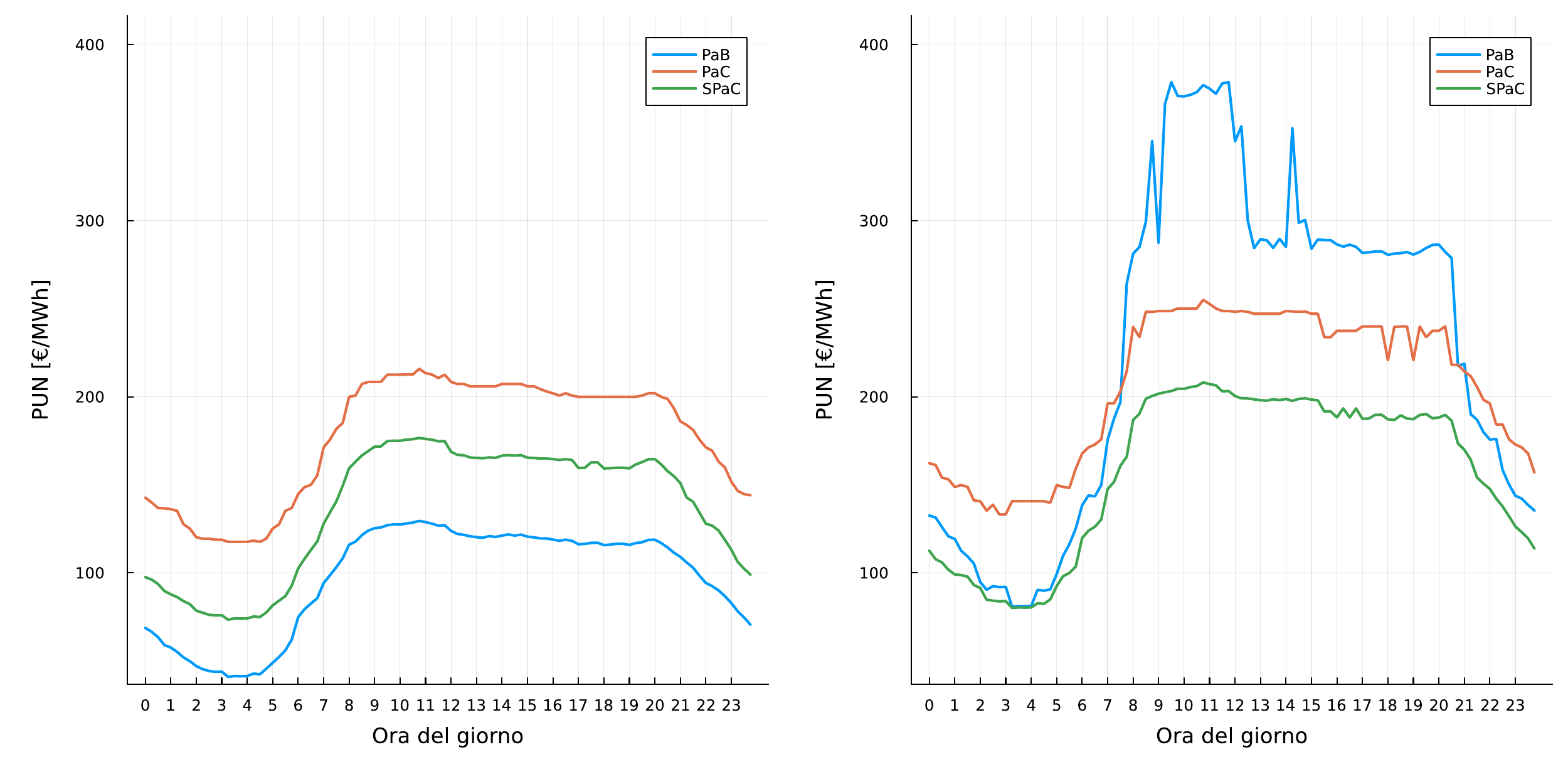}
    \caption{Simulated daily PUN in the three markets for the two strategies (10 operators)}
    \label{fig:pun-giornaliero-top10op}
\end{figure}

Using optimal policies learned via RL, simulation of a typical day in the MGP is performed to evaluate the behavior of the ten operators in operational conditions that may occur in the Italian electricity market. Simulation uses the same load curve of September 12, 2025, scaled between 25\% and 80\% of total market capacity (30,999~MW) to represent daily demand variability. Figure~\ref{fig:pun-giornaliero-top10op} illustrates the PUN evolution over 24 hours for the three markets, comparing marginal cost strategies with those optimized via RL.

\smallskip
\noindent
In the marginal cost case, the already observed efficiency ordering is confirmed: PaB yields the lowest PUN (varying between 52 and 129~\EUR{}/MWh), followed by SPaC (84-177~\EUR{}/MWh) and PaC (128-216~\EUR{}/MWh). With RL strategies, all markets show significant PUN increases, particularly pronounced for PaB which reaches peaks of 377~\EUR{}/MWh during maximum demand hours. SPaC shows greater stability, with PUN varying between 98 and 208~\EUR{}/MWh, confirming its ability to limit market power exercise thanks to technological segmentation.

\paragraph{Total Costs}
        
Analysis of total costs for different demand levels (Table~\ref{tab:costi-tot-domanda-top10op}) reveals the impact of oligopolistic structure on economic results. In the marginal cost regime, costs follow the expected ordering with significantly higher absolute values compared to the PNIEC scenario, reflecting the presence of plants with higher marginal costs in the operator portfolios. PaB maintains maximum efficiency with daily costs of 18.0~M\EUR{}, while PaC reaches 31.8~M\EUR{} (+76\%) and SPaC stands at 25.0~M\EUR{} (+39\%). With RL strategies, cost increases are particularly marked for PaB (43.7~M\EUR{}, +142\%), highlighting how operators manage to effectively exploit the pay-as-bid mechanism to maximize profits.
PaC shows a more contained increase (37.2~M\EUR{}, +17\%)---but, again, from a much higher baseline---, while SPaC registers a 16\% increase (29.2~M\EUR{}). PaB/SPaC and PaC/SPaC ratios confirm the superiority of the segmented mechanism: in the RL case, SPaC results systematically more efficient, with PaB costing up to 1.8 times more during peak hours and PaC up to 1.5 times more.

\begin{table}[H]
    \caption{10-operator scenario: total costs in the three markets for different demand levels.}
    \label{tab:costi-tot-domanda-top10op}
    \centering
    \resizebox{\textwidth}{!}{%
    \begin{tabular}{@{}lllllllllllll@{}}
    \toprule
    Market   Capacity {[}MW{]} &
    31000 &
    31000 &
    31000 &
    31000 &
    31000 &
    31000 &
    31000 &
    31000 &
    31000 &
    31000 &
    31000 &
    31000 \\
    D {[}MW{]} &
    9009 &
    10236 &
    12986 &
    14230 &
    16970 &
    19212 &
    21579 &
    23229 &
    23747 &
    24236 &
    24388 &
    24800 \\
    D/Market Capacity {[}\%{]} &
    29 &
    33 &
    42 &
    46 &
    55 &
    62 &
    70 &
    75 &
    77 &
    78 &
    79 &
    80 \\ \midrule
    \multicolumn{13}{c}{\textbf{Strategy: marginal-cost offers}} \\ 
    Cost tot PaB {[}\euro{]} &
    470,813.56 &
    635,850.53 &
    1,030,207.80 &
    1,218,241.87 &
    1,674,332.08 &
    2,079,937.79 &
    2,540,726.74 &
    2,879,586.53 &
    2,987,637.99 &
    3,090,663.10 &
    3,122,970.49 &
    3,211,129.77 \\
    Cost tot PaC {[}\euro{]} &
    1,149,695.74 &
    1,400,734.15 &
    1,930,693.11 &
    2,210,497.11 &
    2,981,279.89 &
    3,555,758.76 &
    4,331,120.38 &
    4,842,334.37 &
    4,950,385.84 &
    5,151,585.75 &
    5,185,836.76 &
    5,354,055.94 \\
    Cost tot SPaC {[}\euro{]} &
    758,179.51 &
    949,376.32 &
    1,402,087.68 &
    1,676,708.17 &
    2,279,524.27 &
    2,872,357.07 &
    3,519,287.34 &
    3,929,896.54 &
    4,080,836.35 &
    4,242,223.15 &
    4,284,131.55 &
    4,380,787.01 \\ \midrule
    \multicolumn{13}{c}{\textbf{Strategy: RL policy}} \\ 
    Cost tot PaB {[}\euro{]} &
    988,148.84 &
    1,281,222.41 &
    1,869,979.46 &
    2,132,940.55 &
    3,182,855.29 &
    5,078,735.48 &
    6,152,267.68 &
    8,023,091.55 &
    8,696,349.43 &
    8,987,342.22 &
    9,060,148.89 &
    9,350,667.36 \\
    Cost tot PaC {[}\euro{]} &
    1,341,074.44 &
    1,629,760.94 &
    2,224,108.60 &
    2,501,829.08 &
    3,330,114.07 &
    4,119,496.76 &
    5,047,413.58 &
    5,766,317.56 &
    5,906,551.70 &
    6,062,662.21 &
    6,100,674.56 &
    6,325,766.62 \\
    Cost tot SPaC {[}\euro{]} &
    883,522.05 &
    1,060,613.58 &
    1,609,214.82 &
    1,854,197.89 &
    2,573,167.56 &
    3,189,539.29 &
    4,109,372.88 &
    4,657,196.52 &
    4,810,690.49 &
    4,957,932.92 &
    5,012,078.81 &
    5,162,281.34 \\
    Cost PaB/Cost SPaC &
    \cellcolor[HTML]{63BE7B}1.1 &
    \cellcolor[HTML]{83C77C}1.2 &
    \cellcolor[HTML]{72C27B}1.2 &
    \cellcolor[HTML]{6EC17B}1.2 &
    \cellcolor[HTML]{8ECA7D}1.2 &
    \cellcolor[HTML]{FED480}1.6 &
    \cellcolor[HTML]{EDE582}1.5 &
    \cellcolor[HTML]{FB9574}1.7 &
    \cellcolor[HTML]{F96C6C}1.8 &
    \cellcolor[HTML]{F8696B}1.8 &
    \cellcolor[HTML]{F96C6C}1.8 &
    \cellcolor[HTML]{F96A6C}1.8 \\
    Cost PaC/Cost SPaC &
    \cellcolor[HTML]{F9726D}1.5 &
    \cellcolor[HTML]{F8696B}1.5 &
    \cellcolor[HTML]{FCB37A}1.4 &
    \cellcolor[HTML]{FDC37D}1.3 &
    \cellcolor[HTML]{FFDD82}1.3 &
    \cellcolor[HTML]{FFDF82}1.3 &
    \cellcolor[HTML]{87C87D}1.2 &
    \cellcolor[HTML]{A7D17E}1.2 &
    \cellcolor[HTML]{85C87D}1.2 &
    \cellcolor[HTML]{75C37C}1.2 &
    \cellcolor[HTML]{63BE7B}1.2 &
    \cellcolor[HTML]{7DC57C}1.2 \\ \bottomrule
    \end{tabular}%
    }
\end{table}

\paragraph{Operator Profits}
        
Table~\ref{tab:profitti-domanda-top10op} highlights how asymmetric oligopolistic structure significantly influences profit distribution among operators. In the marginal cost regime, PaB generates no profits (by definition), while PaC and SPaC produce respectively 13.7~M\EUR{} and 6.4~M\EUR{} of overall profits. Profit concentration reflects market structure: in PaC, OP3 benefits most (3.4~M\EUR{}) thanks to its dominance in the NMCS segment, followed by OP1 (2.5~M\EUR{}) which exploits its broad overall capacity. With RL strategies, profits grow substantially in all markets, reaching 24.9~M\EUR{} in PaB, 19.1~M\EUR{} in PaC and 10.6~M\EUR{} in SPaC. OP1 emerges as the main beneficiary in PaB (9.4~M\EUR{}), exploiting its dominant position and ability to influence price through aggressive markup strategies. In PaC, OP3 maintains leadership (4.3~M\EUR{}) followed by OP1 (4.0~M\EUR{}), while in SPaC OP1 recovers first position (3.5~M\EUR{}). PaB/SPaC and PaC/SPaC ratios confirm that SPaC effectively limits rent extraction: profits in PaB result up to 3.3 times higher compared to SPaC, while PaC generates profits up to 2.6 times greater.

\begin{table}[H]
    \caption{10-operator scenario: operator profits in the three markets for different demand levels.}
    \label{tab:profitti-domanda-top10op}
    \centering
    \resizebox{\textwidth}{!}{%
    \begin{tabular}{@{}lllllllllllll@{}}
    \toprule
    Market   Capacity {[}MW{]} &
    31000 &
    31000 &
    31000 &
    31000 &
    31000 &
    31000 &
    31000 &
    31000 &
    31000 &
    31000 &
    31000 &
    31000 \\
    D {[}MW{]} &
    9009 &
    10236 &
    12986 &
    14230 &
    16970 &
    19212 &
    21579 &
    23229 &
    23747 &
    24236 &
    24388 &
    24800 \\
    D/Market Capacity {[}\%{]} &
    29 &
    33 &
    42 &
    46 &
    55 &
    62 &
    70 &
    75 &
    77 &
    78 &
    79 &
    80 \\ \midrule
    \multicolumn{13}{c}{\textbf{Strategy: marginal-cost offers}} \\
    Operator profits PaB   {[}\euro{]} &
    - &
    - &
    - &
    - &
    - &
    - &
    - &
    - &
    - &
    - &
    - &
    - \\
    Operator profits PaC   {[}\euro{]} &
    678,882.18 &
    764,883.62 &
    900,485.31 &
    992,255.24 &
    1,306,947.81 &
    1,475,820.97 &
    1,790,393.64 &
    1,962,747.85 &
    1,962,747.85 &
    2,060,922.65 &
    2,062,866.27 &
    2,142,926.16 \\
    Operator profits SPaC   {[}\euro{]} &
    257,334.91 &
    278,885.44 &
    329,786.68 &
    409,115.25 &
    542,885.94 &
    722,590.38 &
    897,046.66 &
    964,526.49 &
    1,005,446.69 &
    1,061,858.91 &
    1,070,916.69 &
    1,078,349.01 \\ \midrule
    \multicolumn{13}{c}{\textbf{Strategy: RL policy}} \\
    Operator profits PaB   {[}\euro{]} &
    462,897.99 &
    612,986.63 &
    797,482.05 &
    872,096.31 &
    1,449,990.68 &
    2,833,648.00 &
    3,450,171.10 &
    5,135,968.65 &
    5,701,721.43 &
    5,890,923.86 &
    5,932,053.18 &
    6,135,752.41 \\
    Operator profits PaC   {[}\euro{]} &
    866,363.23 &
    990,219.12 &
    1,176,719.91 &
    1,272,088.51 &
    1,641,952.46 &
    2,028,185.15 &
    2,500,941.34 &
    2,876,284.87 &
    2,902,233.21 &
    2,963,372.59 &
    2,963,437.26 &
    3,110,895.03 \\
    Operator profits SPaC   {[}\euro{]} &
    380,721.96 &
    386,117.28 &
    518,987.35 &
    591,458.90 &
    827,109.30 &
    1,024,785.21 &
    1,523,855.42 &
    1,670,174.29 &
    1,727,764.20 &
    1,771,386.61 &
    1,825,825.72 &
    1,876,088.15 \\
    Profits PaB/Profits SPaC &
    \cellcolor[HTML]{63BE7B}1.2 &
    \cellcolor[HTML]{8FCA7D}1.6 &
    \cellcolor[HTML]{89C97D}1.5 &
    \cellcolor[HTML]{82C67C}1.5 &
    \cellcolor[HTML]{A3D07E}1.8 &
    \cellcolor[HTML]{FDC37D}2.8 &
    \cellcolor[HTML]{E0E282}2.3 &
    \cellcolor[HTML]{FB9273}3.1 &
    \cellcolor[HTML]{F96E6C}3.3 &
    \cellcolor[HTML]{F8696B}3.3 &
    \cellcolor[HTML]{F9766E}3.2 &
    \cellcolor[HTML]{F9726D}3.3 \\
    Profits PaC/Profits SPaC &
    \cellcolor[HTML]{FB9E76}2.3 &
    \cellcolor[HTML]{F8696B}2.6 &
    \cellcolor[HTML]{FBA076}2.3 &
    \cellcolor[HTML]{FDB57A}2.2 &
    \cellcolor[HTML]{FED380}2.0 &
    \cellcolor[HTML]{FED480}2.0 &
    \cellcolor[HTML]{6FC17B}1.6 &
    \cellcolor[HTML]{A6D17E}1.7 &
    \cellcolor[HTML]{89C97D}1.7 &
    \cellcolor[HTML]{85C77C}1.7 &
    \cellcolor[HTML]{63BE7B}1.6 &
    \cellcolor[HTML]{7BC47C}1.7 \\ \bottomrule
    \end{tabular}%
    }
\end{table}

\paragraph{PUN Evolution}
        
Table~\ref{tab:pun-domanda-top10op} quantifies PUN evolution varying demand levels, highlighting impact of oligopolistic strategies on market prices. In the marginal cost regime, PUN grows progressively with demand reflecting technological merit order: PaB varies from 52.3 to 129.5~\EUR{}/MWh, SPaC from 84.2 to 176.7~\EUR{}/MWh, and PaC from 127.6 to 215.9~\EUR{}/MWh. These values result systematically higher compared to the PNIEC scenario, confirming influence of oligopolistic structure and higher marginal costs of the realistic portfolios. With RL strategies, price increases are dramatic in PaB, which reaches peaks of 377~\EUR{}/MWh during maximum demand hours (+191\% compared to marginal cost). PaC shows more contained increases, varying from 148.9 to 255.1~\EUR{}/MWh (+18\% on average), while SPaC maintains a similar stability with PUN between 98.1 and 208.2~\EUR{}/MWh (+18\% on average) and lower costs.

\smallskip
\noindent
SPaC's ability to limit price volatility clearly emerges from the comparison: while in PaB PUN increases 7.2 times from minimum to maximum demand with RL strategies, in SPaC the increase is limited to 2.1 times, guaranteeing greater predictability for consumers and producers.

\begin{table}[H]
    \caption{10-operator scenario: PUN in the three markets for different demand levels.}
    \label{tab:pun-domanda-top10op}
    \centering
    \resizebox{\textwidth}{!}{%
    \begin{tabular}{@{}lllllllllllll@{}}
    \toprule
    Market   Capacity {[}MW{]} &
    31000 &
    31000 &
    31000 &
    31000 &
    31000 &
    31000 &
    31000 &
    31000 &
    31000 &
    31000 &
    31000 &
    31000 \\
    D {[}MW{]} &
    9009 &
    10236 &
    12986 &
    14230 &
    16970 &
    19212 &
    21579 &
    23229 &
    23747 &
    24236 &
    24388 &
    24800 \\
    D/Market Capacity {[}\%{]} &
    29 &
    33 &
    42 &
    46 &
    55 &
    62 &
    70 &
    75 &
    77 &
    78 &
    79 &
    80 \\ \midrule
    \multicolumn{13}{c}{\textbf{Strategy: marginal-cost offers}} \\
    PUN PaB {[}\euro/MWh{]} &
    \cellcolor[HTML]{FCFCFF}52.26 &
    \cellcolor[HTML]{FCF4F7}62.12 &
    \cellcolor[HTML]{FCE4E7}79.33 &
    \cellcolor[HTML]{FCDFE1}85.61 &
    \cellcolor[HTML]{FBD3D6}98.66 &
    \cellcolor[HTML]{FBCACD}108.26 &
    \cellcolor[HTML]{FBC2C4}117.74 &
    \cellcolor[HTML]{FBBCBF}123.96 &
    \cellcolor[HTML]{FBBABD}125.81 &
    \cellcolor[HTML]{FBB9BB}127.52 &
    \cellcolor[HTML]{FBB8BB}128.05 &
    \cellcolor[HTML]{FBB7BA}129.48 \\
    PUN PaC {[}\euro/MWh{]} &
    \cellcolor[HTML]{FBB9BB}127.61 &
    \cellcolor[HTML]{FAB1B3}136.84 &
    \cellcolor[HTML]{FAA6A8}148.67 &
    \cellcolor[HTML]{FAA0A2}155.34 &
    \cellcolor[HTML]{F98E90}175.68 &
    \cellcolor[HTML]{F98587}185.08 &
    \cellcolor[HTML]{F97779}200.71 &
    \cellcolor[HTML]{F97072}208.46 &
    \cellcolor[HTML]{F97072}208.46 &
    \cellcolor[HTML]{F96C6F}212.56 &
    \cellcolor[HTML]{F96C6E}212.64 &
    \cellcolor[HTML]{F8696B}215.89 \\
    PUN SPaC {[}\euro/MWh{]} &
    \cellcolor[HTML]{FCE0E3}84.15 &
    \cellcolor[HTML]{FCD8DB}92.75 &
    \cellcolor[HTML]{FBCACD}107.97 &
    \cellcolor[HTML]{FBC2C4}117.83 &
    \cellcolor[HTML]{FAB3B5}134.33 &
    \cellcolor[HTML]{FAA5A8}149.51 &
    \cellcolor[HTML]{FA999B}163.09 &
    \cellcolor[HTML]{FA9396}169.18 &
    \cellcolor[HTML]{FA9193}171.84 &
    \cellcolor[HTML]{F98E90}175.04 &
    \cellcolor[HTML]{F98E90}175.67 &
    \cellcolor[HTML]{F98D8F}176.65 \\ \midrule
    \multicolumn{13}{c}{\textbf{Strategy: RL policy}} \\
    PUN PaB {[}\euro/MWh{]} &
    \cellcolor[HTML]{FCF6F9}109.68 &
    \cellcolor[HTML]{FCEEF1}125.16 &
    \cellcolor[HTML]{FCE4E7}143.99 &
    \cellcolor[HTML]{FCE1E4}149.89 &
    \cellcolor[HTML]{FBCDD0}187.56 &
    \cellcolor[HTML]{FAA5A7}264.35 &
    \cellcolor[HTML]{FA9A9C}285.10 &
    \cellcolor[HTML]{F97A7C}345.39 &
    \cellcolor[HTML]{F96F71}366.20 &
    \cellcolor[HTML]{F96D6F}370.83 &
    \cellcolor[HTML]{F96C6E}371.50 &
    \cellcolor[HTML]{F8696B}377.04 \\
    PUN PaC {[}\euro/MWh{]} &
    \cellcolor[HTML]{FCE2E5}148.85 &
    \cellcolor[HTML]{FCDCDF}159.21 &
    \cellcolor[HTML]{FBD6D9}171.26 &
    \cellcolor[HTML]{FBD4D6}175.81 &
    \cellcolor[HTML]{FBC9CB}196.24 &
    \cellcolor[HTML]{FBBFC2}214.42 &
    \cellcolor[HTML]{FBB5B7}233.90 &
    \cellcolor[HTML]{FAADB0}248.24 &
    \cellcolor[HTML]{FAADB0}248.72 &
    \cellcolor[HTML]{FAACAF}250.15 &
    \cellcolor[HTML]{FAACAF}250.15 &
    \cellcolor[HTML]{FAAAAC}255.07 \\
    PUN SPaC {[}\euro/MWh{]} &
    \cellcolor[HTML]{FCFCFF}98.07 &
    \cellcolor[HTML]{FCFAFD}103.61 &
    \cellcolor[HTML]{FCEFF2}123.92 &
    \cellcolor[HTML]{FCECEE}130.30 &
    \cellcolor[HTML]{FCE0E3}151.63 &
    \cellcolor[HTML]{FCD9DB}166.02 &
    \cellcolor[HTML]{FBCCCE}190.43 &
    \cellcolor[HTML]{FBC7C9}200.49 &
    \cellcolor[HTML]{FBC5C8}202.58 &
    \cellcolor[HTML]{FBC4C7}204.57 &
    \cellcolor[HTML]{FBC4C6}205.52 &
    \cellcolor[HTML]{FBC2C5}208.16 \\ \bottomrule
    \end{tabular}%
    }
\end{table}

\subsubsection{Synthetic Results}

Table~\ref{tab:risultati-sintetici-top10op} summarizes overall results of the daily simulation for the ten-operator scenario. The data confirm the competitive dynamics already observed in previous analyses. The PaB market maintains maximum efficiency at marginal costs (18.0~M\EUR{}), but becomes the least efficient with optimized strategies (43.7~M\EUR{}), registering a 142\% increase. The PaC market presents high costs in both configurations (from 31.8~M\EUR{} to 37.2~M\EUR{}, +17\%), confirming relative stability of the uniform mechanism even in presence of strategic behaviors. The SPaC market confirms itself as the most balanced compromise solution, with costs increasing from 25.0~M\EUR{} to 29.2~M\EUR{} (+16\%).

\smallskip
\noindent
The Profit/Cost ratio reveals the intensity of rent extraction in each mechanism. With marginal cost strategies, PaB presents a ratio of 0\% (no profit by definition), PaC a ratio of 43\% and SPaC a ratio of 26\%. Introduction of RL strategies substantially amplifies these ratios: PaB reaches 57\%, PaC 51\% and SPaC 36\%. Compared to the PNIEC scenario, greater rent extraction capability is observed in all mechanisms, likely reflecting the influence of market concentration.

\begin{table}[H]
    \caption{Summary results of the representative-day simulation (10-operator scenario).}
    \label{tab:risultati-sintetici-top10op}
    \centering
    \resizebox{\textwidth}{!}{%
    \begin{tabular}{@{}lllllll@{}}
        \toprule
        & \multicolumn{3}{c}{\textbf{Strategy: marginal-cost offers}}             & \multicolumn{3}{c}{\textbf{Strategy: policy RL}}                   \\ \midrule
        \textbf{Market} & Market cost {[}\euro{]} \footnotemark & Operators' profits  {[}\euro{]} & Profit/Cost & Market cost {[}\euro{]} & Operators' profits {[}\euro{]} & Profit/Cost \\
        PaB  & 180,336,057.75 & -              & 0\%  & 437,156,222.87 & 249,299,586.43 & 57\% \\
        PaC  & 317,507,254.16 & 137,171,196.41 & 43\% & 372,141,894.80 & 191,242,362.14 & 51\% \\
        SPaC & 250,410,346.73 & 63,962,890.97  & 26\% & 291,782,876.47 & 105,885,305.32 & 36\% \\ \bottomrule
    \end{tabular}%
    }
\end{table}
\footnotetext{By \emph{Market cost} we mean the total expenditure sustained by consumers for purchasing energy in the market which, given demand inelasticity, coincides with aggregated producer revenues.}

\noindent
Profit distribution among operators (Table~\ref{tab:dettaglio-profitti-top10op}) reveals the impact of asymmetric oligopolistic structure and technological segmentation on individual results. In the marginal cost regime, OP3 emerges as the main beneficiary in PaC (3.4~M\EUR{}) and SPaC (1.2~M\EUR{} for the NMCS segment), exploiting its leadership in the flexible technologies segment. OP1, despite overall dominance, obtains more contained profits (2.5~M\EUR{} in PaC, 2.1~M\EUR{} in SPaC) due to higher marginal costs of its diversified portfolio. With RL strategies, OP1 recovers leadership in profit generation, reaching 9.4~M\EUR{} in PaB, 4.0~M\EUR{} in PaC and 3.5~M\EUR{} in SPaC. This trend reversal highlights how size and technological diversification allow OP1 to effectively exploit strategic opportunities offered by different market mechanisms. OP3 still maintains a competitive position (2.8~M\EUR{} in PaB, 4.3~M\EUR{} in PaC, 1.8~M\EUR{} in SPaC), confirming advantages of specialization in the NMCS segment. Medium-sized operators (OP2, OP4) benefit significantly from the PaB mechanism, where they manage to generate profits comparable or superior to those of PaC, while minor operators (OP10, OP8) show more limited but still positive performance in all markets.

\begin{table}[H]
    \caption{Operator profit breakdown, 10-operator scenario\protect\footnotemark}
    \label{tab:dettaglio-profitti-top10op}
    \centering
    \resizebox{0.5\textwidth}{!}{%
    \begin{tabular}{@{}llll@{}}
        \toprule
        \multicolumn{1}{c}{\multirow{2}{*}{Operator}} & \multicolumn{3}{c}{Profits {[}\euro{]}}                                         \\
        \multicolumn{1}{c}{}                           & \multicolumn{1}{c}{PaB} & \multicolumn{1}{c}{PaC} & \multicolumn{1}{c}{SPaC} \\ \midrule
        \multicolumn{4}{c}{\textbf{Strategy: marginal-cost offers}}              \\
        OP3  & -             & 33,553,430.19 & 11,883,888.13 \\
        OP1  & -             & 25,264,287.89 & 20,986,500.80 \\
        OP5  & -             & 19,848,563.86 & 5,096,571.41  \\
        OP10 & -             & 5,493,251.45  & 1,644,156.57  \\
        OP6  & -             & 14,898,113.24 & 2,071,088.24  \\
        OP2  & -             & 9,028,923.37  & 4,721,999.79  \\
        OP7  & -             & 4,254,651.89  & 2,885,291.68  \\
        OP8  & -             & 1,931,484.21  & 1,775,052.12  \\
        OP9  & -             & 8,568,553.55  & 3,814,566.76  \\
        OP4  & -             & 14,329,936.76 & 9,083,775.48  \\ \midrule
        \multicolumn{4}{c}{\textbf{Strategy: RL policy}}                      \\
        OP3  & 28,260,462.69 & 42,951,576.46 & 18,250,480.35 \\
        OP1  & 93,647,903.34 & 39,839,953.16 & 34,505,630.62 \\
        OP5  & 6,411,564.68  & 23,745,390.95 & 7,120,343.96  \\
        OP10 & 6,775,066.51  & 7,298,250.29  & 2,408,818.62  \\
        OP6  & 12,557,217.19 & 18,704,649.37 & 3,915,226.72  \\
        OP2  & 30,681,323.79 & 14,458,462.35 & 9,675,729.93  \\
        OP7  & 13,530,948.98 & 7,076,196.41  & 5,355,600.15  \\
        OP8  & 13,151,044.68 & 4,229,094.14  & 3,911,201.15  \\
        OP9  & 10,638,968.53 & 11,307,888.31 & 5,471,703.65  \\
        OP4  & 33,645,086.05 & 21,630,900.71 & 15,270,570.15 \\ \bottomrule
    \end{tabular}%
    }
\end{table}
\footnotetext{Note that the results are based on the public offers available on the GME website, but processed using PU-type clustering and marginal-cost estimates. Therefore, the values reported are limited to the simulations considered in this report and do not necessarily correspond to operators' actual profits in the real market.}

\begin{figure}[htb]
    \centering
    \includegraphics[width=0.75\textwidth]{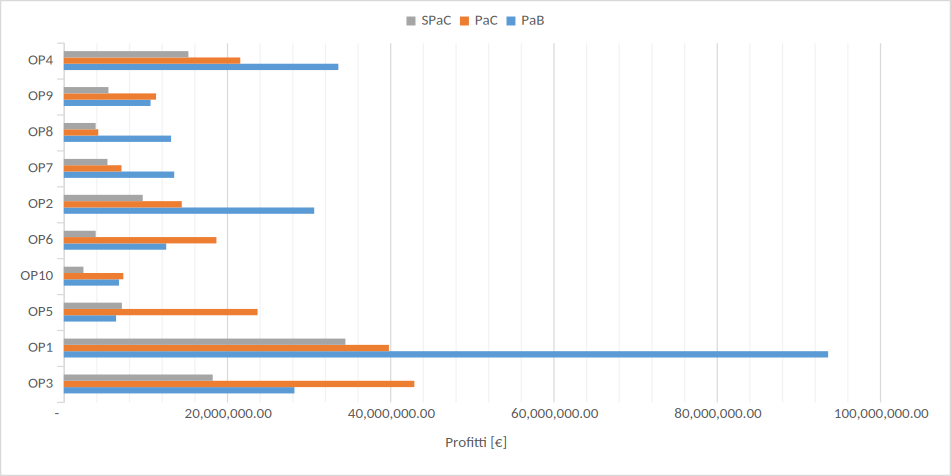}
    \caption{Comparison of operators' profits among the three markets with optimal strategy learned via RL (10-operator scenario)}
    \label{fig:profitti-top10op-mercati}
\end{figure}

\noindent
Figure~\ref{fig:profitti-top10op-mercati} visually summarizes profit distribution among operators in the three mechanisms with RL strategies, confirming the previous analytical observations: OP1 dominates in PaB thanks to its overall capacity, while in PaC competition results in more balanced profits among major operators. SPaC emerges as the most equitable mechanism, limiting both absolute profits and relative disparities, suggesting a potential benefit in terms of political acceptability and market competitive sustainability.

\smallskip
\noindent
In summary, results of the ten-operator scenario confirm effectiveness of the SPaC mechanism as a compromise solution between economic efficiency and market power control, resulting particularly advantageous in oligopolistic contexts characterized by high concentration and significant competitive asymmetries. SPaC's ability to limit excessive rent extraction while maintaining participation incentives for all operators makes it a promising proposal for reform of the Italian electricity market.

\subsubsection{Simulation of a Year in the MGP}

Extending the reasoning developed for typical day simulations, actual load data for the entire year 2024 publicly available on the Terna website, has been used to simulate electricity market evolution on an annual basis. For each operator, the same policy learned previously is adopted (Figure~\ref{fig:policy-rl-multidemand-top10}). Additional hypotheses introduced in the simulation are:
\begin{itemize}
 \item annual load curve scaled between 25\% and 80\% of overall market capacity;
 \item offered capacities and marginal costs of production units kept constant over the entire year, in absence of seasonal variability of various nature.
\end{itemize}
The simulation is conducted for each of the three considered market configurations, on intra-daily basis with quarter-hour resolution, for a total of 35,136 market outcomes for each market over the year.

\begin{figure}[htb]
    \centering
    \includegraphics[width=0.5\linewidth]{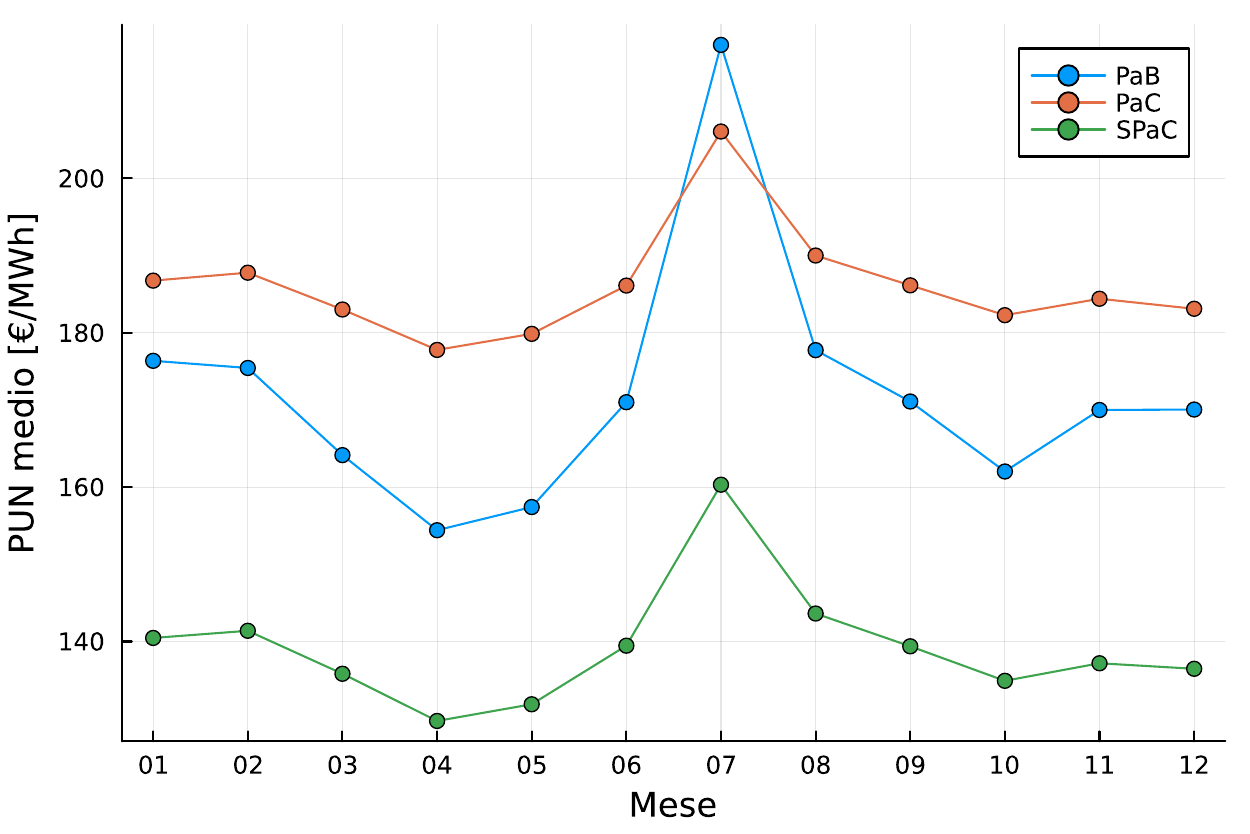}
    \caption{Average monthly PUN trend in different simulation configurations}
    \label{fig:PUN-2024-sim}
\end{figure}

\smallskip
\noindent
In Figure~\ref{fig:PUN-2024-sim}, the average monthly PUN trend in the different configurations is reported. In all cases, a clear seasonal dynamic emerges, with lower values in spring and autumn months and a pronounced peak in July, consistent with summer demand increase. The PaB configuration presents the highest average PUN levels and the greatest variability. The PaC configuration shows an analogous trend, with lower average values and more contained variability. The SPaC configuration instead evidences systematically lower PUN levels throughout the entire year, while maintaining the same seasonal structure.

\noindent
To provide a quantitative synthesis of effects of the SPaC configuration compared to reference cases PaB and PaC, the following indicators are used.

\paragraph{PUN Reduction (min/average/max)}
    Unique National Price variation is calculated quarter hour by quarter hour as:
    \[
    \Delta PUN_t = 
    \frac{PUN_{\text{SPaC},t} - PUN_{\text{base},t}}{PUN_{\text{base},t}} \times 100\%,
    \quad \text{with } \text{base} \in \{\text{PaB}, \text{PaC}\}.
    \]
    Starting from these values, average, minimum and maximum values over the entire year are defined:
    \[
    \overline{\Delta PUN} = \frac{1}{N}\sum_{t=1}^{N} \Delta PUN_t,
    \qquad
    \Delta PUN_{\min} = \min_t \Delta PUN_t,
    \qquad
    \Delta PUN_{\max} = \max_t \Delta PUN_t.
    \]

\paragraph{System Cost Reduction (year)}
    Percentage variation of total annual cost, i.e., the cost sustained by consumers, is defined as:
    \[
    \Delta C_{\text{year}} = 
    \frac{C^{\text{SPaC}}_{\text{tot}} - C^{\text{base}}_{\text{tot}}}
         {C^{\text{base}}_{\text{tot}}} \times 100\%, \quad \text{with } \text{base} \in \{\text{PaB}, \text{PaC}\}.
    \]

\paragraph{Operator Profit Reduction (year)}
    Similarly, variation of total operator profits is:
    \[
    \Delta \Pi_{\text{year}} = 
    \frac{\Pi^{\text{SPaC}}_{\text{tot}} - \Pi^{\text{base}}_{\text{tot}}}
         {\Pi^{\text{base}}_{\text{tot}}} \times 100\%, \quad \text{with } \text{base} \in \{\text{PaB}, \text{PaC}\}.
    \]

\paragraph{Average Operator Markup}
    The markup represents an estimate of the aggregated system average overcharge: it indicates, in percentage terms, by how much revenues obtained by operators exceed respective production costs. Being a calculation based on total sum of profits and costs of all actors, the value reflects average profitability of the entire generator fleet rather than that of a single plant:
    \[
        \text{Markup} = 
        \frac{\Pi_{\text{tot}}}{C_{\text{tot}} - \Pi_{\text{tot}}} \times 100\%, 
    \]
    where $C_{\text{tot}}$ represents total system cost (inclusive of profits), therefore the denominator ($C_{\text{tot}} - \Pi_{\text{tot}}$) corresponds to the sum of production costs only.

\begin{table}[htb]
    \centering
    \begin{tabular}{lccc}
    \hline
    \textbf{} &  & \textbf{SPaC vs PaB} & \textbf{SPaC vs PaC} \\
    \hline
    Average PUN reduction (min/max) &   & $-17.7\%$ ($-0.7\%$ / $-46.4\%$)  & $-25.7\%$ ($-14.5\%$ / $-43.1\%$) \\
    System cost reduction & $\Delta C_{\text{year}}$         & $-20.4\%$                         & $-24.3\%$                         \\
    Operators profit reduction &  $\Delta \Pi_{\text{year}}$  & $-42.8\%$                         & $-50.2\%$                         \\
    Average markup &                & $51.0\%$ (vs $88.5\%$)            & $51.0\%$ (vs $105.3\%$)           \\
    \hline
    \end{tabular}
    \caption{Comparison SPaC vs PaB and PaC cases – Annual simulation with learned RL policy}
    \label{tab:confronto-spac-2024}
\end{table}

\noindent
Results in Table~\ref{tab:confronto-spac-2024} show that the SPaC mechanisms results a significant reduction of average PUN and system costs in a year, accompanied by a contraction of operator profits and a significantly lower markup compared to the alternative market configurations.

\section{Conclusions}
\label{sec:conclusioni}

This report has analyzed the Segmented Pay-as-Clear (SPaC) mechanism as an alternative to traditional regimes, the currently in force Pay-as-Clear (PaC) classical Pay-as-Bid (PaB), for the Italian electricity market, through simulations based on Reinforcement Learning to evaluate strategic behavior of operators.

\smallskip
\noindent
The simulation results confirm that the PaB mechanism, while presenting maximum efficiency under ideal, and unrealistic, perfect competition conditions, becomes the most expensive when operators adopt strategies optimized via Q-Learning, highlighting an intrinsic vulnerability to strategic abuse. PaC maintains high costs in both configurations, confirming systematic generation of inframarginal rents independently of operator behavior. SPaC emerges as a balanced intermediate solution, with contained cost increases and significantly limited rent extraction compared to other mechanisms.

\smallskip
\noindent
In the realistic scenario based on portfolios of ten operators, built through clustering of public offers of 2024, SPaC confirms its effectiveness in an oligopolistic context, limiting rent extraction and reducing PUN volatility. A particularly relevant result is SPaC's ability to maintain participation incentives for all operators while limiting overall rents, with a more equitable distribution suggesting greater political sustainability of the mechanism.

\subsection{Limitations and Future Developments}

Main limitations of this study concern:
\begin{itemize}
 \item simplification introduced by K-means clustering in portfolio reconstruction, which does not completely capture operational technical constraints;
 \item absence of a public database associating production units to respective technologies, making necessary to infer such classification indirectly;
 \item choice of a discrete action space for Q-Learning, representing a compromise between computational tractability and representativeness;
 \item assumption of rigid demand, which may overestimate impact of bidding strategies on equilibrium price.
\end{itemize}
Future developments of this research could include: implementation of more sophisticated Reinforcement Learning algorithms (Deep Q-Networks, Policy Gradient Methods) to explore continuous action spaces; integration of elastic demand functions and inter-temporal constraints; consideration of network constraints; systematic sensitivity analysis with respect to key parameters and market structure, as well as use of a database through which it is possible to precisely associate to each production unit its technology, as implicitly foreseen by a real implementation of SPaC.

\subsection{Policy Implications}

Despite the highlighted limitations, the results provide robust evidence on the potential effectiveness of the SPaC mechanism as a tool for electricity market reform, confirming the simplified simulations contained in the original \cite{frangioniBilevelProgrammingApproach2024}. Technological segmentation, together with the mechanism of variable demand split, provides a design principle that allows reconciling reduction of inframarginal while safeguarding sufficient incentives for investments in new capacity. SPaC effectiveness results particularly pronounced in asymmetric oligopolistic contexts, where segmentation introduces an endogenous control element that limits unilateral exercise of market power.

\smallskip
\noindent
The developed simulation framework constitutes a tool to support policy makers, regulators and eventually Nominated Energy Market Operators (NEMO) for ex-ante evaluation/proposal of market design reforms, allowing what-if analyses on different regimes and allowing quantification of impacts on PUN, demand cost, aggregated profits and their distribution. The tool can be extended to test \textit{complementarity} with other policy mechanisms (e.g., two-way contracts for difference, which by their nature have a long-term effect) and to explore trade-offs between efficiency, cost for consumers and investment signals.

\smallskip
\noindent
In conclusion, the SPaC mechanism represents a promising proposal for reform of the Italian (and European) electricity market, capable of reducing allocative and distributive inefficiencies of the currently used PaC mechanism while maintaining participation and investment incentives for all market operators.


\clearpage
\bibliographystyle{ieeetr}
\bibliography{bibliography}

@misc{acerACERGuidanceREMIT2024,
  title = {{{ACER}} Guidance on {{REMIT}} Application},
  author = {{Acer}},
  date = {2024-12-18},
  url = {https://www.acer.europa.eu/remit/about-remit/remit-guidance},
  urldate = {2025-07-07},
  langid = {english}
}

@article{aliabadiAgentbasedSimulationPower2017,
  title = {An Agent-Based Simulation of Power Generation Company Behavior in Electricity Markets under Different Market-Clearing Mechanisms},
  author = {Aliabadi, Danial Esmaeili and Kaya, Murat and Şahin, Güvenç},
  date = {2017-01},
  journaltitle = {Energy Policy},
  shortjournal = {Energy Policy},
  volume = {100},
  pages = {191--205},
  issn = {0301-4215},
  doi = {10.1016/j.enpol.2016.09.063},
  url = {https://www.sciencedirect.com/science/article/pii/S0301421516305419},
  urldate = {2025-07-11},
  langid = {english}
}

@online{areraDelibera30220252025,
  title = {Delibera 302/2025 - Rapporto sugli esiti del mercato elettrico del giorno prima nel biennio 2023-2024},
  author = {{Arera}},
  date = {2025-07-01},
  url = {https://www.arera.it/atti-e-provvedimenti/dettaglio/25/302-25},
  urldate = {2025-10-17},
  langid = {italian}
}

@article{bezansonJuliaFreshApproach2017,
  title = {Julia: {{A Fresh Approach}} to {{Numerical Computing}}},
  shorttitle = {Julia},
  author = {Bezanson, Jeff and Edelman, Alan and Karpinski, Stefan and Shah, Viral B.},
  date = {2017-01},
  journaltitle = {SIAM Review},
  shortjournal = {SIAM Rev.},
  volume = {59},
  number = {1},
  pages = {65--98},
  publisher = {{Society for Industrial and Applied Mathematics}},
  issn = {0036-1445},
  doi = {10.1137/141000671},
  url = {https://epubs.siam.org/doi/10.1137/141000671},
  urldate = {2025-08-26}
}

@article{boscoStrategicBiddingVertically2012,
  title = {Strategic Bidding in Vertically Integrated Power Markets with an Application to the {{Italian}} Electricity Auctions},
  author = {Bosco, Bruno and Parisio, Lucia and Pelagatti, Matteo},
  date = {2012-11},
  journaltitle = {Energy Economics},
  shortjournal = {Energy Economics},
  volume = {34},
  number = {6},
  pages = {2046--2057},
  issn = {0140-9883},
  doi = {10.1016/j.eneco.2011.11.005},
  url = {https://www.sciencedirect.com/science/article/pii/S0140988311002799},
  urldate = {2025-08-11},
  langid = {english}
}

@inproceedings{durvasuluClassificationGeneratorsParticipating2017,
  title = {Classification of Generators Participating in the Bulk-Power Market},
  booktitle = {2017 {{IEEE International Conference}} on {{Industrial Technology}} ({{ICIT}})},
  author = {Durvasulu, Venkat and Hansen, Timothy M. and Tonkoski, Reinaldo},
  date = {2017-03},
  pages = {575--579},
  doi = {10.1109/ICIT.2017.7915422},
  url = {https://ieeexplore.ieee.org/document/7915422},
  urldate = {2025-08-28},
  eventtitle = {2017 {{IEEE International Conference}} on {{Industrial Technology}} ({{ICIT}})}
}

@article{ferrariCompetitionElectricityMarkets2005,
  title = {Competition in Electricity Markets: International Experience and the Case of {{Italy}}},
  shorttitle = {Competition in Electricity Markets},
  author = {Ferrari, Alessandra and Giulietti, Monica},
  date = {2005-09},
  journaltitle = {Utilities Policy},
  shortjournal = {Utilities Policy},
  series = {Special {{Issue}} on {{Utility Service Quality}}},
  volume = {13},
  number = {3},
  pages = {247--255},
  issn = {0957-1787},
  doi = {10.1016/j.jup.2004.07.003},
  url = {https://www.sciencedirect.com/science/article/pii/S0957178704000724},
  urldate = {2025-07-05},
  langid = {english}
}

@article{frangioniBilevelProgrammingApproach2024,
  title = {A Bilevel Programming Approach to Price Decoupling in {{Pay-as-Clear}} Markets, with Application to Day-Ahead Electricity Markets},
  author = {Frangioni, Antonio and Lacalandra, Fabrizio},
  date = {2024-11},
  journaltitle = {European Journal of Operational Research},
  shortjournal = {European Journal of Operational Research},
  volume = {319},
  number = {1},
  pages = {316--331},
  issn = {03772217},
  doi = {10.1016/j.ejor.2024.06.018},
  url = {https://linkinghub.elsevier.com/retrieve/pii/S0377221724004661},
  urldate = {2025-07-01},
  langid = {english}
}

@online{garciaBilevelJuMPjlModelingSolving2022,
  title = {{{BilevelJuMP}}.Jl: Modeling and Solving Bilevel Optimization in Julia},
  shorttitle = {{{BilevelJuMP}}.Jl},
  author = {Garcia, Joaquim Dias and Bodin, Guilherme and Street, Alexandre},
  date = {2022-06-07},
  eprint = {2205.02307},
  eprinttype = {arXiv},
  eprintclass = {math},
  doi = {10.48550/arXiv.2205.02307},
  url = {http://arxiv.org/abs/2205.02307},
  urldate = {2025-08-12},
  langid = {english},
  pubstate = {prepublished}
}

@online{GestoreMercatiEnergetici,
  title = {Gestore Dei {{Mercati Energetici S}}.p.{{A}}.},
  url = {https://www.mercatoelettrico.org/it-it/},
  urldate = {2025-09-23}
}

@article{harderASSUMEAgentbasedSimulation2025,
  title = {{{ASSUME}}: {{An}} Agent-Based Simulation Framework for Exploring Electricity Market Dynamics with Reinforcement Learning},
  shorttitle = {Assume},
  author = {Harder, Nick and Miskiw, Kim K. and Khanra, Manish and Maurer, Florian and Patil, Parag and Qussous, Ramiz and Weinhardt, Christof and Klobasa, Marian and Ragwitz, Mario and Weidlich, Anke},
  date = {2025-05},
  journaltitle = {SoftwareX},
  shortjournal = {SoftwareX},
  volume = {30},
  pages = {102176},
  issn = {23527110},
  doi = {10.1016/j.softx.2025.102176},
  url = {https://linkinghub.elsevier.com/retrieve/pii/S2352711025001438},
  urldate = {2025-08-04},
  langid = {english}
}

@book{kirschenFundamentalsPowerSystem2004,
  title = {Fundamentals of {{Power System Economics}}},
  author = {Kirschen, Daniel and Strbac, Goran},
  date = {2004-03-26},
  edition = {1},
  publisher = {Wiley},
  location = {Chichester},
  doi = {10.1002/0470020598},
  url = {https://onlinelibrary.wiley.com/doi/book/10.1002/0470020598},
  urldate = {2025-08-21},
  isbn = {978-0-470-84572-1 978-0-470-02059-3},
  langid = {english},
  pagetotal = {284}
}

@online{KmeansClusteringjl2025,
  title = {K-Means · {{Clustering}}.Jl},
  date = {2025-01-06},
  url = {https://juliastats.org/Clustering.jl/dev/kmeans.html},
  urldate = {2025-10-02},
  langid = {english}
}

@inproceedings{linAutomatedClassificationPower2020,
  title = {Automated {{Classification}} of {{Power Plants}} by {{Generation Type}}},
  booktitle = {Proceedings of the {{Eleventh ACM International Conference}} on {{Future Energy Systems}}},
  author = {Lin, Liuzixuan and Chien, Andrew A.},
  date = {2020-06-18},
  series = {E-{{Energy}} '20},
  pages = {86--96},
  publisher = {Association for Computing Machinery},
  location = {New York, NY, USA},
  doi = {10.1145/3396851.3397708},
  url = {https://doi.org/10.1145/3396851.3397708},
  urldate = {2025-08-28},
  isbn = {978-1-4503-8009-6}
}

@book{schweppeSpotPricingElectricity1988,
  title = {Spot Pricing of Electricity},
  author = {Schweppe, Fred C. and Caramanis, Michael C. and Tabors, Richard D. and Bohn, Roger E.},
  date = {1988},
  publisher = {Springer US},
  location = {Boston, MA},
  doi = {10.1007/978-1-4613-1683-1},
  url = {http://link.springer.com/10.1007/978-1-4613-1683-1},
  urldate = {2025-07-08},
  isbn = {978-1-4612-8950-0 978-1-4613-1683-1},
  langid = {english}
}

@online{StudyLevelizedCost,
  title = {Study: {{Levelized Cost}} of {{Electricity}} - {{Renewable Energy Technologies}} - {{Fraunhofer ISE}}},
  shorttitle = {Study},
  url = {https://www.ise.fraunhofer.de/en/publications/studies/cost-of-electricity.html},
  urldate = {2025-10-01},
  langid = {english},
  organization = {Fraunhofer Institute for Solar Energy Systems ISE}
}

@book{suttonReinforcementLearningIntroduction2018,
  title = {Reinforcement {{Learning}}: {{An Introduction}}},
  shorttitle = {Reinforcement {{Learning}}},
  author = {Sutton, Richard S. and Barto, Andrew G.},
  editor = {Bach, Francis},
  editortype = {redactor},
  date = {2018-11-13},
  series = {Adaptive {{Computation}} and {{Machine Learning}} Series},
  edition = {2},
  publisher = {MIT Press},
  location = {Cambridge, MA, USA},
  url = {https://mitpress.mit.edu/9780262039246/reinforcement-learning/},
  urldate = {2025-08-21},
  isbn = {978-0-262-03924-6},
  langid = {english},
  pagetotal = {552}
}

@article{watkinsQlearning1992,
  title = {Q-Learning},
  author = {Watkins, Christopher J. C. H. and Dayan, Peter},
  date = {1992-05-01},
  journaltitle = {Machine Learning},
  shortjournal = {Mach Learn},
  volume = {8},
  number = {3},
  pages = {279--292},
  issn = {1573-0565},
  doi = {10.1007/BF00992698},
  url = {https://doi.org/10.1007/BF00992698},
  urldate = {2025-08-21},
  langid = {english}
}

@misc{report_code,
  author = {Altamura, Andrea and Lacalandra, Fabrizio},
  title = {{report-electricity-market-comparison-rl: Codice utilizzato per le simulazioni oggetto del report}},
  year = {2025},
  publisher = {GitHub},
  howpublished = {\url{https://github.com/andyalt23/report-electricity-market-comparison-rl}}
}

\end{document}